\documentclass[12pt]{amsart}
 \usepackage{mathtools}

\usepackage{amsmath,bbm, color}
\usepackage{latexsym,mathrsfs,bm}
\usepackage{url}
\usepackage{paralist}
\usepackage{amsthm,amsfonts,amssymb,amsmath}
\usepackage{comment}

\usepackage[latin1]{inputenc}

\DeclareMathOperator{\sign}{sign}

\DeclareMathOperator{\GL}{GL}

\DeclareMathOperator{\SL}{SL}

\DeclareMathOperator{\sgn}{sgn}

\newcommand\R{\mathbb{R}}
\newcommand\Z{\mathbb{Z}}
\newcommand\C{\mathbb{C}}

\newtheorem{thm}{Theorem}[section]

\newtheorem{prop}[thm]{Proposition}

\newtheorem{conj}[thm]{Conjecture}
\theoremstyle{remark}

\newtheoremstyle{named}{}{}{\itshape}{}{\bfseries}{.}{.5em}{\thmnote{#3}}
\theoremstyle{named}
\newtheorem*{namedtheorem}{Theorem}

\usepackage{color}

\let\originalleft\left
\let\originalright\right
\renewcommand{\left}{\mathopen{}\mathclose\bgroup\originalleft}
\renewcommand{\right}{\aftergroup\egroup\originalright}

\allowdisplaybreaks

\numberwithin{equation}{section}

\begin{document}
\title[Whittaker Fourier series in string theory]{Whittaker Fourier type solutions to differential equations arising from string theory}
\author{Ksenia Fedosova \& Kim Klinger-Logan}
\begin{abstract}
In this article, we find the full Fourier expansion for the generalized non-holomorphic Eisenstein series for certain values of parameters. We give a connection of the boundary condition on such Fourier series with convolution formulas on the divisor functions. Additionally, we discuss a possible relation with the differential Galois theory.
\end{abstract}
\maketitle

\section{Introduction}
The goal of this paper is to examine the Fourier expansion of the solutions to inhomogeneous eigenvalue equations involving of a product of two non-holomorphic 
Eisenstein series. Explicitly, for certain $\alpha,\beta\in \mathbb{Z}_{>0}+1/2$, we find a Whittaker Fourier expansion for a solution $f(z)$ to equations of the form \begin{equation}
\label{eq:DE}
	(\Delta - \lambda) f(z) = E_\alpha(z) E_\beta(z),  \quad z = x + i y \in \mathfrak{H},
\end{equation}
where the Eisenstein series, $E_s(z)$, is defined as 
\begin{equation}\label{def:Eisenstein_series}
E_s(z) := \sum_{\gamma\in P\cap \Gamma\backslash \Gamma} \text{Im}(\gamma z)^s 
\end{equation}
for $\Gamma =\SL_2(\Z)$ and $P$ the subgroup of upper triangular matrices. We recall the non-holomorphic Eisenstein series as  $E_s(z)$  converge absolutely for $\text{Re}(s)>1$ and are eigenfunctions for the Laplace operator
$
-\Delta :=- y^2(\partial^2_x + \partial^2_y)
$ with eigenvalue $s(1-s)$.

Solutions to equations of the form \eqref{eq:DE} have been studied in \cite{BKP2020, CGPWW2021, DKS2021_1, DKS2021_2, GMV2015, KKL, SDK, Pioline}. If they satisfy the $\SL_2(\mathbb{Z})$-automorphicity condition, these solutions are sometimes referred to as \textit{generalised non-holomorphic Eisenstein series}~\cite{CGPWW2021}. Such functions arise in the maximally supersymmetric $\mathcal{N}= 4$ super-Yang-Mills (SYM) theory when studying the duality properties of certain correlation functions in the $1/N$ expansion. For $\alpha,\beta\in \mathbb{Z}_{>0}+1/2$ with $\alpha+\beta = q +2,q, \dots$, the generalized non-holomorphic Eisenstein series appear in the the order $\frac{1}{N^{q}}$ contributions with $q\in\mathbb{Z}_{>0}$ to the SYM free energy $F=-\log Z$ 
\cite[p.\,6]{CGPWW2021}. At low orders there is an explicit connection between the correlator functions of the $SU(N)$ $\mathcal{N}=4$ super Yang-Mills theory in the $1/N$ expansion
 \cite{CGPWW2021} and the 10-dimensional type IIB superstring scattering amplitude of gravitons. 
The $D^6R^4$ interactions in the low energy expansion of the 4-loop supergraviton is given by the the solution to (\ref{eq:DE}) where $\lambda=12$ and $\alpha=\beta=3/2$ \cite{GMRV2010, GMV2015} and an explicit Fourier solution has been given this case in~\cite{GMV2015}.  More generally, solutions to (\ref{eq:DE}) for $\alpha, \beta\in \mathbb{Z}_{>0}$ are examples of modular graph functions and solutions have been found in \cite{DK2019, DKS2021_1, DKS2021_2}.

We use the following method to investigate solutions of \eqref{eq:DE}: for $z=x+iy$ we start by assuming  $f(z)$ is periodic in~$x$ and expand it in corresponding Fourier series. From \eqref{eq:DE}, we deduce an ordinary differential equation on every Fourier mode of $f$. Each of these differential equations is an inhomogeneous differential equation of the second order; the homogeneous part coincides with the modified Bessel equation, while the inhomogeneous part comprises an infinite sum involving polynomials and $K$-Bessel functions. The indices of the $K$-Bessel function in the inhomogeneous part are independent on the parameters in the homogeneous part. Assuming that the solution has this form, we introduce a system of linear equations on parameters upon which this special form depends. For certain physically relevant $(\lambda, \alpha, \beta)$  we solve this system of linear equations using a symbolic algebra system to obtain the searched parameters. In addition to finding all Fourier modes for such solutions, we are able to determine conditions on $(\lambda, \alpha, \beta)$ when solutions of this form do not exist. 
Experimentally, we are able to find explicit Fourier solutions in many instances outside those contained here; however, we have chosen to only include these for brevity. 

The method we present in this paper is motivated by the exact expression of the solution in~\cite{GMV2015}. In the former article, Green, Miller, and Vanhove found the explicit expressions for the Fourier modes of the function, satisfying \eqref{eq:DE} with $\lambda = 12$ and $\alpha=\beta=3/2$.
The Fourier modes of the solution in \cite{GMV2015} are exactly of the form \eqref{eq_sums_of_bessels_general_form}. We also note that the full spectral solution for $f$ in terms of $L^2(\Gamma\backslash\mathfrak{H})$-eigenfunctions is given in \cite{KKL}. The Fourier expansion of the solution $f$ to (\ref{eq:DE}) for $\alpha=\beta=3/2$ and $\lambda =12$ was also explicitly computed \cite{SDK} using the Poincar\'{e} series solution found in \cite[Appendix A]{GMV2015}. The method used in \cite{SDK} is different from that used in \cite{GMV2015} and outlined below. Importantly, in \cite{SDK}, the authors are not able to extend their method outside of the case where $\alpha=\beta$ in (\ref{eq:DE}); however, the method outlined in this paper does not require such a dependence.

In \cite{CGPWW2021}, Chester, Green, Pufu, Wang, and Wen generalized Eisenstein series were studied for certain values of $\alpha, \beta, \lambda$. Although each full Fourier expansion was not obtained in \cite{CGPWW2021}, the authors provided many important properties to the solution. We would like to note that they have expressed the solution to the zeroth Fourier term not in terms of $K_0$ and $K_1$ as we did, but rather in terms of modified Bessel functions of integer index. These representations are related to the ones found below via a recurrence relation of $K$-Bessel functions.

Inhomogeneous differential equations of Bessel type with inhomogeneous parts involving Bessel functions appear not only in string theory, but also in the theory of vector-valued automorphic functions. More precisely, in \cite{FPR}, Fedosova, Pohl, and Rowlett considered functions $\varphi: \mathfrak{H}\to V$ for some complex finite-dimensional vector space~$V$ that are Laplace eigenfunctions with eigenvalue $s(1-s)$ for $s \in \mathbb{C}$, thus 
\begin{equation}\label{eq:Delta_eigen_intro}
(\Delta - s(1-s) )\varphi= 0\,.
\end{equation}
Additionally, they required that $\varphi$ satisfies the \emph{twist-periodicity} condition 
\begin{equation}\label{eq:twist_intro}
\varphi(z+1) = A \varphi(z)
\end{equation}
for all~$z\in\mathfrak{H}$ for some $A \in \GL(V)$. The Fourier expansion of $\varphi$ is well-known for~$A$ being a unitary matrix: in this case, after considering the Fourier expansion of~$\varphi$, we arrive to a modified Bessel equation, depending on $A$. For diagonalizable~$A$, the modified Bessel equation is a homogeneous differential equation. Interestingly enough, if we allow a non-diagonalizable matrix, then some entries of the Fourier coefficient of~$\varphi$ satisfy the differential equation
\[
(y^2 \partial_y^2-\lambda - 4 \pi^2 n^2 y^2) f(y) = g(y), \quad n \in \Z,
\]
where $g$ is a combination of the modified Bessel function of the second kind and a certain polynomial in $y$.

\subsection{Discussion of main results}\label{Section:main_results_discussion}

We denote by $\mathscr{S}$ the set containing all  $(\lambda, \alpha, \beta)$ such that either 
\begin{enumerate}[(i)]
    \item $\lambda=12,30, 56$ and $\alpha=\beta=\tfrac{3}{2}$, or
    \item $\lambda=20$ and $\alpha=\tfrac{3}{2},\beta=\tfrac{5}{2}$, or
    \item $\lambda=30$ and $\alpha=\beta=\tfrac{5}{2}$, or
    \item $\lambda=30$ and $\alpha=\tfrac{3}{2},\beta=\tfrac{7}{2}$.
\end{enumerate}

Of these case, in Appendix C of \cite{CGPWW2021} Chester, Green, Pufu, Wang and Wen examine the zero mode of solutions (i) for $\lambda = 12, 30, 56, 90$ and~(ii) for $\lambda = 20, 42$ and (iii) and~(iv) for $\lambda = 30, 56, 90$. However, for the nonzero modes the full Fourier coefficients were not explicitly given. The method outlined in this paper gives all Fourier modes in these cases as well; however, we leave them out of this paper for the interest of space.

We are also able to find solutions of the form in Theorem A for $\alpha=\beta =3/2$ and $\lambda = 2$ (see Section \ref{sec:3/2,2}); $\alpha=3/2,\beta=5/2$ and $\lambda =6$ (see Section \ref{sec:3/2,5/2,6}); $\alpha=\beta =5/2$ and $\lambda = 2$  and $12$ (see Sections \ref{sec:5/2,12} and \ref{sec:5/2,12}); and $\alpha=3/2,\beta=7/2$ and $\lambda =12$ (see Section \ref{sec:3/2,7/2,12}). However, while the principle part of these solutions does converge,  the homogeneous part does not converge for the necessary choice of $\alpha_{n_1,n_2}$.

\begin{namedtheorem}[Theorem A]
Let $(\lambda=r(r+1), \alpha, \beta) \in \mathscr{S}$ and let $f: \mathfrak{H} \to \mathbb{C}$ be a $1$-periodic function in the $x$-variable that satisfies 
	\begin{equation*}
	(\Delta - \lambda) f(z) = E_\alpha(z) E_\beta(z), \quad z = x + i y \in \mathfrak{H}
	\end{equation*}
	for   $E_\bullet(z)$, $\bullet \in \{ \alpha, \beta \}$ as in \eqref{def:Eisenstein_series}.
	Then $f$ can be expressed in the form \[f(z) = \sum_{n \in \Z}  \hat{f}_n(y) e^{2 \pi i n x} = \sum_{n \in \Z} \sum_{n_1+n_2=n} \hat{f}_{n_1,n_2}(y) e^{2 \pi i n x}.\]
	For $n \neq 0$, \begin{align}\label{eq_sums_of_bessels_general_form}
	\hat{f}_{n_1,n_2}(y) &= \alpha_{n_1, n_2} \sqrt{y} K_{r+1/2} (2 \pi |n_1+n_2|y)  + \beta_{n_1, n_2} \sqrt{y} I_{r+1/2} (2 \pi |n_1+n_2|y) \\ &+  \sum_{i,j=0}^1 q^{i,j}(y) K_i(2 \pi |n_1| y) K_j(2 \pi |n_2| y), \nonumber 
	\end{align}
	where, for $\eta \in \mathbb{C}$,  $I_\eta$ and $K_\eta$ denote the modified Bessel function of the first and second kind of index $\eta$, respectively,  and $q^{i,j}=q^{i,j}_{n_1, n_2, \lambda, \alpha, \beta}$ are polynomials in $y$ and $\frac{1}{y}$. The case $n \neq 0$ can be obtained by a suitable limiting procedure. If we impose the requirement  
		\begin{equation}\label{condition:fourier_of_f_grows_not_too_fast}
	|\hat{f}_{n_1,n_2}(y)| = o(e^y), \quad y \to \infty,
	\end{equation}
	then $\beta_{n_1,n_2}$ vanishes; additionally, for each $(\alpha, \beta, \lambda)$, 
	there exists a unique choice of~$\alpha_{n_1,n_2}$ so that $\hat{f}_{n_1,n_2}(y)=o(y^{-r})$ as $y \to 0$.
\end{namedtheorem}

More precisely, we obtain the degrees of  polynomials $q^{i,j}$ in Figure \ref{figure}. We denote by 
$m_{i,j}$ and  $M_{i,j}$ the lowest and highest power of $y$ present in $q^{i,j}(y)$, respectively.

\begin{figure}[h]
{\large
\begin{center}
\begin{tabular}{|c|c|c|c|c|c|c|}
    \hline
     $(\alpha,\beta)$  & $m_{0,0}$ & $M_{0,0}$ &   $m_{0,1}$ & $M_{0,1}$ &  $m_{1,1}$ & $M_{1,1}$\\ \hline
  $(3/2,3/2)$  & $-r+2$ & $1$ &   $-r+1$& $0$ & $-r+2$ & $1$\\ \hline
    $(3/2,5/2)$  & $-r+2$ & $0$ &   $-r+1$& $1$ & $-r+2$ & $0$\\ \hline
        $(5/2,5/2)$  & $-r+2$ & $1$ &   $-r+1$& $0$ & $\min\{-r+1,-1\}$ & $1$\\ \hline
        $(3/2,7/2)$  & $-r+2$ & $1$ &   $-r+1$& $0$ & $-r+2$ & $1$\\ \hline
\end{tabular}
\end{center}\caption{For $\lambda=r(r+1)$, let $m_{i,j}$ and  $M_{i,j}$  be the lowest and highest power of $y$ present in $q^{i,j}(y)$, respectively. Note that 
$m_{0,1}=m_{1,0}$ and  $M_{0,1}=M_{1,0}$.}}\label{figure}\end{figure}

In \cite[Section C.1\,(a)]{CGPWW2021}, 
Chester, Green, Pufu, Wang, and Wen  conjectured, based on ideas from the AdS-CFT correspondence and Yang-Mills theory, that the total sum of the Fourier coefficients corresponding to the homogeneous solution vanishes. 
In \cite{SDK}, the authors provided an argument in support of this conjecture for every non-zero Fourier term for $\lambda=12$, $\alpha=\beta=3/2$ (the zeroth term can be dealt with with the help of Ramanujan summation formulas). 
We do not prove this conjecture in this article. However, the special choices of $\alpha_{n_1,n_2}$ made in each case in order to obtain a unique boundary condition imply that the conjecture would follow from a certain convolution series on divisor functions. 
  Following methods similar to those of \cite{GMV2015} and \cite{CGPWW2021}, we show at least for one choice of parameters and the zeroth coefficient (Section \ref{Section6.askojdoakjsdo}), the formal vanishing of the homogeneous part follows from a certain derivative of the Ramanujan identity.  
 If we want to deal with the non-zero Fourier coefficient, we would have to prove a more general version of the Ramanujan identity.

\subsection{Application to large $N$ expansion of integrated correlators}

In \cite{CGPWW2021}, the authors gave an evidence that the first few terms of the large-$N$ expansion of $\partial^4_m \log Z|_{m=0, b=1}$ can be expressed with the help of the generalized Eisenstein series. If we denote by $\mathcal{E}(r, \alpha, \beta; z)$ the modular functions that satisfy the inhomogeneous Laplace equation
\[
(\Delta - r(r+1)) \mathcal{E}(r, \alpha, \beta; z) = - 4 \zeta(2 \alpha) \zeta(2 \beta) E_{\alpha}(z) E_{\beta}(z),
\]
where $\zeta$ denotes the Riemann zeta function,
then $1/N^2$ contribution from \cite[(2.13)]{CGPWW2021} is conjectured to be equal to 
\begin{align*}
T_{-2}(y) =    C_1 + \frac{14175}{704\pi^4} \mathcal{E}(6,\tfrac{5}{2},\tfrac{3}{2})-\frac{1215}{88\pi^4} \mathcal{E}(4,\tfrac{5}{2},\tfrac{3}{2})
\end{align*}
for some constant $C_1$.  

With the help of the method described in the article it is possible to show that 
\[
C_1 + \frac{14175}{704\pi^4} \mathcal{E}(6,\tfrac{5}{2},\tfrac{3}{2})-\frac{1215}{88\pi^4} \mathcal{E}(4,\tfrac{5}{2},\tfrac{3}{2}) = C_1 +  \sum_{n \in \mathbb{Z}} T_{-2,n}(y)e^{2 \pi i n x},
\]
where for $n \neq 0$,
\begin{align*}
T_{-2,n}(y) &=  C_{2,n} \sqrt{y} K_{9/2}(2 \pi |n|y) +  C_{3,n} \sqrt{y} K_{13/2}(2 \pi |n|y) + \sum_{n_1+n_2 = n} T_{-2,n_1,n_2}(y),
\end{align*}
for some $C_{2,n}$ and $C_{3,n}$, where for $n_1+n_2 \neq 0$ and $n_1 n_2 \neq 0$,
\begin{align}
T_{-2,n_1,n_2}(y) = -& \left| n_1\right| {}^2 \left| n_2\right|  \sigma _{-4}\left(\left| n_1\right| \right) \sigma _{-2}\left(\left| n_2\right|\right)  \\
& \times \sum_{i, j=0}^1  w_{i,j}(y) K_i(2 \pi |n_1| y) K_j( 2 \pi |n_2| y) \nonumber
\end{align} 
for 
\begin{align*}
    w_{0,0}(y) &= \sgn(n_1) \frac{27 n_2}{ \pi ^6 \left(n_1+n_2\right){}^{12} y^4} \\ & \times (25200 n_1^6-201600 n_1^5 n_2+441000 n_1^4 n_2^2 -352800 n_1^3 n_2^3+88200 n_1^2 n_2^4 \\ &+(7630 \pi ^2 n_1^8-70840 \pi ^2 n_2 n_1^7+72460 \pi ^2 n_2^2 n_1^6+210200 \pi ^2 n_2^3 n_1^5-89240 \pi ^2 n_2^4 n_1^4\\ &-119080 \pi ^2 n_2^5 n_1^3+29620 \pi ^2 n_2^6 n_1^2+200 \pi ^2 n_2^7 n_1+10 \pi ^2 n_2^8) y^2 \\
    &+ (235 \pi ^4 n_1^{10}-4180 \pi ^4 n_2 n_1^9-1899 \pi ^4 n_2^2 n_1^8+26848 \pi ^4 n_2^3 n_1^7\\ & +35838 \pi ^4 n_2^4 n_1^6-3624 \pi ^4 n_2^5 n_1^5-23270 \pi ^4 n_2^6 n_1^4-6784 \pi ^4 n_2^7 n_1^3\\& +1383 \pi ^4 n_2^8 n_1^2+28 \pi ^4 n_2^9 n_1+\pi ^4 n_2^{10}) y^4),\\
    w_{0,1}(y) &= \sgn(n_1) \sgn(n_2)\frac{9}{ \pi ^{7} \left(n_1+n_2\right){}^{13} y^5} \\ & 
    \times (-75600 n_2^2 n_1^5+529200 n_2^3 n_1^4-793800 n_2^4 n_1^3+264600 n_2^5 n_1^2 \\
    &+ (210 \pi ^2 n_1^9+9450 \pi ^2 n_2 n_1^8-351720 \pi ^2 n_2^2 n_1^7+587160 \pi ^2 n_2^3 n_1^6\\
    &+1159200 \pi ^2 n_2^4 n_1^5-675360 \pi ^2 n_2^5 n_1^4-665280 \pi ^2 n_2^6 n_1^3+221760 \pi ^2 n_2^7 n_1^2\\&+630 \pi ^2 n_2^8 n_1+30 \pi ^2 n_2^9) y^2
    \\
    & + (45 \pi ^4 n_1^{11}+2385 \pi ^4 n_2 n_1^{10}-63642 \pi ^4 n_2^2 n_1^9-28218 \pi ^4 n_2^3 n_1^8+381438 \pi ^4 n_2^4 n_1^7 \\& +479622 \pi ^4 n_2^5 n_1^6-102672 \pi ^4 n_2^6 n_1^5-349392 \pi ^4 n_2^7 n_1^4-84243 \pi ^4 n_2^8 n_1^3\\&+26913 \pi ^4 n_2^9 n_1^2+402 \pi ^4 n_2^{10} n_1+18 \pi ^4 n_2^{11}) y^4 \\
    & + (20 \pi ^6 n_2 n_1^{12}-1140 \pi ^6 n_2^2 n_1^{11}-1940 \pi ^6 n_2^3 n_1^{10}+10004 \pi ^6 n_2^4 n_1^9+34632 \pi ^6 n_2^5 n_1^8\\
    &+37752 \pi ^6 n_2^6 n_1^7+7384 \pi ^6 n_2^7 n_1^6-15960 \pi ^6 n_2^8 n_1^5-11676 \pi ^6 n_2^9 n_1^4\\ &-1988 \pi ^6 n_2^{10} n_1^3+252 \pi ^6 n_2^{11} n_1^2+4 \pi ^6 n_2^{12} n_1) y^6
    ),
    \end{align*}
    
    \begin{align*}
    w_{1,0}(y) &= \frac{9 n_2}{\pi ^{7} n_1 \left(n_1+n_2\right){}^{13} y^5} (75600 n_1^7-529200 n_1^6 n_2+793800 n_1^5 n_2^2-264600 n_1^4 n_2^3\\
    &+(60690 \pi ^2 n_1^9-454230 \pi ^2 n_2 n_1^8+311040 \pi ^2 n_2^2 n_1^7\\&+1310400 \pi ^2 n_2^3 n_1^6-357840 \pi ^2 n_2^4 n_1^5-761040 \pi ^2 n_2^5 n_1^4\\&+88200 \pi ^2 n_2^6 n_1^3+7560 \pi ^2 n_2^7 n_1^2+630 \pi ^2 n_2^8 n_1+30 \pi ^2 n_2^9) y^2\\
    &+(5850 \pi ^4 n_1^{11}-62550 \pi ^4 n_2 n_1^{10}-53577 \pi ^4 n_2^2 n_1^9  +352947 \pi ^4 n_2^3 n_1^8 \\& +508248 \pi ^4 n_2^4 n_1^7-50088 \pi ^4 n_2^5 n_1^6-342822 \pi ^4 n_2^6 n_1^5-102702 \pi ^4 n_2^7 n_1^4\\&+21222 \pi ^4 n_2^8 n_1^3+1398 \pi ^4 n_2^9 n_1^2+87 \pi ^4 n_2^{10} n_1+3 \pi ^4 n_2^{11}) y^4 \\
    &+(20 \pi ^6 n_1^{13}-1140 \pi ^6 n_2 n_1^{12}-1940 \pi ^6 n_2^2 n_1^{11}+10004 \pi ^6 n_2^3 n_1^{10}\\&+34632 \pi ^6 n_2^4 n_1^9+37752 \pi ^6 n_2^5 n_1^8+7384 \pi ^6 n_2^6 n_1^7-15960 \pi ^6 n_2^7 n_1^6\\&-11676 \pi ^6 n_2^8 n_1^5-1988 \pi ^6 n_2^9 n_1^4+252 \pi ^6 n_2^{10} n_1^3+4 \pi ^6 n_2^{11} n_1^2) y^6) \\
    w_{1,1}(y) &= \sgn(n_2) \frac{9}{ \pi ^6 n_1 \left(n_1+n_2\right){}^{12} y^4} \\ & \times (210 n_1^8+9240 n_1^7 n_2-451680 n_1^6 n_2^2+1633560 n_1^5 n_2^3-1157280 n_1^4 n_2^4\\&+81240 n_1^3 n_2^5+6960 n_1^2 n_2^6+600 n_1 n_2^7+30 n_2^8\\
    &+ (150 \pi ^2 n_1^{10}+6960 \pi ^2 n_2 n_1^9-218082 \pi ^2 n_2^2 n_1^8+315504 \pi ^2 n_2^3 n_1^7+757824 \pi ^2 n_2^4 n_1^6\\&-285072 \pi ^2 n_2^5 n_1^5-446400 \pi ^2 n_2^6 n_1^4+60528 \pi ^2 n_2^7 n_1^3+4794 \pi ^2 n_2^8 n_1^2\\&+384 \pi ^2 n_2^9 n_1+18 \pi ^2 n_2^{10}) y^2 \\
    &+ (5 \pi ^4 n_1^{12}+420 \pi ^4 n_2 n_1^{11}-13025 \pi ^4 n_2^2 n_1^{10}-3196 \pi ^4 n_2^3 n_1^9+89202 \pi ^4 n_2^4 n_1^8\\&+116952 \pi ^4 n_2^5 n_1^7-9026 \pi ^4 n_2^6 n_1^6-73800 \pi ^4 n_2^7 n_1^5-23271 \pi ^4 n_2^8 n_1^4\\&+3652 \pi ^4 n_2^9 n_1^3+147 \pi ^4 n_2^{10} n_1^2+4 \pi ^4 n_2^{11} n_1) y^4 ).
\end{align*}
The cases $n_1+n_2=0$ or $n_1 n_2 = 0$ can be obtained by a certain limiting procedure.

The $1/N^3$ contribution from \cite[(2.13)]{CGPWW2021} is conjectured to be equal to 
\[
\alpha_r \mathcal{E}(3,\tfrac{3}{2},\tfrac{3}{2})+\sum_{r=5,7,9} \alpha_r \mathcal{E}(r,\tfrac{3}{2},\tfrac{3}{2}) + \beta_r \mathcal{E}(r,\tfrac{5}{2},\tfrac{5}{2}) + \gamma_r \mathcal{E}(r,\tfrac{7}{2},\tfrac{3}{2}) 
\]
for $\alpha_r, \beta_r, \gamma_r$ defined in \cite[(2.14)]{CGPWW2021}. We can write the expression above as 
\[
\sum_{n \in \mathbb{Z}}  T_{-3,n}(y)e^{2 \pi i n x},
\]
where for $n \neq 0$,
\begin{align*}
T_{-3,n}(y) &=  C_{4,n} \sqrt{y} K_{7/2}(2 \pi |n|y) +  C_{5,n} \sqrt{y} K_{11/2}(2 \pi |n|y) \\ & +  C_{6,n} \sqrt{y} K_{15/2}(2 \pi |n|y) + C_{7,n} \sqrt{y} K_{19/2}(2 \pi |n|y) \\ &+ \sum_{n_1+n_2 = n} T_{-3,n_1,n_2}(y),
\end{align*}
for some $C_{j,n}$, $j=4,5,6,7$, where for $n_1+n_2 \neq 0$ and $n_1 n_2 \neq 0$,
\begin{align}
T_{-3,n_1,n_2}(y) = \sum_{i, j=0}^1  v_{i,j}(y) K_i(2 \pi |n_1| y) K_j( 2 \pi |n_2| y), \nonumber
\end{align} 
where $v_{i,j}(y)$ is some rational function on $y$.

\subsection{Acknowledgements}
We would like to thank Michael Green, Mark van Hoeij, Remi Jaoui, Askold Khovanskii, Stephen D.\,Miller,	Alexey Ovchinnikov, Wolfgang Soergel, Congkao Wen, Katrin Wendland and Don Zagier for discussions. K.\,K-L. acknowledges support from NSF Grant number DMS-2001909.

\section{Method of solution}

In this section, we outline a method for finding  the Fourier expansions of solutions~$f(z)$ to equations of the form 
\begin{equation}\label{eq:main_equation}
(\Delta - \lambda) f(z) = c_{\alpha,\beta}\,  \zeta(2 \alpha) \zeta(2 \beta) E_\alpha(z) E_\beta(z) 
\end{equation}
where $\lambda=r(r+1)$ for $r \in \mathbb{N}$
and~$\alpha,\beta \in \Z _{>0}+1/2$, and $c_{\alpha,\beta}$ is some constant depending on $\alpha$ and~$\beta$. 
The constants $c_{\alpha,\beta} \in \mathbb{C}$  are chosen for the convenience purpose and to shorten the outcome. For the particular $\alpha,\beta$ for which we write down the exact solutions, we let
\[
c_{\alpha, \beta} := \begin{cases}
4,& \quad \alpha = \beta = 3/2, \\
6,& \quad \alpha = 3/2, \beta = 5/2,\\
9, & \quad \alpha = \beta = 5/2, \\
30, & \quad \alpha = 3/2, \beta = 7/2.
\end{cases}
\]

We start with recalling that for $\text{Re}(s) > 1$,  
\[
E_s(z) = \sum_{n \in \mathbb{Z}} a_{n,s}(y) e^{2 \pi i n x},
\]
where 
\begin{align*}\label{eq:zeroth_coefficient}
a_{0,s}(y)= y^s + \frac{\sqrt{\pi} \Gamma(s-\tfrac{1}{2}) \zeta(2s-1)}{\Gamma(s) \zeta(2s)} y^{1-s}
\end{align*}
and for $n \neq 0$,
\begin{equation}\label{eq:non-zeroth_coefficient}
a_{n,s}(y) = \frac{2 \pi^s}{\Gamma(s) \zeta(2s)} |n|^{s-\tfrac{1}{2}} \sigma_{1-2s}(|n|) \sqrt{y} K_{s-\tfrac{1}{2}}(2 \pi |n|y)
\end{equation}
where for $z \in \mathbb{C}$ and $n \in \mathbb{N}$,  $\sigma_z(n) := \sum_{d | n} d^z$  is the divisor function \cite[p. 278]{Zagier81}. We note that in the notations of \cite[(2.10)]{CGPWW2021},
\[
E_s(z) = \frac{1}{2 \zeta(2s)} E(s, z, \overline{z}),
\]
and thus 
\begin{equation}\label{our_solution_vs_Chester}
4  \zeta(2 \alpha) \zeta(2 \beta) E_\alpha(z) E_\beta(z) = - E(\alpha, z, \overline{z}) E(\beta, z, \overline{z}).
\end{equation}

  This expansion implies that for $\alpha,\beta>1$
\begin{equation}\label{eq:Fourier_of_product}
 c_{\alpha,\beta}\,  \zeta(2 \alpha) \zeta(2 \beta)E_\alpha (z) E_\beta(z) = \sum_{n \in \mathbb{Z}} S_n(y) e^{2 \pi i n x},
\end{equation}
where 
\begin{equation}\label{Sn:def}
S_n(y) = \sum_{\substack{n_1, n_2 \in \Z\\ n_1+n_2=n}}  s_{n_1, n_2}(y)
\end{equation}
for 
\[
s_{n_1, n_2}(y) = c_{\alpha,\beta}\,  \zeta(2 \alpha) \zeta(2 \beta) a_{n_1,\alpha}(y) a_{n_2,\beta}(y).
\]
Explicitly, each $s_{n_1, n_2}(y)$ can be written as follows: 
\begin{enumerate}
	\item For $n_1=n_2=0$:
	\begin{align*}
	s_{0,0}(y)
	&=  c_{\alpha,\beta}\,  \zeta(2 \alpha) \zeta(2 \beta)y^{\alpha+\beta} + c_{\alpha,\beta}\, \zeta(2 \beta) y^{1-\alpha+\beta}  \frac{\sqrt{\pi} \Gamma(\alpha-\tfrac{1}{2}) \zeta(2\alpha-1)}{\Gamma(\alpha) } \\
	&+  c_{\alpha,\beta}\,  \zeta(2 \alpha)y^{1+\alpha-\beta} \frac{\sqrt{\pi} \Gamma(\beta-\tfrac{1}{2}) \zeta(2\beta-1)}{\Gamma(\beta) } \\ & +  c_{\alpha,\beta}\,  y^{2-\alpha-\beta} \frac{\pi \Gamma(\alpha-\tfrac{1}{2}) \Gamma(\beta-\tfrac{1}{2}) \zeta(2\alpha-1)  \zeta(2\beta-1)}{\Gamma(\alpha) \Gamma(\beta)  }.   
	\end{align*}
	\item For  $n_1 = 0$, $n_2 \neq 0$:
	\begin{align*}
		s_{0,n}(y)
		&= \frac{2c_{\alpha,\beta} \pi^\beta \zeta(2 \alpha) }{\Gamma(\beta)\zeta(2 \beta)}  |n|^{\beta-\tfrac{1}{2}} \sigma_{1-2\beta}(|n|) y^{\alpha+1/2} K_{\beta-\tfrac{1}{2}}(2 \pi |n|y) \\
		&+ \frac{ 2c_{\alpha,\beta} \pi^{\beta+1/2} \Gamma(\alpha-\tfrac{1}{2}) \zeta(2\alpha-1)}{\Gamma(\alpha)  \Gamma(\beta) \zeta(2 \beta)} 
		 |n|^{\beta-\tfrac{1}{2}} \sigma_{1-2\beta}(|n|) y^{3/2-\alpha}  K_{\beta-\tfrac{1}{2}}(2 \pi |n|y).
	\end{align*}
	\item For  $n_1 \neq 0$, $n_2 = 0$:
	\begin{align*}
	s_{n,0}(y) &=  \frac{2 c_{\alpha,\beta} \pi^\alpha \zeta(2 \beta)}{\Gamma(\alpha) \zeta(2 \alpha)} |n|^{\alpha-\tfrac{1}{2}} \sigma_{1-2\alpha}(|n|) y^{\beta+1/2} K_{\alpha-\tfrac{1}{2}}(2 \pi |n|y) \\
	&+ \frac{ 2c_{\alpha,\beta} \pi^{\alpha+1/2} \Gamma(\beta-\tfrac{1}{2}) \zeta(2\beta-1)}{\Gamma(\beta)  \Gamma(\alpha) \zeta(2 \alpha)} 
	|n|^{\alpha-\tfrac{1}{2}} \sigma_{1-2\alpha}(|n|) y^{3/2-\beta}  K_{\alpha-\tfrac{1}{2}}(2 \pi |n|y).
	\end{align*}
	\item For $n_1 n_2 \neq 0$:
	\begin{align}\label{sn1n2:non-zero_case}
		s_{n_1,n_2}(y) &= \frac{4c_{\alpha,\beta} \pi^{\alpha+\beta}}{\Gamma(\alpha)  \Gamma(\beta) } |n_1|^{\alpha-\tfrac{1}{2}} |n_2|^{\beta-\tfrac{1}{2}} \sigma_{1-2\alpha}(|n_1|)    \sigma_{1-2\beta}(|n_2|)  \\ &  \ \ \ \ \ \ \ \times  y K_{\alpha-\tfrac{1}{2}}(2 \pi |n_1|y) K_{\beta-\tfrac{1}{2}}(2 \pi |n_2|y). \nonumber
	\end{align}
	\end{enumerate}

To solve \eqref{eq:main_equation}, note that the Fourier expansion of 
the right
is given by \eqref{eq:Fourier_of_product}.
 Although we do not assume the $SL_2(\mathbb{Z})$-invariance of $f$, we do require that, for $z=x+iy$, $f(z)$ is periodic in the $x$-direction with period~$1$. Given this assumption, the differential equation, \eqref{eq:main_equation}, can be equivalently stated as the simultaneous differential equations on the Fourier coefficients, $\hat{f}_n(y)$, of $f(z)$:
\begin{equation}\label{key}
(y^2 \partial_y^2-\lambda - 4 \pi^2 n^2 y^2) \hat{f}_n(y) = S_n(y), \quad n \in \Z.
\end{equation}
We express 
\[
\hat{f}_n(y) = \sum_{\stackrel{n_1,n_2 \in \mathbb{Z}}{n_1+n_2=n}} \hat{f}_{n_1,n_2}(y),
\]
for $\hat{f}_{n_1,n_2}(y)$ satisfying
\begin{equation}\label{eq:fn1n2}
(y^2 \partial_y^2-\lambda - 4 \pi^2 (n_1+n_2)^2 y^2) \hat{f}_{n_1,n_2}(y) = s_{n_1,n_2}(y).
\end{equation}
Each solution of \eqref{eq:fn1n2} can be written as a sum of a solution, $\hat{f}_{n_1,n_2}^H(y)$, of the homogeneous equation
\begin{equation}\label{eq:fn1n2hom}
(y^2 \partial_y^2-\lambda - 4 \pi^2 (n_1+n_2)^2 y^2) \hat{f}^H_{n_1,n_2}(y) = 0,
\end{equation}
and a particular solution, $\hat{f}_{n_1,n_2}^P(y)$, of \eqref{eq:fn1n2}. Thus,
\begin{equation}\label{homogeneous_equation_solutioniuz2iu3j4h}
\hat{f}_{n_1,n_2}(y) = \hat{f}_{n_1,n_2}^P(y) + \hat{f}_{n_1,n_2}^H(y).
\end{equation}

We note that for $\lambda=r(r+1)$, $r \in \mathbb{R}$ and $n_1+n_2 \neq 0$,
\[
\hat{f}^H_{n_1,n_2}(y) = 
\alpha_{n_1, n_2} \sqrt{y} K_{r+1/2} (2 \pi |n_1+n_2|y)  + \beta_{n_1, n_2} \sqrt{y} I_{r+1/2} (2 \pi |n_1+n_2|y)
\]
 for some $\alpha_{n_1, n_2}, \beta_{n_1, n_2} \in \C $. However, for $\beta_{n_1, n_2} \neq 0$, the function $\hat{f}^H_{n_1,n_2}(y)$ grows exponentially as $y \to \infty$, that contradicts \eqref{condition:fourier_of_f_grows_not_too_fast}. From this we deduce that $\beta_{n_1, n_2} = 0$
and thus 
\begin{equation}\label{eq:homogeneous_part123sdf}
\hat{f}^H_{n_1,n_2}(y) = \alpha_{n_1, n_2} \sqrt{y} K_{r+1/2} (2 \pi |n_1+n_2|y)
\end{equation}
for some $\alpha_{n_1, n_2} \in \C$. We note that for $\text{Re}(r)>-1/2$ and $y \to 0$, \cite[10.30.2]{NIST} implies
\begin{align}\label{asymptotics_of_homogeneous_solutionnasokdijfok}
\hat{f}^H_{n_1,n_2}(y) = \alpha_{n_1, n_2} 	y^{-r} \left(\tfrac{1}{2} |\pi n|^{-r-\tfrac{1}{2}} \Gamma \left(r+\tfrac{1}{2}\right)+O\left(y^2\right)\right).
\end{align} 
In the case where $n_1+n_2 = 0$, for $\lambda=r(r+1)$, $r \in \mathbb{R}$, $r\neq 1/2$,
\[
\hat{f}^H_{n_1,n_2}(y) = 
\alpha_{n_1,n_2}y^{-r} + \beta_{n_1,n_2} y^{r+1}
\]
for some $\alpha_{n_1,n_2}, \beta_{n_1,n_2} \in \C$. If we demand $\hat{f}^H_{n_1,n_2}(y) = o(y^{r+1})$ as $y\to \infty$, we would have to take $\beta_{n_1,n_2}=0$ and thus 
\begin{equation}\label{eq:homogeneous_part123sdf23sdf}
\hat{f}^H_{n_1,n_2}(y) = \alpha_{n_1, n_2} y^{-r}.
\end{equation}

It remains to find a particular solution, $\hat{f}^P_{n_1,n_2}(y)$. In what follows we assume that the solution is a linear combination of special functions multiplied by rational functions. We find the explicit constants which appear in front of these special functions by solving systems of certain linear equations. In Section \ref{sec:part}, we describe in more details which system of linear equations need to be solved depending on the values of $(n_1,n_2)$. More precisely, we will consider the following cases:
\begin{enumerate}
	\item In Section \ref{sec:systemasdfkjh}, we will consider the case  $n_1 n_2 \neq 0$.
	\item In Section \ref{sec:02938409uijeojhbfk}, we will consider the case when exactly one of $n_i$, $i=1,2$, is equal to zero.
	\item In Section \ref{sec:osaiujoiujwnkjh324o}, we will consider the case $n_1 = n_2 = 0$.
\end{enumerate}
Finally, in Section \ref{sec:convergence}, we will show that, at least the the physically relevance cases we have considered, the sums $\displaystyle\sum_{n_1+n_2=n}\hat f_{n_1,n_2}(y)$ converge for each $n$.

After we outline the strategy of finding solutions, we provide in further sections explicit examples of such for some physically relevant $\lambda, \alpha$ and $\beta$. We stress that we are able to find the solutions as functions of $n_1$ and $n_2$ without restricting ourselves to any particular values of $n_1$ and $n_2$. More precisely, we write down the explicit solutions for $f(z)$ in the following cases:
\begin{quote}
\begin{enumerate}[I.]
	\item $\lambda = 30, 56$ and $\alpha=\beta=3/2$ in Section \ref{sec:3/2} (we omit $\lambda=12$, because it has been treated in \cite{GMV2015}), 
	\item $\lambda = 20$ and $\alpha=3/2$, $\beta=5/2$ in Section \ref{sec:3/2,5/2}, 
	\item $\lambda =  30$ and $\alpha=\beta=5/2$  in Section \ref{sec:5/2}, 
	\item $\lambda = 30$ and $\alpha=3/2$, $\beta=7/2$ in Section~\ref{sec:3/2,7/2},
	\item  $\lambda = 2$ and $\alpha=\beta =3/2$; $\lambda =6$ and $\alpha=3/2,\beta=5/2$; $\lambda = 2, 12$ and $\alpha=\beta =5/2$; and $\lambda =12$ and $\alpha=3/2,\beta=7/2$ in Section \ref{sec:Appendix A}.
\end{enumerate}\end{quote}
We note that the given solutions partially cover  generalised Eisenstein series  that arise in coefficients $\mathcal{H}(q, \tau,\overline{\tau})$ of even terms in the $1/N$  up to order $1/N^3$ in \cite[(2.11)]{CGPWW2021}.
To be more precise, \cite[(2.11)]{CGPWW2021} expresses $\partial_m^4 \log Z|_{m = 0, b=1}$, that is a fourth derivative of the squashed sphere partition function of the $N= 2$ theory with respect to the squashing parameter $b=1$ and mass parameter $m=0$.

 We have not included the solutions of the differential equations that would cover the full expansion up to the order $1/N^3$ to keep the length of the article reasonable. 

\section{Particular solutions}\label{sec:part}
In this section, we explicitly describe the system of linear equations that finds a particular solution of \eqref{eq:fn1n2}, depending on the values of $(n_1, n_2)$.
\subsection{Solutions to \eqref{eq:fn1n2} for $n_1 n_2 \neq 0$}\label{sec:systemasdfkjh}
Substituting \eqref{sn1n2:non-zero_case} into \eqref{eq:fn1n2} and denoting 
\[
\hat{f}^P_{n_1,n_2}(y)  = \frac{4 \pi^{\alpha+\beta}}{\Gamma(\alpha) \zeta(2\alpha) \Gamma(\beta) \zeta(2\beta)} |n_1|^{\alpha-\tfrac{1}{2}} |n_2|^{\beta-\tfrac{1}{2}} \sigma_{1-2\alpha}(|n_1|)    \sigma_{1-2\beta}(|n_2|) g(y),
\]
we obtain a differential equation on $g$:
\begin{align}
(-4 \pi^2 y^2 (n_1+n_2)^2 + y^2 \partial_y^2 - \lambda) g(y) &= y K_{\alpha-1/2} (2 \pi | n_1 | y) K_{\beta-1/2} (2 \pi |n_2 | y),
\end{align}
or
\begin{align}\label{key_ksdfjnhjsojoilmmsl}
(-4 \pi^2 y^2 (|n_1|+ \sgn(n_1 n_2) |n_2|)^2 + y^2 \partial_y^2 - \lambda) g(y) = y K_{\alpha-1/2} (2 \pi  |n_1|  y) K_{\beta-1/2} (2 \pi |n_2|  y).
\end{align}
We introduce the notation 
\[
P_{\lambda} := -4 \pi^2 y^2 (|n_1| + \sgn(n_1 n_2) |n_2|)^2 + y^2 \partial_y^2 - \lambda.
\]
In this notation, \eqref{key_ksdfjnhjsojoilmmsl} reads
\begin{equation}\label{P_sldjhkf9ousdofisj}
\begin{aligned}
P_{\lambda} g(y) &= y K_{\alpha-1/2} (2 \pi  |n_1|  y) K_{\beta-1/2} (2 \pi |n_2|  y) .
\end{aligned}
\end{equation}
If $\alpha, \beta \in \Z + \tfrac{1}{2} $, then, using recursive formulas for the modified Bessel functions as in \cite[10.29(i)]{NIST} -- or, for the particular choices of $\alpha$ and $\beta$, as in Section \ref{sec:bessels_and_relations_between_them} --  we can rewrite the right hand side for \eqref{P_sldjhkf9ousdofisj} as 
\begin{equation}\label{wqLkjlkjlkjwewerwempjk}
\sum_{i,j=0}^1 h^{i,j}(y) K_i(2 \pi |n_1| y) K_j(2 \pi |n_2| y)
\end{equation}
where $h^{i,j}$ for each $i,j \in \{0,1\}$ is a polynomial in $y$ and $y^{-1}$.

We note that for any $a,b,c,d \in \mathbb{Z}$,
\begin{align}
P_{\lambda} (y^a K_0( 2 \pi |n_1| y) K_0(2 \pi |n_2| y)) &= (-\sgn(n_1 n_2)  8 \pi^2  |n_1| |n_2| y^{a+2}+a^2 y^a-\lambda  y^a-a y^a)\nonumber\\
 &  \ \ \ \ \ \  \ \ \ \ \ \ \ \ \ \ \ \  \ \ \ \ \ \ \times K_0(2 \pi |n_1| y) K_0(2 \pi |n_2| y) \nonumber\\
& + (2 \pi |n_2| y^{a+1}-4 \pi |n_2| a y^{a+1}) K_0(2 \pi |n_1| y) K_1(2 \pi |n_2| y)\nonumber \\
&+( 2 \pi |n_1| y^{a+1}-4 \pi |n_1| a y^{a+1}) K_1(2 \pi |n_1| y) K_0(2 \pi |n_2| y)\label{eq:La}\\
& + (8 \pi^2 |n_1| |n_2| y^{a+2})  K_1(2 \pi |n_1| y) K_1(2 \pi |n_2| y)\nonumber\\ 
P_{\lambda} (y^b K_0(2 \pi |n_1| y) K_1(2 \pi |n_2| y)) &= ( 2 \pi |n_2| y^{b+1}- 4 \pi  b |n_2| y^{b+1}) K_0(2 \pi |n_1| y) K_0(2 \pi |n_2| y) \nonumber\\
&+ (b^2 y^b -\sgn(n_1 n_2)  8 \pi^2 |n_1| |n_2|  y^{b+2}-\lambda  y^b-3 b y^b+2 y^b)\nonumber\\
 &  \ \ \ \ \ \  \ \ \ \ \ \ \ \ \ \ \ \  \ \ \ \ \ \ \times K_0(2 \pi |n_1| y) K_1(2 \pi |n_2| y)\nonumber\\
&+ (8 \pi^2 |n_1| |n_2| y^{b+2}) K_1(2 \pi |n_1| y) K_0(2 \pi |n_2| y)\label{eq:Lb}\\
&+ (6 \pi  |n_1| y^{b+1}- 6 \pi  b |n_1| y^{b+1})  K_1(2 \pi |n_1| y) K_1(2 \pi |n_2| y)\nonumber \\
P_{\lambda} (y^c K_1(2 \pi_1 |n_1| y) K_0(2 \pi |n_2| y)) &= ( 2 \pi |n_1| y^{c+1}- 4 \pi  |n_1| c y^{c+1})K_0(2 \pi |n_1| y) K_0(2 \pi |n_2| y)\nonumber \\
&+ (8 \pi^2 |n_1| |n_2| y^{c+2}) K_0(2 \pi |n_1| y) K_1(2 \pi |n_2| y) \nonumber\\
& - (\sgn(n_1 n_2) 2 8 \pi^2 |n_1| |n_2| y^{c+2}-c^2 y^c+\lambda  y^c+3 c y^c-2 y^c) \nonumber\\
 &  \ \ \ \ \ \  \ \ \ \ \ \ \ \ \ \ \ \  \ \ \ \ \ \ \times K_1(2 \pi |n_1| y) K_0(2 \pi |n_2| y)\label{eq:Lc} \\
&+ ( 6 \pi |n_2| y^{c+1}- 4 \pi |n_2| c y^{c+1})  K_1(2 \pi |n_1| y) K_1(2 \pi |n_2| y) \nonumber\\
P_{\lambda}( y^d  K_1(2 \pi |n_1| y) K_1(2 \pi |n_2| y)  )  & = (8 \pi^2 |n_1| |n_2| y^{d+2}) K_0(2 \pi |n_1| y) K_0(2 \pi |n_2| y)\nonumber\\
&+ (6 \pi |n_1| y^{d+1}-4 \pi |n_1| d y^{d+1})K_0(2 \pi |n_1| y) K_1(2 \pi |n_2| y) \nonumber\\
& + ( 6 \pi |n_2| y^{d+1}-4 \pi |n_2| d y^{d+1} ) K_1(2 \pi |n_1| y) K_0(2 \pi |n_2| y)\label{eq:Ld} \\
& - ( \sgn(n_1 n_2) 8 \pi^2  |n_1| |n_2| y^{d+2}-d^2 y^d+\lambda  y^d+5 d y^d-6 y^d  ) \nonumber\\
 &  \ \ \ \ \ \  \ \ \ \ \ \ \ \ \ \ \ \  \ \ \ \ \ \ \times K_1(2 \pi |n_1| y) K_1(2 \pi |n_2| y)\nonumber.
\end{align}
If we assume that a solution of \eqref{key_ksdfjnhjsojoilmmsl} can be expressed as a sum 
\begin{equation}\label{beautiful_equation_form}
    g(y) = \sum_{i,j=0}^1 q^{i,j}(y) K_i(2 \pi |n_1| y) K_j(2 \pi |n_2| y),
\end{equation}
where $q^{i,j}$ are some polynomials in $y$ and $y^{-1}$, then $q^{i,j}$ can be found by solving a system of linear equations on the coefficients of $q^{i,j}$.
More precisely, assume that 
\[
\max_{i,j} \deg q^{i,j} = M, \quad \min_{i,j} \deg q^{i,j} = m.
\] 
Then each of $q^{i,j}$ for $i,j \in \{0,1\}$ is parametrized by $(M-m+1)$ complex coefficients.

On the other hand, 
\[
P_{\lambda}(g) = \sum_{i,j=0}^1 \tilde{q}^{i,j}(y) K_i(2 \pi |n_1| y) K_j(2 \pi |n_2| y)
\]
for some polynomials $\tilde{q}^{i,j}$ such that 
\[
\max_{i,j} \deg \tilde{q}^{i,j} = M+2, \quad \min_{i,j} \deg \tilde{q}^{i,j} = m.
\]
Thus, \eqref{P_sldjhkf9ousdofisj} can be equivalently written as 
 $4(M+3-m)$ linear equations with ${4(M-m+1)}$ variables: the variables are exactly the coefficients of $q^{i,j}$, and the linear equations come from the requirement that $h^{i,j}=\tilde{q}^{i,j}$. 

We note that as one can see from \eqref{eq:La}-\eqref{eq:Lc}, the corresponding matrix of the system of linear equation is a band matrix, that simplifies the calculations.

The possibility of writing $g$ in such form depends on $\alpha$, $\beta$, $\lambda$, $M$ and $m$. Below, we write down some elementary limitations on the set of parameters that are needed in order for a solution of such form to exist. Further in the article, we  speculate on possible connection of restriction with the differential Galois theory. 

{\prop A solution to (\ref{eq:DE}) of the form  \eqref{beautiful_equation_form} with the condition 
	\[\min_{i,j\in\{0,1\}} \deg h^{i,j} > \min_{i,j\in\{0,1\}} \deg q^{i,j}\] 
	 does not exist unless $\lambda$ is of the form $r(r+1)$ for $r\in\mathbb{Z}_{>0}$.}
{\proof 
We note that from $h^{i,j}=\tilde{q}^{i,j}$, the equality 
\[
\min_{i,j\in\{0,1\}} \deg h^{i,j} = \min_{i,j\in\{0,1\}} \deg \tilde{q}^{i,j}
\]
must hold, that implies
\[
\min_{i,j\in\{0,1\}} \deg \tilde{q}^{i,j} < \min_{i,j\in\{0,1\}} \deg q^{i,j}.
\]
Together with \eqref{eq:La}-\eqref{eq:Ld}, the inequality above implies that there exist $a,b,c,d\in\Z$ such that at least one of the following equalities holds:
\begin{align*}
a^2-\lambda - a &= 0,\\
b^2-\lambda-3b+2&=0,\\
c^2-\lambda-3c+2&=0,\\
d^2-\lambda-5d+6&=0.
\end{align*}
That implies the statement of the proposition. \qed 
}

Additionally, we prove the following elementary proposition:
{\prop
A solution to (\ref{eq:DE}) of the form  (\ref{beautiful_equation_form}) does not exist unless $\alpha \in \tfrac{1}{2} + \mathbb{Z}$ and $\beta \in \tfrac{1}{2}+\mathbb{Z}$.
}
{\proof 
	We give the proof  for $\alpha$ by contradiction; the proof for $\beta$ is similar. 
Consider the right hand sides of \eqref{eq:La}-\eqref{eq:Ld}. Although we have assumed $n_1, n_2$ to be non-zero integers, the formulas above would hold if we let $n_1, n_2$ be non-zero real numbers. Keeping that in mind, we let $n_2 = n_2(y) = \tfrac{1}{y}$ depend on $y$. Having fixed the product of $n_2$ and $y$ and keeping in mind asymptotic expansions of the modified Bessel function of the second kind, we consider the corresponding asymptotic expansions of $P_{\lambda, \pm} g(y)$ for $y\to 0$ only to find integer powers of $y$ and $\log(y)$.

On the other hand, if $\alpha \neq 0$, the asymptotic expansion of $y K_{\alpha-1/2} (2 \pi |n_1| y) K_{\beta-1/2}(2 \pi)$ as $y \to 0$ contains only terms of the type $y^{\alpha + 1/2 + \star}$ for $\star \in \mathbb{Z}$. Thus, $\alpha \in \mathbb{Z}+1/2$.
\qed 
}

In what follows, we give explicit solutions to \eqref{key_ksdfjnhjsojoilmmsl}  for some physically relevant combinations of $\alpha,\beta$ and $\lambda$.

\subsection{Solutions to \eqref{eq:fn1n2} for $n_1 n_2 = 0$, but not both zero}\label{sec:02938409uijeojhbfk}
Without loss of generality we assume $n_1 = 0, n_2 \neq 0$. We note that if we find $g_1$ and $g_2$ that satisfy 
\begin{align*}
(-4 \pi^2 y^2 |n_2|^2 + y^2 \partial_y^2 - \lambda) g_1(y) = y^{\tfrac{1}{2}+\alpha} K_{\beta-1/2} (2 \pi  |n_2|  y)
\end{align*}
or 
\begin{align*}
(-4 \pi^2 y^2 |n_2|^2 + y^2 \partial_y^2 - \lambda) g_2(y) = y^{\tfrac{3}{2}-\alpha} K_{\beta-1/2} (2 \pi |n_2|  y),
\end{align*}
then the function 
\begin{align*}
	 \frac{2 \pi^\beta}{\Gamma(\beta) \zeta(2\beta)} |n|^{\beta-\tfrac{1}{2}} \sigma_{1-2\beta}(|n|) g_1(y) + \frac{ 2 \pi^{\beta+1/2} \Gamma(\alpha-\tfrac{1}{2}) \zeta(2\alpha-1)}{\Gamma(\alpha) \zeta(2\alpha) \Gamma(\beta) \zeta(2\beta)} 
	|n|^{\beta-\tfrac{1}{2}} \sigma_{1-2\beta}(|n|) g_2(y)
\end{align*}
solves \eqref{eq:fn1n2} for $n_1 = 0$ and $|n_2|=|n| $.

We assume that each of $g_i$ with $i=1,2$ can be represented as the following sum: 
 \[\sum_{j=0}^1 p^{j}(y)  K_j(2 \pi |n_2| y),\]
where $p^j$ is a polynomial in $y$ and $y^{-1}$. We note that for $g,h \in \R$,
\begin{align*}
	L_\lambda(x^g K_0(2 \pi |n_2| y)) =& ((g-1) g-\lambda ) x^g K_0(2 \pi |n_2| x)\\ & +2 \pi |n_2| (1-2 g) x^{g+1} K_1(2 \pi |n_2| x),\\
	L_\lambda(x^h K_1(2 \pi |n_2| y)) =& (-\lambda +(h-3) h+2) x^h K_1(2 \pi |n_2| x) \\ &+2 \pi |n_2| (1-2 h) x^{h+1} K_0(2 \pi |n_2| x).
\end{align*}
We assume 
\[
\max_{j \in \{0,1 \}} \deg p^{j} = M, \quad \min_{j \in \{0,1 \}} \deg p^{j} = m.
\] 
Then, in order to find coefficients of $p^j$ for $j=0,1$, we have to solve a system of linear equations with $2(M-m+1)$ variables, that are coefficients of the polynomials $p^j$ for $j=0,1$, 
and $2(M-m+1)+2$ equalities on coefficients at
\[
\bigcup_{\ell=m}^{M+2} \bigcup_{j=0}^1 \{ y^\ell K_j( 2 \pi |n_2| y) \}.
\]

\subsection{Solutions to \eqref{eq:fn1n2}  for $n_1 = n_2 = 0$}\label{sec:osaiujoiujwnkjh324o}
 We note that in order to solve \eqref{eq:fn1n2} for $n_1 = n_2 = 0$ it is sufficient to find solutions of  
\begin{align*}
(y^2 \partial_y^2 - \lambda) g(y) = y^{j_1+j_2}, \quad  j_1\in \{\alpha,1-\alpha\} \text{ and }j_2\in \{\beta,1-\beta\} .
\end{align*}
A particular solution can be easily constructed as products and sums of $\log(y)$ and polynomials in half-powers of $y$ and $y^{-1}$.

\subsection{Convergence for each Fourier mode}\label{sec:convergence}

We now examine the $n$-th Fourier mode $\widehat f_n(y).$ We note that in all cases, the particular solution $\displaystyle\sum_{n_1+n_2=n}\hat f^P_{n_1,n_2}(y)$ converges. However, the homogeneous part of the solution only converges for large enough $\lambda$. Fortunately, these cases correspond to the physically relevant cases considered in \cite{CGPWW2021}.

\subsubsection{Zero Fourier modes}

The zeroth Fourier mode is given by 
\begin{align*}
\widehat f_0(y)  =\widehat f_{0,0}(y) +\sum_{n_1\neq 0 }\widehat f_{n_1,-n_1}(y).
\end{align*}
Furthermore, the sum above is given by 
\begin{align}
\sum_{n_1\neq 0 }\widehat f_{n_1,-n_1}(y) & =\sum_{n_1\neq 0 }\widehat f^P_{n_1,-n_1}(y) +\sum_{n_1\neq 0 } \widehat f^H_{n_1,-n_1}(y)   \\ 
&= \sum_{n_1\neq 0 }\widehat f^P_{n_1,-n_1}(y) + y^{-r} \sum_{n_1\neq 0 } \alpha_{n_1,-n_1}  \nonumber
\end{align}
assuming both sums are convergent. The second equality follows from \eqref{eq:homogeneous_part123sdf23sdf}.

From \eqref{eq:6.2.4alphasum}, \eqref{eq:6.3.4alphasum},  \eqref{eq:7.9convergence},  \eqref{eq:8.11convergence},  and \eqref{eq:9.5convergence} we see that each
$\sum_{n_1\neq 0 } \alpha_{n_1,-n_1}
$ converges.
Each expression for
$\hat f^P_{n_1,-n_1}(y)$ is given by \eqref{eq:6.4fdef}, \eqref{eq:7.4fdef}, \eqref{eq:8.4fdef}, and \eqref{eq:9.4fdef}  and is exponentially suppressed as $y\to \infty$, as seen from the exponential decay of the modified Bessel functions of the second kind.

\subsubsection{Non-zero Fourier modes}

In order to show that the Fourier series is convergent, we first note that
\[
\widehat f_n(y) =\widehat f_{n,0}(y) + \widehat f_{0,n}(y)+\sum_{n_1=1}^{n-1}\widehat f_{n_1,n-n_1}(y) + \left( \sum_{n_1\geq n+1} +\sum_{n_1 \le -1}\right) \widehat f_{n_1,n-n_1}(y).
\]
We must verify that the last sum is convergent.  Note that by \eqref{eq:homogeneous_part123sdf},
\begin{align}
\sum_{n_1\geq n+1}\widehat f_{n_1,n-n_1}(y)& = \sum_{n_1\geq n+1}\widehat f^P_{n_1,n-n_1}(y) +\sum_{n_1\geq n+1} \widehat f^H_{n_1,n-n_1}(y)   \\
&= \sum_{n_1\geq n+1}\widehat f^P_{n_1,n-n_1}(y) + \sqrt{y} K_{r+1/2} (2 \pi |n|y) \sum_{n_1\geq n+1} \alpha_{n_1,n-n_1}  \nonumber
\end{align}
assuming both sums are convergent. 

From \eqref{asymptott_exp_6.2.3section},  \eqref{alpha_something_asymptwerwot6.3.3}, \eqref{alpha_something_asymptot7.7}, \eqref{alpha_something_asymptot6.3.3}, \eqref{alpha_something_asymptot6.3.sdfsdfsdq3}  we see that 
\[
\alpha_{n_1,n-n_1}=o(|n_1|^{-2}), \quad |n_1| \to \infty.
\]
Estimating the behavior of $\hat{f}^P_{n_1,n-n_1}(y)$  as $|n_1|\to \infty$ using \eqref{eq:6.2.3fasym}, \eqref{eq:6.2.3fasym56}, \eqref{eq:7.2.4fasym},  \eqref{eq:8.3.3fasym},  \eqref{eq:9.2.4fasym}, \eqref{eq:6.1.3fasym}, \eqref{eq:7.1.3fasym}, \eqref{eq:8.1.3fasym}, \eqref{eq:8.2.3fasym},  and \eqref{eq:9.1.4fasym},  we see that the contribution from the modified Bessel functions, ${K_i(2\pi |n_1|y)K_j(2\pi |n-n_1|y)}$, exponentially suppresses these terms as $|n_1|$ gets large. We treat the term $\sum_{n_1 \le -1}$ in the same manner.

\section{Differential Galois theory}

This article gives explicit solutions for specific combination of $\alpha$, $\beta$ and $\lambda = r(r+1)$ with $r>0$ listed in Section \ref{Section:main_results_discussion}; however, experimentally we were able to compute solutions for other combinations of $\alpha$, $\beta$ and $\lambda$ as well. 

Finding solutions for large $\alpha, \beta$ and $\lambda$ involves solving a systems of linear equation for a large number of variables. This becomes  computationally challenging, even though the corresponding matrices are band matrices. We obtained that a particular solution of \eqref{eq:fn1n2} is of the form \eqref{wqLkjlkjlkjwewerwempjk}, at least, in the cases where 
\begin{equation}\label{special_set_of_solvable_parameters}
\alpha, \beta \in \mathbb{Z}+\tfrac{1}{2}, \quad \alpha + \beta + r \in 2 \mathbb{Z}, \quad |\alpha-\beta|  < r
\end{equation}
and
\[
\alpha < 30, \quad \beta < 30, \quad  r < 15.
\]
Note that the functions that appear in  \cite[(2.11)]{CGPWW2021} satisfy the condition above.

Below, we make a conjecture that the solution are ``nice'' if $(\alpha, \beta, \gamma)$ belongs to \eqref{special_set_of_solvable_parameters},  regardless of how large each parameter may be. Discussing in which way they are ``nice" would require some basic facts from the differential Galois theory, that we outline as follows. 

The fundamental system of the homogeneous solution of~\eqref{eq:fn1n2hom} is well-known for any values of $r \in \mathbb{R}$ and involves the modified Bessel functions (see \eqref{homogeneous_equation_solutioniuz2iu3j4h}). Moreover, it is possible to  show\footnote{For $J$-Bessel functions, the proof can be found in \cite[Appendix]{Kolchin}; we can obtain the same statement for $K$ and $I$-Bessel functions by exploiting formulas relating Bessel functions to each other.} that modified Bessel equations, $K_\eta$, can be expressed via elementary functions if and only if their index, $\eta$, belongs to $\tfrac{1}{2} + \mathbb{Z}$. In our notations, this corresponds to demanding $r \in \mathbb{Z}$. 

Recall \cite[Chapter 3]{Khovansky} that a \textit{differential field}, $K$, is a field together with a \textit{derivation} (i.e. an additive map that satisfies the Leibniz rule, $\partial(ab) = \partial(a) b + a \partial(b)$). An elementary example of a differential field would be the field $\mathbb{C}(t)$ of rational functions over $\mathbb{C}$ together with an usual operation of differentiation. Solutions of the type \eqref{wqLkjlkjlkjwewerwempjk} belong to a particular object in differential Galois theory -- namely, they belong to a certain Picard-Vessiot extension of a differential field. A differential field $P$ is called a {\it Picard-Vessiot extension of the field $K$}, if there exists a linear differential equation with coefficients in $K$ such that $P$ is obtained from~$K$ by adjoining a fundamental system of solutions of this differential equation. 

When  $n_1,n_2 \in \mathbb{Z}\setminus \{0\}$ with $n:=n_1+n_2 \neq 0$, the differential field, $P$, that we are interested in can be obtained by adjoining to $\mathbb{C}(t)$ solutions of the equations 
\begin{equation}\label{picard-vessiot1}
(y^2 \partial_y^2 - 4 \pi^2 n_1^2 y^2+1/4) f(y) = 0,
\end{equation}
\begin{equation}\label{picard-vessiot2}
(y^2 \partial_y^2 - 4 \pi^2 n_2^2 y^2+1/4) f(y) = 0,
\end{equation}
and 
\begin{equation}\label{picard-vessiot3}
    (y^2 \partial_y^2-r(r+1) - 4 \pi^2 (n_1+n_2)^2 y^2) f(y) = 0.
\end{equation}
 We note that $P$ contains $\sqrt{y}\,K_0(2 \pi |n_1| y)$ and $\sqrt{y}\, K_0(2 \pi |n_2| y)$ by construction. Since it is an extension of $\mathbb{C}(t)$, it also contains any sum of the type
\[
\sum_{j=0}^{n} a_j y^{j+1/2} K_0(2 \pi |n_1|y), \quad n \in \mathbb{N}, a_j \in \mathbb{C}.
\]
Moreover, the recurrence relation between $K_0$ and $K_1$ implies, that $P$ contains  any sums of type 
\[
\sum_{j=0}^{n} b_j y^{j+1/2} K_1(2 \pi |n_1|y), \quad n \in \mathbb{N}, b_j \in \mathbb{C}.
\]
Thus, we obtain solutions of the type \eqref{wqLkjlkjlkjwewerwempjk} belong to the field $P$ that we have just constructed. 

On the other hand, we can reformulate the inhomogeneous differential equation \eqref{eq:fn1n2} as the following homogeneous differential equation of the third order on $g_{n_1,n_2}(y)$:
\begin{equation}\label{eq:inhom_vs_hom}
s_{n_1,n_2}(y) \frac{\partial}{\partial y}\left(\frac{(y^2 \partial_y^2-\lambda - 4 \pi^2 (n_1+n_2)^2 y^2) g_{n_1,n_2}(y)}{s_{n_1,n_2}(y)} \right) = 0.
\end{equation}
We note that for any solution \eqref{eq:inhom_vs_hom} there is a constant $c$ such that $g_{n_1,n_2}(y)$
is a solution of 
\[
(y^2 \partial_y^2-\lambda - 4 \pi^2 (n_1+n_2)^2 y^2) g_{n_1,n_2}(y) = c s_{n_1,n_2}(y). 
\]
And, on the other hand, every $\hat{f}_{n_1,n_2}(y)$ is also solution of \eqref{eq:inhom_vs_hom}. 

We will prove in this article, that for certain $\alpha, \beta, \lambda$, the function $g$ belongs to $P$ simply by providing an explicit solution, thus, the Galois group of the differential equation is trivial. On the other hand, the numerical simulations  suggest that it is also solvable for all relatively small  $\alpha, \beta, \lambda $ satisfying \eqref{special_set_of_solvable_parameters}. 

\begin{conj}
Let $P$ be a Picard-Vessiot extension of the field of rational functions over  $\mathbb{C}$, obtained by adjoining solutions of \eqref{picard-vessiot1}-\eqref{picard-vessiot3}. Then the Galois group of \eqref{eq:inhom_vs_hom} is trivial in the category of algebraic groups if  $\alpha, \beta, r$, however big they are, satisfy~\eqref{special_set_of_solvable_parameters}.
\end{conj}
However, proving or disproving this conjecture is beyond the scope of this paper. There are certain related results, see \cite{Minchenko}, where the authors present the algorithm that calculates the differential Galois group of a third-order homogeneous linear differential equation. 

\section{$\alpha = \beta=3/2$}\label{sec:3/2}

In this section we solve 
\[
(\Delta - \lambda) f(z) = - (2 \zeta(3) E_{3/2}(z))^2, \quad z=x+iy \in \mathfrak{H}
\]
 for 
 \[f(z) = \sum_{n \in \Z} \sum_{n_1+n_2=n} \hat{f}_{n_1,n_2}(y) e^{2 \pi i n x}\]
 in terms of $\hat{f}_{n_1,n_2}(y) = \hat{f}^P_{n_1,n_2}(y)+\hat{f}^H_{n_1,n_2}(y).$
 
When $n_1=n_2=0$, $\hat{f}_{0,0}(y)$ contains no $K$-Bessel or divisor functions and is given by a a polynomial in $y$ and $1/y$ below.  For $n_1n_1=0$ but not both zero, 
 \begin{equation}\label{eq:6.1fdef}
\hat{f}^P_{0,n}(y) =\hat{f}^P_{n,0}(y)=- 16 \pi \frac{\sigma_2(|n|)}{|n|} \sum_{i=0,1} \nu_{i}(n, y) K_i(2 \pi |n| y),
\end{equation}
for $n_1 n_2\neq 0$ and $n_1+n_2\neq 0$,
\begin{equation}\label{eq:6.3fdef}
\hat{f}^P_{n_1,n_2}(y) = - 64 \pi^2 \frac{\sigma_2(|n_1|)\sigma_2(|n_2|)}{|n_1 n_2|}\sum_{i,j=0,1} \eta_{i,j}(n_1, n_2, y) K_i(2 \pi |n_1| y) K_j(2 \pi |n_2| y),
\end{equation} and 
for $n_1 =-n_2$,
\begin{equation}\label{eq:6.4fdef}
\hat{f}^P_{-n_2,n_2}(y) = - 64 \pi^2 \frac{\sigma_2(|n_2|)\sigma_2(|n_2|)}{|n_2 |^2}\sum_{(i,j)\in \{(0,0),(0,1),(1,1)\} } \mu_{i,j}(n_2, y) K_i(2 \pi |n_2| y) K_j(2 \pi |n_2| y)
\end{equation}
where $\nu_{i}, \eta_{i,j}$ and $\mu_{i,j}$ defined below depending on each value of $\lambda$.

\subsection{$\lambda = 30$}\label{sec:3/2,12}
This case corresponds to \cite[Section C.3.1]{CGPWW2021} with $r=5$.

\subsubsection{$n_1=0$ and $n_2=0$.}
Any solution of \eqref{eq:fn1n2} for $n_1=n_2=0$ is equal to 
\[
\hat{f}_{0,0}(y) = \frac{105 y^4 \zeta (3)^2+56 \pi ^2 y^2 \zeta (3)+10 \pi ^4}{630 y}+ \frac{\alpha_{0,0}}{y^5}+\beta_{0,0}\,y^6
\]
for some $c_1, c_2 \in \C$. Its asymptotic behavior for $y \to 0$ can be described by 
\[
\hat{f}_{0,0}(y) =\alpha_{0,0} y^{-5}.
\]
At this moment of time, we do not choose $\alpha_{0,0}$ -- that will be reserved for Section \ref{n1n2_term_sldkjf98234kjhbnbsoodoooSDD}. Our goal would be to choose $\alpha_{0,0}$ in such a way that 
\[
\sum_{n} \alpha_{n,-n} = 0.
\]
In our notation and after the evaluation of the Riemann zeta function at even integers, the first three summands of the first line of \cite[(C.27)]{CGPWW2021}  read 
\[
\frac{y^3 \zeta (3)^2}{6}+\frac{4}{45} \pi ^2 y \zeta (3)+\frac{\pi ^4}{63 y},
\]
that coincides with our result.
\subsubsection{$n_1n_2=0$ but not both zero} 
Though this case of $\alpha=\beta=3/2$ and $\lambda =30$ is generally addressed in  \cite{CGPWW2021}, we note that the term  $\hat{f}^P_{0,n}(y)$ was not found explicitly. For
$\hat{f}^P_{0,n}(y)$ as in \eqref{eq:6.1fdef}, we have 
\begin{align*}
 \nu_0(n,y)&=\sgn(n)\left[-\zeta(3)\left(
\frac{126y^{-3}}{n^{5} \pi^{5}}  
+\frac{35y^{-1}}{n^{3} \pi^{3}}  
+\frac{y}{2 n \pi}  \right) + 2\zeta(2) \left( \frac{3y^{-3}}{5 n^{3} \pi^{3}}  
+\frac{y^{-1}}{6 n \pi}\right)\right],
\\ \nu_1(n,y)&=
-\zeta(3) \left(\frac{126y^{-4}}{n^{6} \pi^{6}}  
+\frac{98y^{-2}}{n^{4} \pi^{4}}  
+\frac{15}{2 n^{2} \pi^{2}} \right) +2\zeta(2) \left(\frac{3y^{-4}}{5 n^{4} \pi^{4}}  
+\frac{7y^{-2}}{15 n^{2} \pi^{2}}  \right).
\end{align*}
Its asymptotic behavior for $y \to 0$ can be described as follows 
\begin{equation}\label{eq:6.2.2fasym}
\hat{f}^P_{0,n}(y) = -\frac{48  \sigma _2(|n|) \left(\pi ^2 \zeta (2) n^2-105 \zeta (3)\right)}{5 \pi ^6 n^8 y^5}+\frac{32 \left(\pi ^2 \zeta (2) n^2-105 \zeta (3)\right) \sigma _2(|n|)}{15 \pi ^4 n^6 y^3}+O\left(\frac{1}{y^2}\right).
\end{equation}\label{eq:6.2.2alphaasym}
The unique  choice of $\alpha_{n,0}=\alpha_{0,n}$ that gets rid of the $O(y^{-5})$-term in $\hat{f}^P_{0,n}(y)+\hat{f}^H_{0,n}(y)$ is
\begin{align}
\alpha_{n,0}=\alpha_{0,n}= 
-\frac{1024  \left(\pi ^2 \zeta (2) n^2-105 \zeta (3)\right) \sigma _2(|n|)}{1575 \pi  \left|n\right|{}^{5/2}}.
\end{align}

\subsubsection{$n_1n_2\neq 0$ and $n_1+n_2\neq 0$.}

 In \cite[p.\,46]{CGPWW2021}, many terms\footnote{Powers of around $D$-instanton contributions with $n \neq 0$, that include the instanton sectors of $(n_1, n_2) =(2,0),(1,1),(1,-2),(1,-3),(2,-3)$ were found.} in the perturbative expansion  $\hat{f}^P_{n_1,n_2}(y)$ were evaluated. 
 However, these values were not explicitly written or evaluated in full in \cite{CGPWW2021}.
For
$\hat{f}^P_{n_1,n_2}(y) $ as in \eqref{eq:6.3fdef}, we have  
\begin{align*}
	\eta_{0,0}=
	& \sgn(n_1n_2) \Big[\frac{y^{-3}\, n_{1} n_{2}}{ \left(n_{1} + n_{2}\right)^{10}}\frac{126}{\pi^{4}}  \left(n_{1}^{4} - 6 n_{1}^{3} n_{2} + 10 n_{1}^{2} n_{2}^{2} - 6 n_{1} n_{2}^{3} + n_{2}^{4}\right)\\
	& \ \ \ \ \ \ \ \ \   \ \ \ \ \ \ \ +\frac{y^{-1}\, n_{1} n_{2}}{ \left(n_{1} + n_{2}\right)^{8}}\frac{2}{5 \pi^{2}}  \left(89 n_{1}^{4} - 792 n_{1}^{3} n_{2} + 1598 n_{1}^{2} n_{2}^{2} - 792 n_{1} n_{2}^{3} + 89 n_{2}^{4}\right)\\
	&\ \ \ \ \ \ \ \ \  \ \ \ \  \ \ \  +\frac{y\, n_{1} n_{2}}{ \left(n_{1} + n_{2}\right)^{6}}\frac{2}{ 15}  \left(5 n_{1}^{4} - 92 n_{1}^{3} n_{2} + 190 n_{1}^{2} n_{2}^{2} - 92 n_{1} n_{2}^{3} + 5 n_{2}^{4}\right)\Big],\\
	\eta_{0,1}=
	&\sgn(n_1)\Big[\frac{y^{-4}\, n_{1} n_{2}^{2}}{\left(n_{1} + n_{2}\right)^{11}}\frac{126}{\pi^{5} }  \left(- n_{1}^{3} + 5 n_{1}^{2} n_{2} - 5 n_{1} n_{2}^{2} + n_{2}^{3}\right)\\
	&\ \ \ \ \ \  \ \ \ \ \ \ \ \ +\frac{y^{-2}\, n_{1}}{5 \pi^{3} \left(n_{1} + n_{2}\right)^{9}}  \left(3 n_{1}^{5} + 99 n_{1}^{4} n_{2} - 2728 n_{1}^{3} n_{2}^{2} + 6512 n_{1}^{2} n_{2}^{3} \right.\\
	 & \ \ \ \ \ \ \ \ \ \ \ \ \ \ \ \ \ \ \ \ \ \ \ \ \ \ \ \ \ \ \ \ \ \ \ \ \  \ \ \ \  \ \ \ \ \left.- 3611 n_{1} n_{2}^{4} + 493 n_{2}^{5}\right)\\
	&  \ \ \ \ \ \  \ \ \ \  \ \ \ \ + \frac{n_{1}}{30 \pi \left(n_{1} + n_{2}\right)^{7}} \left(5 n_{1}^{5} + 147 n_{1}^{4} n_{2} - 2614 n_{1}^{3} n_{2}^{2} + 5726 n_{1}^{2} n_{2}^{3}\right.\\
	 & \ \ \ \ \ \ \ \ \ \ \ \ \ \ \ \ \ \ \ \ \ \ \ \ \ \ \ \ \ \ \ \ \ \ \ \ \ \ \ \ \  \ \ \ \ \left.- 2799 n_{1} n_{2}^{4} + 239 n_{2}^{5}\right)\Big],
	 \end{align*}
	 
	 \begin{align*}
	\eta_{1,0}=
	&\sgn(n_2)\Big[\frac{y^{-4}\, n_{1}^{2} n_{2}}{\left(n_{1} + n_{2}\right)^{11}}\frac{126}{\pi^{5} }  \left(n_{1}^{3} - 5 n_{1}^{2} n_{2} + 5 n_{1} n_{2}^{2} - n_{2}^{3}\right)\\
	& \ \ \ \ \ \ +\frac{y^{-2}\, n_{2} }{5 \pi^{3} \left(n_{1} + n_{2}\right)^{9}} \left(493 n_{1}^{5} - 3611 n_{1}^{4} n_{2} + 6512 n_{1}^{3} n_{2}^{2} - 2728 n_{1}^{2} n_{2}^{3} \right.\\
	 & \ \ \ \ \ \ \ \ \ \ \ \ \ \ \ \ \ \ \ \ \ \ \ \ \ \ \ \ \ \ \ \ \ \ \ \ \  \ \ \ \  \ \ \ \ \left.+ 99 n_{1} n_{2}^{4} + 3 n_{2}^{5}\right)\\
	& \ \ \ \ \ \ +\frac{n_{2}}{30 \pi \left(n_{1} + n_{2}\right)^{7}} \left(239 n_{1}^{5} - 2799 n_{1}^{4} n_{2} + 5726 n_{1}^{3} n_{2}^{2} - 2614 n_{1}^{2} n_{2}^{3}\right.\\
	 & \ \ \ \ \ \ \ \ \ \ \ \ \ \ \ \ \ \ \ \ \ \ \ \ \ \ \ \ \ \ \ \ \ \ \ \ \  \ \ \ \  \ \ \ \ \left. + 147 n_{1} n_{2}^{4} + 5 n_{2}^{5}\right)\Big],\\
	\eta_{1,1}=
	&\frac{y^{-3}}{5 \pi^{4} \left(n_{1} + n_{2}\right)^{10}}(3 n_{1}^{6} + 102 n_{1}^{5} n_{2} - 3399 n_{1}^{4} n_{2}^{2} + 8124 n_{1}^{3} n_{2}^{3} - 3399 n_{1}^{2} n_{2}^{4}\\
	 & \ \ \ \ \ \ \ \ \ \ \ \ \ \ \ \ \ \ \ \ \ \ \ \ \ \ \ \ \ \ \ \ \ \ \ \ \  \ \ \ \  \ \ \ \ + 102 n_{1} n_{2}^{5} + 3 n_{2}^{6})\\
	&+\frac{y^{-1}}{15 \pi^{2} \left(n_{1} + n_{2}\right)^{8}}\left(7 n_{1}^{6} + 220 n_{1}^{5} n_{2} - 4727 n_{1}^{4} n_{2}^{2} + 10280 n_{1}^{3} n_{2}^{3} - 4727 n_{1}^{2} n_{2}^{4}\right.\\
	 & \ \ \ \ \ \ \ \ \ \ \ \ \ \ \ \ \ \ \ \ \ \ \ \ \ \ \ \ \ \ \ \ \ \ \ \ \  \ \ \ \  \ \ \ \ \left. + 220 n_{1} n_{2}^{5} + 7 n_{2}^{6} \right)\\
	&  +\frac{y\,  n_{1} n_{2} }{ \left(n_{1} + n_{2}\right)^{6}}\frac{2}{15}\left(5 n_{1}^{4} - 92 n_{1}^{3} n_{2} + 190 n_{1}^{2} n_{2}^{2} - 92 n_{1} n_{2}^{3} + 5 n_{2}^{4}\right).
\end{align*}
Note that 
\begin{align}\label{eq:6.2.3fasym}
	\hat{f}^P_{n_1,n_2}(y)=& - 16 y^{-5}  \frac{\sigma_2(|n_1|)\sigma_2(|n_2|)}{5 |n_1 n_2|^2  \pi^{4} \left(n_{1} + n_{2}\right)^{11}} \Big(3 n_{1}^{7} + 105 n_{1}^{6} n_{2} + 1260 n_{1}^{5} n_{2}^{2} \log{\left(|n_{1}| \pi \right)}\nonumber \\& - 1260 n_{1}^{5} n_{2}^{2} \log{\left(|n_{2}| \pi \right)} - 3297 n_{1}^{5} n_{2}^{2} - 6300 n_{1}^{4} n_{2}^{3} \log{\left(|n_{1}| \pi \right)} + 6300 n_{1}^{4} n_{2}^{3} \log{\left(|n_{2}| \pi \right)}\nonumber\\ & + 4725 n_{1}^{4} n_{2}^{3}+ 6300 n_{1}^{3} n_{2}^{4} \log{\left(|n_{1}| \pi \right)} - 6300 n_{1}^{3} n_{2}^{4} \log{\left(|n_{2}| \pi \right)} + 4725 n_{1}^{3} n_{2}^{4} \\ &- 1260 n_{1}^{2} n_{2}^{5} \log{\left(|n_{1}| \pi \right)} + 1260 n_{1}^{2} n_{2}^{5} \log{\left(|n_{2}| \pi \right)} - 3297 n_{1}^{2} n_{2}^{5} + 105 n_{1} n_{2}^{6} + 3 n_{2}^{7}\Big)\nonumber\\
	&+o(y^{-5}).\nonumber
\end{align}

We recall that by \eqref{asymptotics_of_homogeneous_solutionnasokdijfok}, \[\hat{f}^H_{n_1,n_2}(y) = \alpha_{n_1, n_2} \sqrt{y} K_{5+\frac{1}{2}}(2 \pi  |n_1+n_2| y) = \alpha_{n_1,n_2} \frac{945}{64 \pi ^5 |n_1+n_2|^{11/2} y^5}+o(y^{-5}).\]  
Comparing the right hand sides of the  two previous formulas, we obtain that there is a unique choice of $\alpha_{n_1, n_2}$ that guarantees that  $$\hat{f}_{n_1,n_2}(y) = \hat{f}^P_{n_1,n_2}(y)+ \hat{f}^H_{n_1,n_2}(y) = o(y^{-5}),$$ as  $y \to 0$ given by\newpage
\begin{align}\label{alphan1n2_slodikjuf9o8w794r8w7u}
	\alpha_{n_1,n_2}=&  \frac{1024 \pi \sigma_2(|n_1|)\sigma_2(|n_2|) \sgn(n_1+n_2) }{ 4725 |n_1 n_2|^{2}   \left|n_{1} + n_{2}\right|^{11/2}} \nonumber \\ & \times \Big(3 n_{1}^{7} + 105 n_{1}^{6} n_{2} + 1260 n_{1}^{5} n_{2}^{2} \log{\left(|n_{1}| \pi \right)}\nonumber \\& - 1260 n_{1}^{5} n_{2}^{2} \log{\left(|n_{2}| \pi \right)} - 3297 n_{1}^{5} n_{2}^{2} - 6300 n_{1}^{4} n_{2}^{3} \log{\left(|n_{1}| \pi \right)} + 6300 n_{1}^{4} n_{2}^{3} \log{\left(|n_{2}| \pi \right)}\nonumber\\ & + 4725 n_{1}^{4} n_{2}^{3}+ 6300 n_{1}^{3} n_{2}^{4} \log{\left(|n_{1}| \pi \right)} - 6300 n_{1}^{3} n_{2}^{4} \log{\left(|n_{2}| \pi \right)} + 4725 n_{1}^{3} n_{2}^{4} \\ &- 1260 n_{1}^{2} n_{2}^{5} \log{\left(|n_{1}| \pi \right)} + 1260 n_{1}^{2} n_{2}^{5} \log{\left(|n_{2}| \pi \right)} - 3297 n_{1}^{2} n_{2}^{5} + 105 n_{1} n_{2}^{6} + 3 n_{2}^{7}\Big).\nonumber
\end{align}

Moreover, it is not complicated to check that for fixed values of $n$,
\begin{equation}\label{asymptott_exp_6.2.3section}
\alpha_{n-n_1, n_1} = O(|n_1|^{-4}), \quad |n_1| \to \infty.
\end{equation}

\subsubsection{ $n_1= -n_2$.}\label{n1n2_term_sldkjf98234kjhbnbsoodoooSDD} This case has been considered in \cite[(C.29)]{CGPWW2021}. 

For $
\hat{f}^P_{-n_2,n_2}(y) $ as in \eqref{eq:6.4fdef}
we have  
\begin{align*}
\mu_{0,0} &=
\frac{y^{-1}}{110 n_{2}^{2} \pi^{2}}  
+\frac{y}{110}  
-\frac{ 8 n_{2}^{2} \pi^{2}y^{3}}{1155}  
-\frac{ 512 n_{2}^{4} \pi^{4}y^{5}}{17325}  
+\frac{16384 n_{2}^{6} \pi^{6}y^{7}}{51975},  
\\
\mu_{0,1} &= \sgn(n_2) \left( 
\frac{y^{-2}}{55 n_{2}^{3} \pi^{3}}  
+\frac{3}{110 n_{2} \pi}  
+\frac{4 n_{2} \pi y^{2}}{385}  
+\frac{256 n_{2}^{3} \pi^{3}y^{4}}{17325}  
+\frac{8192 n_{2}^{5} \pi^{5}y^{6}}{51975}  \right) ,
\\
\mu_{1,1} &=
\frac{y^{-3}}{110 n_{2}^{4} \pi^{4}}  
+\frac{y^{-1}}{55 n_{2}^{2} \pi^{2}}  
-\frac{17y}{770}  
-\frac{ 8 n_{2}^{2} \pi^{2}y^{3}}{1925}  
-\frac{ 512 n_{2}^{4} \pi^{4}y^{5}}{51975}  
-\frac{16384 n_{2}^{6} \pi^{6}y^{7}}{51975}.
\end{align*}
We note that 
\begin{equation}\label{eq:6.2.4fasym}
\hat{f}^P_{-n_2,n_2}(y) = - 8 \frac{\sigma_2(|n_2|)^2}{|n_2|^2} \cdot \left(\frac{1}{55 \pi ^4 n_2^6 y^5}+O\left(\frac{1}{y}\right) \right).
\end{equation}
The unique choice of $\alpha_{-n_2,n_2}$ that gets rid of the $O(y^{-5})$ term in the expression above is 
\begin{equation}\label{alpha_n2-n2in610}
\alpha_{-n_2,n_2}= \frac{8 \sigma_2(|n_2|)^2}{55 \pi ^4  |n_2|^8}.
\end{equation}
Summing it up and using \eqref{Ramanujan2} for $a=2, b=2$ and $s=8$, we obtain 
\[
\left. \sum_{n_2=-\infty, n_2 \neq 0}^\infty \frac{\sigma_{2}(|n_2|)^2}{|n_2|^8}= 2 \frac{\zeta(s) \zeta(s-a) \zeta(s-b) \zeta(s-a-b)}{\zeta(2s-a-b)} \right|_{a=2,b=2,s=8}=\frac{143 \pi ^{12}}{58769550},
\]
and thus 
\begin{equation}\label{eq:6.2.4alphasum}
\sum_{n_2=-\infty, n_2 \neq 0}^\infty \alpha_{-n_2,n_2} = \frac{52 \pi ^8}{146923875} =  \frac{104 \zeta (8)}{31095}.
\end{equation}
This means, that $\sum_{n_2=-\infty, n_2 \neq 0}^\infty \hat{f}_{-n_2,n_2}^H(y) =   \frac{104 \zeta (8)}{31095} y^{-5} $. Motivated by the desire to have the contribution from the homogeneous elements to be equal to zero, we obtain $\alpha_{0,0} = -\frac{104 \zeta (8)}{31095}$. This matches\footnote{Up to sign} the last summand in the first line of  \cite[(C.27)]{CGPWW2021}.

\subsection{$\lambda = 56$}
This case corresponds to \cite[Section C.3.1]{CGPWW2021} with $r=7$.

\subsubsection{$n_1=0$ and $n_2=0$.}
The solution of \eqref{eq:fn1n2} for $n_1=n_2=0$ is equal to 
\[
\hat{f}_{0,0}(y) = \frac{3402 y^4 \zeta (3)^2+2025 \pi ^2 y^2 \zeta (3)+350 \pi ^4}{42525 y}+\frac{c_1}{y^7}+c_2y^8 \
\]
for some $c_1, c_2 \in \C$. Its asymptotic behavior for $y \to 0$ can be described by 
$
\hat{f}_{0,0}(y) = \frac{c_1}{y^7}.
$
We do not specify the choice of $c_1$ for the moment, but we can set $c_2=0$ so that the $O(y^8)$-term vanishes. We note that the first three summands\footnote{However, our choice of $c_1$ does not coincide with \cite[(C.27)]{CGPWW2021} coefficient at $y^{-7}$. } of $\hat{f}_{0,0}(y)$ coincide with the first three terms of the second line of \cite[(C.27)]{CGPWW2021}.

\subsubsection{$n_1n_2=0$ but not both zero} Though this case of $\alpha=\beta=3/2$ and $\lambda =56$ is generally addressed in  \cite{CGPWW2021}, we note that the term  $\hat{f}^P_{0,n}(y)$ was not found explicitly. This term is given by
\[
\hat{f}^P_{0,n}(y) =\hat{f}^P_{n,0}(y)=- 16 \pi \frac{\sigma_2(|n|)}{|n|} \sum_{i,j=0,1} \nu_{i}(n, y) K_i(2 \pi |n| y) ,
\]
with 
\begin{align*}
\nu_0(n, y)&=
\sgn(n)\Big[-\zeta(3)\left(\frac{30888y^{-5}}{5 n^{7} \pi^{7}}  
+\frac{10692y^{-3}}{5 n^{5} \pi^{5}}  
+\frac{126y^{-1}}{n^{3} \pi^{3}}  
+\frac{y}{2 n \pi} \right)\\
 & \ \ \ \ \ \ \ \ \ \ \ \ \ \ \ \ \ \ \ \ \ \ \ \ \ \ \ \ \ \ \ \ \ \ +2\zeta(2)\left( \frac{286y^{-5}}{35 n^{5} \pi^{5}}  
+\frac{99y^{-3}}{35 n^{3} \pi^{3}}  
+\frac{y^{-1}}{6 n \pi}  \right) \Big],
\\ \nu_1(n, y)&=
-\zeta(3) \left( -\frac{30888y^{-6}}{5 n^{8} \pi^{8}}  
+\frac{26136y^{-4}}{5 n^{6} \pi^{6}}  
+\frac{3402y^{-2}}{5 n^{4} \pi^{4}}  
+\frac{14}{n^{2} \pi^{2}} \right)\\
 & \ \ \ \ \ \ \ \ \ \ \ \ \ \ \ \ \ \ \ \ \ \ \ \ \ \ \ \ \ \ \ \ \ \ + 2\zeta(2) \left(\frac{286y^{-6}}{35 n^{6} \pi^{6}}  
+\frac{242y^{-4}}{35 n^{4} \pi^{4}}  
+\frac{9y^{-2}}{10 n^{2} \pi^{2}} \right).
\end{align*}
The asymptotic expansion 
is
\begin{equation}
\hat{f}^P_{0,n}(y) =\hat{f}^P_{n,0}(y)=- 16 \frac{\sigma_2(|n|)}{|n|} \left( 
\frac{286 \left(\pi ^2 \zeta (2) n^2-378 \zeta (3)\right)}{35 \pi ^8 |n|^9 y^7}+O\left(\frac{1}{y^5}\right) \right).
\end{equation}
There exists a unique choice of $\alpha_{0,n}=\alpha_{n,0}$ such that $\hat{f}^P_{0,n}(y) + \hat{f}^H_{0,n}(y) $ and $\hat{f}^P_{n,0}(y) + \hat{f}^H_{n,0}(y) $  are of order $o(y^{-7})$. \
More precisely, 
\begin{align*}
   \hat{f}^H_{0,n}(y) =  \sqrt{y} \alpha_{0,n} K_{7+\frac{1}{2}}(2 \pi  |n| y) = \frac{135135 \alpha _{n,0}}{256 \pi ^7 |n|^{15/2} y^7} + O(y^{-5}),
\end{align*}
and thus we may set 
\begin{align*}
    \alpha_{0,n} = \alpha_{n,0}  = \frac{8192 \left(\pi ^2 \zeta (2) n^2-378 \zeta (3)\right) \sigma _2(|n|)}{33075 \pi  |n|^{5/2}}.
\end{align*}

\subsubsection{$n_1n_2\neq 0$ and $n_1+n_2\neq 0$.}
Though this case of $\alpha=\beta=3/2$ and $\lambda =56$ is generally addressed in  \cite{CGPWW2021}, we note that the term  $\hat{f}^P_{0,n}(y)$ was not found explicitly. For
$\hat{f}^P_{n_1,n_2}(y) $ as in \eqref{eq:6.3fdef}, we have  
\begin{align*}
\eta_{0,0}=
&\sgn(n_1n_2)\Big[\frac{y^{-5}\,n_{1} n_{2}}{ \left(n_{1} + n_{2}\right)^{14}}\frac{10296  }{5 \pi^{6}} (3 n_{1}^{6}- 38 n_{1}^{5} n_{2}+ 140 n_{1}^{4} n_{2}^{2} - 210 n_{1}^{3} n_{2}^{3} + 140 n_{1}^{2} n_{2}^{4} \\
 &   \ \ \ \   \ \ \ \ \ \ \ \  \ \ \ \ \ \ \ \   \ \ \ \ \ \ \ \ \ \ \ \ \ \ \ \   \ \ \ \ \ \ \ \ \ \ \ \ \ \ \ \  - 38 n_{1} n_{2}^{5} + 3 n_{2}^{6})\\
&\ \ \ +\frac{y^{-3}\,n_{1} n_{2} }{ \left(n_{1} + n_{2}\right)^{12}}\frac{22}{175 \pi^{4}}  \left(17075 n_{1}^{6} - 260700 n_{1}^{5} n_{2} + 1170813 n_{1}^{4} n_{2}^{2} - 1907624 n_{1}^{3} n_{2}^{3} \right.\\
 &   \ \ \ \   \ \ \ \ \ \ \ \  \ \ \ \ \ \ \ \   \ \ \ \ \ \ \ \ \ \ \ \ \ \ \ \   \ \ \ \ \ \ \ \ \ \ \ \     \left.
+ 1170813 n_{1}^{2} n_{2}^{4} - 260700 n_{1} n_{2}^{5} + 17075 n_{2}^{6}\right)\\
& \ \ \ +\frac{y^{-1}\,n_{1} n_{2}}{175 \pi^{2} \left(n_{1} + n_{2}\right)^{10}}  \left(22545 n_{1}^{6} - 445070 n_{1}^{5} n_{2} + 2255599 n_{1}^{4} n_{2}^{2} - 3778788 n_{1}^{3} n_{2}^{3} \right. \\
 & \ \ \ \   \ \ \ \ \ \ \ \  \ \ \ \ \ \ \ \   \ \ \ \ \ \ \ \ \ \ \ \ \ \ \ \   \ \ \ \ \ \ \ \ \ \ \ \ \left.   + 2255599 n_{1}^{2} n_{2}^{4} - 445070 n_{1} n_{2}^{5} + 22545 n_{2}^{6}\right)\\
&\ \ \  +\frac{2\,y\, n_{1} n_{2}}{525 \left(n_{1} + n_{2}\right)^{8}}   \left(175 n_{1}^{6} - 6510 n_{1}^{5} n_{2} + 35745 n_{1}^{4} n_{2}^{2} - 61572 n_{1}^{3} n_{2}^{3} + 35745 n_{1}^{2} n_{2}^{4}  \right.\\
 &  \ \ \ \   \ \ \ \ \ \ \ \  \ \ \ \ \ \ \ \   \ \ \ \ \ \ \ \ \ \ \ \ \ \ \ \   \ \ \ \ \ \ \ \ \ \ \ \ \  \ \ \ \ \ \ \    \left. - 6510 n_{1} n_{2}^{5} + 175 n_{2}^{6}\right))\Big]\\
\eta_{0,1}=
&\sgn(n_1)\Big[\frac{y^{-6}\,n_{1} n_{2}^{2}}{ \left(n_{1} + n_{2}\right)^{15}}\frac{10296 }{5 \pi^{7}} \left(- 3 n_{1}^{5} + 35 n_{1}^{4} n_{2} - 105 n_{1}^{3} n_{2}^{2} + 105 n_{1}^{2} n_{2}^{3} - 35 n_{1} n_{2}^{4} + 3 n_{2}^{5}) \right. \\
&\ \ \ +\frac{y^{-4}\,n_{1} }{ \left(n_{1} + n_{2}\right)^{13}}\frac{22}{175 \pi^{5}}  (65 n_{1}^{7} + 3965 n_{1}^{6} n_{2} - 283107 n_{1}^{5} n_{2}^{2} + 1865089 n_{1}^{4} n_{2}^{3}  \\
 & \ \ \ \   \ \ \ \ \ \ \ \  \ \ \ \ \ \ \ \   \ \ \ \ \ \ \ \ \ \   - 3515111 n_{1}^{3} n_{2}^{4}+ 2360373 n_{1}^{2} n_{2}^{5} - 579415 n_{1} n_{2}^{6} + 41645 n_{2}^{7})\\
&\ \ \ +\frac{y^{-2}\, n_{1}}{ 175 \pi^{3}\left(n_{1} + n_{2}\right)^{11}} (495 n_{1}^{7} + 28765 n_{1}^{6} n_{2} - 1501506 n_{1}^{5} n_{2}^{2} + 8955786 n_{1}^{4} n_{2}^{3} \\
 & \ \ \ \   \ \ \ \ \ \ \ \ \ \ \ \ \   \ \  \ \ \  - 15910621 n_{1}^{3} n_{2}^{4} + 9855297 n_{1}^{2} n_{2}^{5} - 2066640 n_{1} n_{2}^{6} + 120280 n_{2}^{7})\\
&\ \ \ +\frac{n_{1}}{1050 \pi\left(n_{1} + n_{2}\right)^{9}}\left(175 n_{1}^{7} + 9135 n_{1}^{6} n_{2} - 348825 n_{1}^{5} n_{2}^{2} + 1980687 n_{1}^{4} n_{2}^{3} \right. \\
 & \ \ \ \   \ \ \ \ \ \ \ \  \ \ \ \ \ \ \ \   \ \ \ \ \ \ \  - 3425931 n_{1}^{3} n_{2}^{4} + 2035749 n_{1}^{2} n_{2}^{5} - 384395 n_{1} n_{2}^{6} + 15645 n_{2}^{7})\Big]	 \end{align*}\newpage
	 
	 \begin{align*}
\eta_{1,0}=
&\sgn(n_2)\Big[\frac{y^{-6}n_{1}^{2} n_{2}}{\left(n_{1} + n_{2}\right)^{15}}\frac{10296}{5 \pi^{7} }  (3 n_{1}^{5} - 35 n_{1}^{4} n_{2} + 105 n_{1}^{3} n_{2}^{2} - 105 n_{1}^{2} n_{2}^{3} + 35 n_{1} n_{2}^{4} - 3 n_{2}^{5})\\
&\ \ \ +\frac{y^{-4}\,n_{2} }{ \left(n_{1} + n_{2}\right)^{13}}\frac{22}{175 \pi^{5}}  (41645 n_{1}^{7} - 579415 n_{1}^{6} n_{2} + 2360373 n_{1}^{5} n_{2}^{2} - 3515111 n_{1}^{4} n_{2}^{3} \\
 & \ \ \ \   \ \ \ \ \ \ \ \  \ \ \ \ \ \ \ \   \ \ \ \ \ \ \
 + 1865089 n_{1}^{3} n_{2}^{4} - 283107 n_{1}^{2} n_{2}^{5} + 3965 n_{1} n_{2}^{6} + 65 n_{2}^{7})\\
& \ \ \ +\frac{y^{-2}\, n_{2} }{175 \pi^{3} (n_{1} + n_{2})^{11}} \left(120280 n_{1}^{7} - 2066640 n_{1}^{6} n_{2} + 9855297 n_{1}^{5} n_{2}^{2} - 15910621 n_{1}^{4} n_{2}^{3} \right.\\
 & \ \ \ \   \ \ \ \ \ \ \ \  \ \ \ \ \ \ \ \   \ \ \ \ \ \ \ + 8955786 n_{1}^{3} n_{2}^{4} - 1501506 n_{1}^{2} n_{2}^{5} + 28765 n_{1} n_{2}^{6} + 495 n_{2}^{7})\\
& \ \ \ +\frac{n_{2}}{1050 \pi \left(n_{1} + n_{2}\right)^{9}}(15645 n_{1}^{7} - 384395 n_{1}^{6} n_{2} + 2035749 n_{1}^{5} n_{2}^{2} - 3425931 n_{1}^{4} n_{2}^{3} \\
 & \ \ \ \   \ \ \ \ \ \ \ \  \ \ \ \ \ \ \ \   \ \ \ \ \ \ \ + 1980687 n_{1}^{3} n_{2}^{4} - 348825 n_{1}^{2} n_{2}^{5} + 9135 n_{1} n_{2}^{6} + 175 n_{2}^{7})\Big]\\
\eta_{1,1}=
&\frac{y^{-5}}{175 \pi^{6} \left(n_{1} + n_{2}\right)^{14}}(1430 n_{1}^{8} + 88660 n_{1}^{7} n_{2} - 7388524 n_{1}^{6} n_{2}^{2} + 50271364 n_{1}^{5} n_{2}^{3}\\
 & \ \ \ \   \ \ \ \ \ \ \ \  \ \ \ \ \ \ \ \   \ \ \ \ \ \ \ - 90631684 n_{1}^{4} n_{2}^{4} + 50271364 n_{1}^{3} n_{2}^{5} - 7388524 n_{1}^{2} n_{2}^{6} \\
 & \ \ \ \   \ \ \ \ \ \ \ \  \ \ \ \ \ \ \ \   \ \ \ \ \ \ \ + 88660 n_{1} n_{2}^{7} + 1430 n_{2}^{8})\\
&+\frac{y^{-3}}{175 \pi^{4} \left(n_{1} + n_{2}\right)^{12}}(1210 n_{1}^{8} + 72160 n_{1}^{7} n_{2} - 4261268 n_{1}^{6} n_{2}^{2} + 25948736 n_{1}^{5} n_{2}^{3} \\
 & \ \ \ \   \ \ \ \ \ \ \ \  \ \ \ \ \ \ \ \   \ \ \ \ \ \ \ - 45143692 n_{1}^{4} n_{2}^{4} + 25948736 n_{1}^{3} n_{2}^{5} - 4261268 n_{1}^{2} n_{2}^{6} \\
 & \ \ \ \   \ \ \ \ \ \ \ \  \ \ \ \ \ \ \ \   \ \ \ \ \ \ \ + 72160 n_{1} n_{2}^{7} + 1210 n_{2}^{8})\\
&+\frac{y^{-1}}{350 \pi^{2} \left(n_{1} + n_{2}\right)^{10}}(315 n_{1}^{8} + 17550 n_{1}^{7} n_{2} - 785712 n_{1}^{6} n_{2}^{2} + 4522866 n_{1}^{5} n_{2}^{3}\\
 & \ \ \ \   \ \ \ \ \ \ \ \  - 7798806 n_{1}^{4} n_{2}^{4} + 4522866 n_{1}^{3} n_{2}^{5} - 785712 n_{1}^{2} n_{2}^{6} + 17550 n_{1} n_{2}^{7} + 315 n_{2}^{8})\\
&+\frac{2  y n_{1} n_{2}}{525 \left(n_{1} + n_{2}\right)^{8}}(175 n_{1}^{6} - 6510 n_{1}^{5} n_{2} + 35745 n_{1}^{4} n_{2}^{2}- 61572 n_{1}^{3} n_{2}^{3} \\
 & \ \ \ \   \ \ \ \ \ \ \ \  \ \ \ \ \ \ \ \   \ \ \ \ \ \ \  + 35745 n_{1}^{2} n_{2}^{4} - 6510 n_{1} n_{2}^{5} + 175 n_{2}^{6}).
\end{align*}
The asymptotic expansion is
\begin{align}\label{eq:6.2.3fasym56}
	\hat{f}^P_{n_1,n_2}(y)=& -\frac{4576\sigma _2\left(|n_1|\right) \sigma _2\left(|n_2|\right)  }{175 \pi ^6 n_1^2 n_2^2 \left(n_1+n_2\right){}^{15}
   y^7}\left(5 n_1^9+315 n_2
   n_1^8-25524 n_2^2 n_1^7+149940 n_2^3 n_1^6\right.\nonumber\\
    & -141120 n_2^4 n_1^5-141120 n_2^5 n_1^4+149940 n_2^6
   n_1^3-25524 n_2^7 n_1^2+315 n_2^8 n_1+5 n_2^9\nonumber\\
    & +7560 n_2^2 n_1^7 \log \left(\pi  n_1\right)-7560
   n_2^2 n_1^7 \log \left(\pi  n_2\right)-88200 n_2^3 n_1^6 \log \left(\pi  n_1\right)\nonumber\\
    & +88200 n_2^3
   n_1^6 \log \left(\pi  n_2\right)+264600 n_2^4 n_1^5 \log \left(\pi  n_1\right)-264600 n_2^4 n_1^5
   \log \left(\pi  n_2\right)\\
    & -264600 n_2^5 n_1^4 \log \left(\pi  n_1\right)+264600 n_2^5 n_1^4 \log
   \left(\pi  n_2\right)+88200 n_2^6 n_1^3 \log \left(\pi  n_1\right)\nonumber\\
    & \left.-88200 n_2^6 n_1^3 \log \left(\pi
    n_2\right)-7560 n_2^7 n_1^2 \log \left(\pi  n_1\right)+7560 n_2^7 n_1^2 \log \left(\pi 
   n_2\right)\right)+O(y^{-5}).\nonumber
\end{align}
We recall that by \eqref{asymptotics_of_homogeneous_solutionnasokdijfok}, \[\hat{f}^H_{n_1,n_2}(y) = \alpha_{n_1, n_2} \sqrt{y} K_{7+\frac{1}{2}}(2 \pi  |n_1+n_2| y) = \alpha_{n_1,n_2} \frac{135135}{256 \pi ^7 |n_1+n_2|^{15/2} y^7}+o(y^{-7}).\]  
Comparing the right hand sides of the  two previous formulas, we obtain that there is a unique choice of $\alpha_{n_1, n_2}$ that guarantees that  $\hat{f}_{n_1,n_2}(y) = \hat{f}^P_{n_1,n_2}(y)+ \hat{f}^H_{n_1,n_2}(y) = o(y^{-7})$, $y \to 0$:
\begin{align}\label{alphan1n2_slodikjuf9o8w7947u2}
	\alpha_{n_1,n_2}=&  \frac{1217536 \sign(n_1+n_2) \sigma _2\left(|n_1|\right) \sigma _2\left(|n_2|\right)  }{23648625 \pi  n_1^2 n_2^2 \left|n_1+n_2\right|{}^{15/2}
   }\left(5 n_1^9+315 n_2
   n_1^8-25524 n_2^2 n_1^7+149940 n_2^3 n_1^6\right.\nonumber\\
    & -141120 n_2^4 n_1^5-141120 n_2^5 n_1^4+149940 n_2^6
   n_1^3-25524 n_2^7 n_1^2+315 n_2^8 n_1+5 n_2^9\nonumber\\
    & +7560 n_2^2 n_1^7 \log \left(\pi  n_1\right)-7560
   n_2^2 n_1^7 \log \left(\pi  n_2\right)-88200 n_2^3 n_1^6 \log \left(\pi  n_1\right)\nonumber\\
    & +88200 n_2^3
   n_1^6 \log \left(\pi  n_2\right)+264600 n_2^4 n_1^5 \log \left(\pi  n_1\right)-264600 n_2^4 n_1^5
   \log \left(\pi  n_2\right)\\
    & -264600 n_2^5 n_1^4 \log \left(\pi  n_1\right)+264600 n_2^5 n_1^4 \log
   \left(\pi  n_2\right)+88200 n_2^6 n_1^3 \log \left(\pi  n_1\right)\nonumber\\
    & \left.-88200 n_2^6 n_1^3 \log \left(\pi
    n_2\right)-7560 n_2^7 n_1^2 \log \left(\pi  n_1\right)+7560 n_2^7 n_1^2 \log \left(\pi 
   n_2\right)\right)+O\left(\frac{1}{y^5}\right).\nonumber
\end{align}

Comparing the formula above with the leading terms in the asymptotic expansion of $\hat{f}^H_{n_1,n_2}(y)$ as $y\to 0$, we get that there  exists a unique choice of $\alpha_{n_1,n_2}$ such that $\hat{f}_{n_1,n_2}(y)=\hat{f}^H_{n_1,n_2}(y)+\hat{f}^P_{n_1,n_2}(y)=o(y^{-7})$. Moreover, it follows that for fixed $n \neq 0$,
\begin{equation}\label{alpha_something_asymptwerwot6.3.3}
    \alpha_{n-n_1,n_1} = O(|n_1|^{-8}), \quad  |n_1|\to \infty.
\end{equation}

\subsubsection{ $n_1= -n_2$.}This case has been considered in \cite[(C.30)]{CGPWW2021}.
For 
$\hat{f}^P_{-n_2,n_2}(y)$  defined in \eqref{eq:6.4fdef}
we have
\begin{align*}
\mu_{0,0} =&
\frac{524288 \pi ^8 n_2^8 y^9}{7441875}-\frac{16384 \pi ^6 n^6 y^7}{2480625}-\frac{256 \pi ^4 n_2^4 y^5}{165375}-\frac{8 \pi ^2 n_2^2 y^3}{4725} +\frac{2}{175 \pi ^4 n_2^4 y^3}\\&  \ \ \ \  \ \ \ \ +\frac{2}{175 \pi ^2 n_2^2 y}+\frac{y}{210},
\\
\mu_{0,1} = &\sgn(n_2) \left( 
\frac{262144 \pi ^7 n_2^7 y^8}{7441875}+\frac{8192 \pi ^5 n_2^5 y^6}{2480625}+\frac{128 \pi ^3 n_2^3 y^4}{55125} \right. \\ &  \ \ \ \  \ \ \ \  \ \ \ \  \ \ \ \ \left. +\frac{6}{175 \pi ^3 n_2^3 y^2}+\frac{4}{175 \pi ^5 n_2^5 y^4}+\frac{4}{945} \pi  n_2 y^2+\frac{2}{105 \pi  n_2} \right),
\\   
\mu_{1,1} =&
-\frac{524288 \pi ^8 n_2^8 y^9}{7441875}-\frac{16384 \pi ^6 n_2^6 y^7}{7441875}-\frac{256 \pi ^4 n_2^4 y^5}{275625}-\frac{8 \pi ^2 n_2^2 y^3}{6615}+\frac{4}{175 \pi ^4 n_2^4 y^3} \\& \ \ \ \  \ \ \ \ +\frac{2}{175 \pi ^6 n_2^6 y^5}+\frac{3}{175 \pi ^2 n_2^2 y}-\frac{23 y}{1890} .
\end{align*}
The asymptotic expansion is 
\begin{equation}\label{eq:6.3.4fasym}
\hat{f}^P_{-n_2,n_2}(y) =-\frac{32 \sigma _2(|n_2|){}^2}{175 \pi ^6 |n_2|^{10} y^7}+O\left(\frac{1}{y^4}\right).
\end{equation}
Since $\hat{f}^H_{-n_2,n_2}(y) = \alpha_{n_1,n_2} y^{-7}$, in order to have $\hat{f}^P_{-n_2,n_2}(y)+\hat{f}^H_{-n_2,n_2}(y)=o(y^{-7})$, $y\to 0$, we must set  $\alpha_{-n_2,n_2} = \frac{32 \sigma _2(|n_2|){}^2}{175 \pi ^6 |n_2|^{10}}$. Now it is time to decide on $\alpha_{0,0}$. In order to do it, we note from \eqref{Ramanujan2},
\begin{align}\label{eq:6.3.4alphasum}
    \sum_{n_2 \neq 0} \alpha_{n_2,-n_2} = \frac{32}{175 \pi^6}\sum_{n_2 \neq 0} \frac{\sigma_2(|n_2|)^2}{n_2^{10}} = \frac{7072 \pi ^{16}}{1695787498125}, \end{align}
    and thus $\alpha_{0,0} = -\frac{7072 \pi ^{16}}{1695787498125}$.

\section{$\alpha =3/2$ and $ \beta=5/2$}\label{sec:3/2,5/2}

In this section, for $z=x+iy \in \mathfrak{H}$ we solve 
\[
(\Delta - \lambda) f(z) = - 6 \zeta(3) \zeta(5)E_{3/2}(z)E_{5/2}(z) .
\]
 for 
 \[f(z) = \sum_{n \in \Z} \sum_{n_1+n_2=n} \hat{f}_{n_1,n_2}(y) e^{2 \pi i n x}\]
 in terms of $\hat{f}_{n_1,n_2}(y) = \hat{f}^P_{n_1,n_2}(y)+\hat{f}^H_{n_1,n_2}(y).$
In order to obtain $\mathcal{E}(\lambda, 5/2, 3/2, z, \overline{z})$ from \cite[(2.13)]{CGPWW2021}, we use \eqref{our_solution_vs_Chester} to note that $\mathcal{E}(\lambda, 5/2, 3/2, z, \overline{z})$ can be obtained from $f$ solving the equation above by multiplying $f$ by $\tfrac{4}{6}=\tfrac{2}{3}$. Thus, instead of $-64$ in \eqref{eq:7.3fdef}, we need to take $-\tfrac{128}{3}$.

When $n_1=n_2=0$, $\hat{f}_{0,0}(y)$ contains no $K$-Bessel or divisor functions and is given by a a polynomial in $y$ and $1/y$ below.
 For $n_1n_1=0$ but not both zero, 
\begin{equation}\label{eq:7.1afdef}
\hat{f}^P_{0,n}(y)=- 8 \pi^2 \sigma_{-4}(|n|)|n|^2\sum_{i=0,1} \nu_{i}(n, y) K_i(2 \pi |n| y) 
\end{equation} and 
\begin{equation}\label{eq:7.1bfdef}
\hat{f}^P_{n,0}(y) =- 8 \pi \sigma_{-2}(|n|)|n|\sum_{i,j=0,1} \nu_{i}(n, y) K_i(2 \pi |n| y)\end{equation}

for $n_1 n_2\neq 0$ and $n_1+n_2\neq 0$,
\begin{align}\label{eq:7.3fdef}
\hat{f}^P_{n_1,n_2}(y) =  - 64 & \pi^3  |n_1||n_2|^2\sigma_{-2}(|n_1|) \sigma_{-4}(|n_2|) \\ &\times\sum_{i,j=0,1} \eta_{i,j}(n_1, n_2, y) K_i(2 \pi |n_1| y) K_j(2 \pi |n_2| y)\nonumber,
\end{align} and 
for $n_1 =-n_2$,
\begin{align}\label{eq:7.4fdef}
\hat{f}^P_{-n_2,n_2}(y) = - 64& \pi^3  |n_2|^3\sigma_{-2}(|n_2|) \sigma_{-4}(|n_2|) \\& \times \sum_{(i,j)\in \{(0,0),(0,1),(1,1)\} } \mu_{i,j}( n_2, y) K_i(2 \pi |n_2| y) K_j(2 \pi |n_2| y)\nonumber,
\end{align}
where $\nu_{i}, \eta_{i,j}$ and $\mu_{i,j}$ defined below depending on each value of $\lambda$.

\subsection{$\lambda =20$} This subsection corresponds to \cite[(C.19)]{CGPWW2021} with $r=4$. 

In what follows we find $\hat{f}^P_{n_1,n_2}(y)$ for different values of $n_1$ and $n_2$. 
\subsubsection{$n_1=0$ and $n_2=0$.}
The particular solution for $n_1=n_2=0$ is equal to 
\[
\hat{f}_{0,0}(y) = \frac{14175 y^6 \zeta (3) \zeta (5)+2100 \pi ^2 y^4 \zeta (5)+84 \pi ^4 y^2 \zeta (3)+40 \pi
   ^6}{18900 y^2}+\frac{c_1}{y^4}+c_2 y^5
\]
for some $c_1, c_2 \in \C$. Its asymptotic behavior for $y \to 0$ can be described by 
$
\hat{f}_{0,0}(y) = \frac{2 \pi ^6}{945 y^2}
$
and the leading term of the asymptotic behavior doesn't depend on $c_1$ and $c_2$. 
If we want to get rid of the $O(y^5)$-asymptotic, we need to set $c_2=0$. We refrain for a moment from defining $c_1$. However, we notice that the first four summands in the r.h.s. of the formula above coincide with the first four summands in the first line in \cite[(C.21)]{CGPWW2021}.

\subsubsection{$n_1=0$ and $n_2 \neq 0$} 
For $\hat{f}^P_{0,n}(y)$ as in \eqref{eq:7.1afdef}, we have 
\begin{align*}
\nu_0(n, y)&=
\frac{\pi ^2 n^2 \left(2 \zeta (2)-27 \zeta (3) y^2\right)-126 \zeta (3)}{3 \pi ^4 n^4 y^2},
\\ \nu_1(n, y)&=\sgn(n)\left(
\frac{-126 \zeta (3)+\pi ^4 n^4 y^2 \left(\zeta (2)-3 \zeta (3) y^2\right)+2 \pi ^2 n^2 \left(\zeta (2)-45 \zeta (3) y^2\right)}{3 \pi ^5 n^5 y^3}  \right).
\end{align*}
The asymptotic expansion of $\hat{f}^P_{0,n}(y)$ as $y\to 0$ is 
\begin{align*}
-\frac{8 \sigma _{-4}(|n|) \left(\pi ^2 \zeta (2) n^2-63 \zeta (3)\right)}{3 \pi ^4 n^4 y^4}+O\left(\frac{1}{y^2}\right)
\end{align*}
There is a unique choice of $\alpha_{0,n}$ that gets rid of the $y^{-4}$-term in the asymptotic expansion of $\hat{f}^P_{0,n}(y)+\hat{f}^H_{0,n}(y)$.
\subsubsection{$n_1\neq 0$ and $n_2 = 0$} 
For $\hat{f}^P_{n,0}(y)$ as in \eqref{eq:7.1bfdef}, we have 
\begin{align*}
	\nu_0(n, y)&=\sgn(n)\left(
	\frac{45 \zeta (5)}{\pi ^3 n^3}+\frac{33600 \zeta (5)+64 \pi ^4 \zeta (4) n^4}{160 \pi ^5 n^5 y^2}-\frac{\zeta (5) y^2}{2 \pi  n}  \right),
	\\ \nu_1(n, y)&=
	\frac{(128 n^{4} \pi^{4} \zeta(4)  + 67200 \zeta(5))y^{-3}}{320 n^{6} \pi^{6}}
	+\frac{150  \zeta(5)y^{-1}}{n^{4} \pi^{4}}  
	+\frac{9 \zeta(5)y}{2 n^{2} \pi^{2}} .
\end{align*}
The asymptotic expansion of $\hat{f}^P_{n,0}(y)$ as $y\to 0$ is 
\begin{align*}
- 8  \sigma_{-2}(|n|) \left( \frac{525 \zeta (5)+\pi ^4 \zeta (4) n^4}{5 \pi ^6 |n|^6 y^4}+O\left(y^{-2}\right) \right).
\end{align*}
There is a unique choice of $\alpha_{n,0}$ that gets rid of the $y^{-4}$-term in the asymptotic expansion of $\hat{f}^P_{n,0}(y)+\hat{f}^H_{n,0}(y)$.

\subsubsection{$n_1n_2\neq 0$ and $n_1+n_2\neq 0$.}

For
$\hat{f}^P_{n_1,n_2}(y)$ as in \eqref{eq:7.3fdef}, we have
\begin{align*}
\eta_{0,0}=
& \sgn(n_1)  \Big[ \frac{y^{-2}\, 7 n_{1} n_{2}^{2}}{\pi^{3} \left(n_{1} + n_{2}\right)^{8}} \left(5 n_{1}^{2} - 8 n_{1} n_{2} + 3 n_{2}^{2}\right) \\
&\ \ \ \ \ \ \ \ \ \ \ \ \ \ +\frac{n_{1} }{30 \pi \left(n_{1} + n_{2}\right)^{6}} \left(3 n_{1}^{4} + 26 n_{1}^{3} n_{2} + 348 n_{1}^{2} n_{2}^{2} - 810 n_{1} n_{2}^{3} + 145 n_{2}^{4}\right) \Big] ,\\
\eta_{0,1}=
& \sgn(n_1) \sgn(n_2) \Big[ \frac{y^{-3}\, 7 n_{1} n_{2}^{3}}{\pi^{4} \left(n_{1} + n_{2}\right)^{9}} \left(- 5 n_{1} + 3 n_{2}\right) \\
&\ \ \ \ +\frac{y^{-1}\,n_{1} }{30 n_{2} \pi^{2} \left(n_{1} + n_{2}\right)^{7}} \left(3 n_{1}^{5} + 29 n_{1}^{4} n_{2} + 149 n_{1}^{3} n_{2}^{2} + 723 n_{1}^{2} n_{2}^{3} - 1820 n_{1} n_{2}^{4} + 460 n_{2}^{5}\right)\\
& \ \ \ \ \ \ \ \ \ \ \ \ \ \ \ \ \ \ \  \ \ \ \ \ \ \ \ \ \ \ \ \ \ \ \ \ \ +\frac{y\, 2 n_{1} n_{2}}{15 \left(n_{1} + n_{2}\right)^{5}} \left(n_{1}^{3} + 15 n_{1}^{2} n_{2} - 45 n_{1} n_{2}^{2} + 5 n_{2}^{3}\right) \Big], \\
\eta_{1,0}=
&\frac{y^{-3}\, n_{1}^{2} n_{2}^{2}}{\pi^{4} \left(n_{1} + n_{2}\right)^{9}}  \left(35 n_{1} - 21 n_{2}\right)\\
&  +\frac{y^{-1}}{30 \pi^{2} \left(n_{1} + n_{2}\right)^{7}}(3 n_{1}^{5} + 29 n_{1}^{4} n_{2} + 899 n_{1}^{3} n_{2}^{2} - 1827 n_{1}^{2} n_{2}^{3} + 190 n_{1} n_{2}^{4} + 10 n_{2}^{5})\\
&  +\frac{y\, 2 n_{1} n_{2} }{15 \left(n_{1} + n_{2}\right)^{5}}\left(n_{1}^{3} + 15 n_{1}^{2} n_{2} - 45 n_{1} n_{2}^{2} + 5 n_{2}^{3}\right),\\
\eta_{1,1}=
& \sgn(n_2) \left[ \frac{y^{-2}}{30 n_{2} \pi^{3} \left(n_{1} + n_{2}\right)^{8}}(3 n_{1}^{6} + 32 n_{1}^{5} n_{2} + 178 n_{1}^{4} n_{2}^{2} + 872 n_{1}^{3} n_{2}^{3}\right.\\ &  \ \ \ \ \ \ \ \ \ \ \ \ \ \  \  \ \ \ \ \ \ \ \ \ \ \  \  \ \ \ \ \ \ \ \ \ \ \  \  \ \ \ \ \ \ \   - 2447 n_{1}^{2} n_{2}^{4} + 200 n_{1} n_{2}^{5} + 10 n_{2}^{6}) \\
&\ \ \ \ \ \ \ \ \ \ \ \ \ \left. +\frac{1}{30 \pi \left(n_{1} + n_{2}\right)^{6}}(4 n_{1}^{5} + 43 n_{1}^{4} n_{2} + 334 n_{1}^{3} n_{2}^{2} - 880 n_{1}^{2} n_{2}^{3} + 110 n_{1} n_{2}^{4} + 5 n_{2}^{5}) \right].
\end{align*}

We note that 
\begin{align}\label{eq:7.2.4fasym}
	\hat{f}^P_{n_1,n_2}(y) & = - \sgn(n_1)y^{-4} \frac{8   \sigma_{-2}(|n_1|) \sigma_{-4}(|n_2|) }{15  \pi^{2} \left(n_{1} + n_{2}\right)^{9}
} (3 n_{1}^{7} + 35 n_{1}^{6} n_{2} + 210 n_{1}^{5} n_{2}^{2} + 1050 n_{1}^{4} n_{2}^{3}\nonumber \\ &+ 2100 n_{1}^{3} n_{2}^{4} \log{\left(|n_{1}| \pi \right)} - 2100 n_{1}^{3} n_{2}^{4} \log{\left(|n_{2}| \pi \right)} - 1575 n_{1}^{3} n_{2}^{4} - 1260 n_{1}^{2} n_{2}^{5} \log{\left(|n_{1}| \pi \right)} \\ & + 1260 n_{1}^{2} n_{2}^{5} \log{\left(|n_{2}| \pi \right)} - 2247 n_{1}^{2} n_{2}^{5} + 210 n_{1} n_{2}^{6} + 10 n_{2}^{7})+o(y^{-4}).\nonumber
\end{align}
Requiring $\hat{f}_{n_1,n_2}(y)=o(y^{-4})$ gives us a unique $\alpha_{n_1, n_2}$ that cancels with the $O(y^{-4})$-term.
Comparing the formula above with the leading terms in the asymptotic expansion of $\hat{f}^H_{n_1,n_2}(y)$ as $y\to 0$, we get that there  exists a unique choice of $\alpha_{n_1,n_2}$ such that $\hat{f}_{n_1,n_2}(y)=\hat{f}^H_{n_1,n_2}(y)+\hat{f}^P_{n_1,n_2}(y)=o(y^{-4})$. Moreover, it follows that for fixed $n \neq 0$,
\begin{equation}\label{alpha_something_asymptot7.7}
    \alpha_{n-n_1,n_1} = O( |n_1|^{-5}), \quad  |n_1|\to \infty.
\end{equation}

\subsubsection{ $n_1= -n_2$.}

For $
\hat{f}^P_{-n_2,n_2}(y)$ as in \eqref{eq:7.4fdef}, we have
\begin{align*}
	\mu_{0,0}&= \sgn(n_2) \left(  
	\frac{1 }{36 n_{2} \pi}  
	+\frac{2 n_{2} \pi y^{2}}{315}  
	+\frac{128 n_{2}^{3} \pi^{3} y^{4}}{4725}  
	-\frac{ 4096 n_{2}^{5} \pi^{5} y^{6}}{14175}  \right)
	\\
	\mu_{0,1}&=
	\frac{ y^{-1}}{18 n_{2}^{2} \pi^{2}}  
	-\frac{ y}{105}  
	-\frac{ 64 n_{2}^{2} \pi^{2} y^{3}}{4725}  
	-\frac{ 2048 n_{2}^{4} \pi^{4} y^{5}}{14175}  
	\\
	\mu_{1,1}&= \sgn(n_2) \left( 
	\frac{ y^{-2}}{36 n_{2}^{3} \pi^{3}}  
	-\frac{1 }{63 n_{2} \pi}  
	+\frac{2 n_{2} \pi y^{2}}{525}  
	+\frac{128 n_{2}^{3} \pi^{3} y^{4}}{14175}  
	+\frac{4096 n_{2}^{5} \pi^{5} y^{6}}{14175} \right).
\end{align*}
The asymptotic expansion is 
\begin{align}\label{eq:7.2.5fasym}
\hat{f}^P_{-n_2,n_2}(y) = - 64& y \sigma_{-2}(|n_2|) \sigma_{-4}(|n_2|)  \left( \frac{1}{144 \pi ^2 |n_2|^2 y^4}-\frac{1}{56y^2}+O\left(1\right) \right).
\end{align}
There is a unique choice of $\alpha_{-n_2,n_2}$ such that $\hat{f}^P_{-n_2,n_2}(y)+\hat{f}^H_{-n_2,n_2}(y)=o(y^{-4})$:
\begin{equation}\label{eq:7.9convergence}
    \alpha_{-n_2,n_2}=  \frac{4 y \sigma_{-2}(|n_2|) \sigma_{-4}(|n_2|)}{9 \pi ^2 |n_2|^2}.
\end{equation}

\section{$\alpha = \beta=5/2$}\label{sec:5/2}
In this section we solve 
\[
(\Delta - \lambda) f(z) = - \left(3 \zeta(5)E_{5/2}(z)\right)^2, \quad z=x+iy \in \mathfrak{H}
\]
 for 
 \[f(z) = \sum_{n \in \Z} \sum_{n_1+n_2=n} \hat{f}_{n_1,n_2}(y) e^{2 \pi i n x}\]
 in terms of $\hat{f}_{n_1,n_2}(y) = \hat{f}^P_{n_1,n_2}(y)+\hat{f}^H_{n_1,n_2}(y).$

When $n_1=n_2=0$, $\hat{f}_{0,0}(y)$ contains no $K$-Bessel or divisor functions and is given by a a polynomial in $y$ and $1/y$ below. For $n_1n_1=0$ but not both zero, 
\begin{equation}\label{eq:8.1fdef}
\hat{f}^P_{0,n}(y) =\hat{f}^P_{n,0}(y)= -8\pi^2 \sigma_{-4}(|n|) |n|^2\sum_{i=0 } \nu_{i}(n, y) K_i(2 \pi |n| y),
\end{equation}
for $n_1 n_2\neq 0$ and $n_1+n_2\neq 0$,
\begin{equation}\label{eq:8.3fdef}
\hat{f}^P_{n_1,n_2}(y)= - 64 \pi^4  |n_1|^2|n_2|^2\sigma_{-4}(|n_1|) \sigma_{-4}(|n_2|) \sum_{i,j=0,1}  \eta_{i,j}(n_1, n_2, y) K_i(2 \pi |n_1| y) K_j(2 \pi | n_2 | y) ,
\end{equation} and 
for $n_1 =-n_2$,
\begin{align}\label{eq:8.4fdef}
\hat{f}^P_{-n_2,n_2}(y) = - 64 &\pi^4  |n_1|^2|n_2|^2\sigma_{-4}(|n_2|) \sigma_{-4}(|n_2|)\\
& \times\sum_{(i,j) \in \{(0,0),(0,1),(1,1)\}}  \eta_{i,j}(n_2, y) K_i(2 \pi |n_2| y) K_j(2 \pi |n_2| y)\nonumber
\end{align}
where $\nu_{i}, \eta_{i,j}$ and $\mu_{i,j}$ are defined below depending on each value of $\lambda$.

\subsection{$\lambda = 30$}

This case corresponds to \cite[C.3.2]{CGPWW2021} with $r=5$.

\subsubsection{$n_1=0$ and $n_2=0$.}

It is not complicated to show that 
\[
\hat{f}_{0,0}(y) = \frac{80 \zeta (4)^2+81 \zeta (5)^2 y^8+72 \zeta (4) \zeta (5) y^4}{90 y^3}+c_2 y^6+\frac{c_1}{y^5}
\]
for some $c_1, c_2 \in \C$.

\subsubsection{$n_1 n_2=0$ but not both zero} 
For $\hat{f}^P_{0,n}(y) $ as in \eqref{eq:8.1fdef}, we have
\begin{align*}
\nu_0&=
\frac{4\zeta(4)\pi ^4 n^4 +7560\zeta(5)}{5 \pi ^6 n^6 y^3} 
+\frac{420\zeta(5)y^{-1}}{n^{4} \pi^{4}}  
+\frac{6\zeta(5)y}{n^{2} \pi^{2}}  
\\ \nu_1&=\sgn(n)\left(
\frac{4\zeta(4)\pi^4 n^4 +7560\zeta(5)}{5 \pi ^7 n^7 y^4} 
+\frac{4\zeta(4)\pi ^4 n^4 +11760\zeta(5)}{10 \pi ^5 n^5 y^2}
+\frac{90\zeta(5)}{n^{3} \pi^{3}}  
-\frac{ \zeta(5)y^{2}}{2 n \pi} \right).
\end{align*}

Its asymptotic behavior as $y \to 0$ can be described as 
\begin{equation}
-8  \sigma_{-4}(|n|)\left( \frac{2 \left(1890 \zeta (5)+\pi ^4 \zeta (4) n^4\right)}{5 \pi ^6 n^8 y^5}-\frac{840 \zeta (5)+\pi ^4 \zeta (4) n^2}{5 \pi ^4 n^4 y^3}+O\left(\frac{1}{y^2}\right) \right).
\end{equation}

\subsubsection{$n_1 n_2 \neq 0$ and $n_1+n_2\neq 0$}
For $\hat{f}^P_{n_1,n_2}(y) $ as in \eqref{eq:8.3fdef}
we have 
\begin{align*}
	\eta_{0,0}=
	&\frac{y^{-3}}{\pi^{4} \left(n_{1} + n_{2}\right)^{10}}(252 n_{1}^{2} n_{2}^{2} \left(n_{1}^{2} - 2 n_{1} n_{2} + n_{2}^{2}\right))\\
	&+\frac{y^{-1}}{5 \pi^{2} \left(n_{1} + n_{2}\right)^{8}}(n_{1}^{6} + 12 n_{1}^{5} n_{2} + 431 n_{1}^{4} n_{2}^{2} - 1400 n_{1}^{3} n_{2}^{3} + 431 n_{1}^{2} n_{2}^{4} + 12 n_{1} n_{2}^{5} + n_{2}^{6})\\
	&+\frac{y}{15 \left(n_{1} + n_{2}\right)^{6}}(2 n_{1} n_{2} \left(n_{1}^{4} + 20 n_{1}^{3} n_{2} - 90 n_{1}^{2} n_{2}^{2} + 20 n_{1} n_{2}^{3} + n_{2}^{4}\right))\\
	\eta_{0,1}=
	& \sgn(n_2) \left(\frac{y^{-4}}{\pi^{5} \left(n_{1} + n_{2}\right)^{11}}(252 n_{1}^{2} n_{2}^{3} \left(- n_{1} + n_{2}\right)) \right. \\
	&+\frac{y^{-2}}{5 n_{2} \pi^{3} \left(n_{1} + n_{2}\right)^{9}}(n_{1}^{7} + 13 n_{1}^{6} n_{2} + 93 n_{1}^{5} n_{2}^{2} + 641 n_{1}^{4} n_{2}^{3} - 2859 n_{1}^{3} n_{2}^{4} + 1073 n_{1}^{2} n_{2}^{5}\\ & \ \ \ \ \ \ \ \ \ \ \ \ \ \ \ \ \ \ \ \ \ \ \ \ \ \ \ \ \ \ \ \ \ \ \ \ \ \ \ \ \ \ \ \ \ \  \ \ \ \ \ \  + 13 n_{1} n_{2}^{6} + n_{2}^{7})\\
	& \left. +\frac{1}{30 \pi \left(n_{1} + n_{2}\right)^{7}}(4 n_{1}^{6} + 57 n_{1}^{5} n_{2} + 623 n_{1}^{4} n_{2}^{2} - 2590 n_{1}^{3} n_{2}^{3} + 714 n_{1}^{2} n_{2}^{4} + 37 n_{1} n_{2}^{5} + 3 n_{2}^{6}) \right)\\
	\eta_{1,0}=
	& \sgn(n_1) \left( \frac{y^{-4}}{\pi^{5} \left(n_{1} + n_{2}\right)^{11}}(252 n_{1}^{3} n_{2}^{2} \left(n_{1} - n_{2}\right)) \right. \\
	&+\frac{y^{-2}}{5 n_{1} \pi^{3} \left(n_{1} + n_{2}\right)^{9}}(n_{1}^{7} + 13 n_{1}^{6} n_{2} + 1073 n_{1}^{5} n_{2}^{2} - 2859 n_{1}^{4} n_{2}^{3} + 641 n_{1}^{3} n_{2}^{4} + 93 n_{1}^{2} n_{2}^{5} \\ & \ \ \ \ \ \ \ \ \ \ \ \ \ \ \ \ \ \ \ \ \ \ \ \ \ \ \ \ \ \ \ \ \ \ \ \ \ \ \ \ \ \ \ \ \ \  \ \ \ \ \ \ + 13 n_{1} n_{2}^{6} + n_{2}^{7})\\
	& \left. +\frac{1}{30 \pi \left(n_{1} + n_{2}\right)^{7}}(3 n_{1}^{6} + 37 n_{1}^{5} n_{2} + 714 n_{1}^{4} n_{2}^{2} - 2590 n_{1}^{3} n_{2}^{3} + 623 n_{1}^{2} n_{2}^{4} + 57 n_{1} n_{2}^{5} + 4 n_{2}^{6}) \right)\end{align*}\newpage
	 
	 \begin{align*}
	\eta_{1,1}=
	& \sgn(n_1n_2) \left( \frac{y^{-3}}{5 n_{1} n_{2} \pi^{4} \left(n_{1} + n_{2}\right)^{10}}(n_{1}^{8} + 14 n_{1}^{7} n_{2} + 106 n_{1}^{6} n_{2}^{2} + 734 n_{1}^{5} n_{2}^{3} - 3758 n_{1}^{4} n_{2}^{4}\right.\\ & \ \ \ \ \ \ \ \ \ \ \ \ \ \ \ \ \ \ \ \ \ \ \ \ \ \ \ \ \ \ \ \ \ \ \ \ \ \ \ \ \ \ \ \ \ \  \ \ \ \ \ \ + 734 n_{1}^{3} n_{2}^{5} + 106 n_{1}^{2} n_{2}^{6} + 14 n_{1} n_{2}^{7} + n_{2}^{8}) \\
	&+\frac{y^{-1}}{30 n_{1} n_{2} \pi^{2} \left(n_{1} + n_{2}\right)^{8}}(3 n_{1}^{8} + 40 n_{1}^{7} n_{2} + 304 n_{1}^{6} n_{2}^{2} + 2360 n_{1}^{5} n_{2}^{3} - 9254 n_{1}^{4} n_{2}^{4} \\ & \ \ \ \ \ \ \ \ \ \ \ \ \ \ \ \ \ \ \ \ \ \ \ \ \ \ \ \ \ \ \ \ \ \ \ \ \ \ \ \ \ \ \ \ \ \  \ \ \ \ \ \ + 2360 n_{1}^{3} n_{2}^{5} + 304 n_{1}^{2} n_{2}^{6} + 40 n_{1} n_{2}^{7} + 3 n_{2}^{8})\\
	& \left. +\frac{y}{15 \left(n_{1} + n_{2}\right)^{6}}(2 n_{1} n_{2} \left(n_{1}^{4} + 20 n_{1}^{3} n_{2} - 90 n_{1}^{2} n_{2}^{2} + 20 n_{1} n_{2}^{3} + n_{2}^{4}\right)) \right).
\end{align*}
We note that 
\begin{align}\label{eq:8.3.3fasym}
\hat{f}^P_{n_1,n_2}(y) =-  \frac{y^{-5} 16\sigma_{-4}(|n_1|) \sigma_{-4}(|n_2|)}{5  \pi^{2} \left(n_{1} + n_{2}\right)^{11}} & (n_{1}^{9} + 15 n_{1}^{8} n_{2} + 120 n_{1}^{7} n_{2}^{2} + 840 n_{1}^{6} n_{2}^{3} 
\\ & + 2520 n_{1}^{4} n_{2}^{4} (n_1-n_2) \log{\left(|n_{1} / n_2| \right)} \nonumber
\\
&- 3024 n_{1}^{5} n_{2}^{4} - 3024 n_{1}^{4} n_{2}^{5} + 840 n_{1}^{3} n_{2}^{6} \nonumber\\&  + 120 n_{1}^{2} n_{2}^{7} + 15 n_{1} n_{2}^{8} + n_{2}^{9})+o(y^{-5}).\nonumber
\end{align}
There exists a choice of $\alpha_{n_1,n_2}$ such that $\hat{f}_{n_1,n_2}(y)=o(y^{-5})$.
Comparing the formula above with the leading terms in the asymptotic expansion of $\hat{f}^H_{n_1,n_2}(y)$ as $y\to 0$, we get that there  exists a unique choice of $\alpha_{n_1,n_2}$ such that $\hat{f}_{n_1,n_2}(y)=\hat{f}^H_{n_1,n_2}(y)+\hat{f}^P_{n_1,n_2}(y)=o(y^{-5})$. Moreover, it follows that for fixed $n \neq 0$,
\begin{equation}\label{alpha_something_asymptot6.3.3}
    \alpha_{n-n_1,n_1} = O(|n_1|^{-6}), \quad  |n_1|\to \infty.
\end{equation}

\subsubsection{ $n_1= -n_2$.}

For $\hat{f}^P_{-n_2,n_2}(y)$ as in \eqref{eq:8.4fdef}, we have 
\begin{align*}
\mu_{0,0}&=
\frac{2 y^{-1}}{55 n_{2}^{2} \pi^{2}}  
-\frac{19 y}{990}  
-\frac{ 8 n_{2}^{2} \pi^{2} y^{3}}{3465}  
-\frac{ 512 n_{2}^{4} \pi^{4} y^{5}}{51975}  
+\frac{16384 n_{2}^{6} \pi^{6} y^{7}}{155925}  ,
\\
\mu_{0,1}&= \sgn(n_2) \left(  
\frac{4 y^{-2}}{55 n_{2}^{3} \pi^{3}}  
-\frac{1}{495 n_{2} \pi}  
+\frac{4 n_{2} \pi y^{2}}{1155}  
+\frac{256 n_{2}^{3} \pi^{3} y^{4}}{51975}  
+\frac{8192 n_{2}^{5} \pi^{5} y^{6}}{155925}  \right),
\\
\mu_{1,1}&=
\frac{2 y^{-3}}{55 n_{2}^{4} \pi^{4}}  
+\frac{17 y^{-1}}{990 n_{2}^{2} \pi^{2}}  
+\frac{103 y}{6930}  
-\frac{ 8 n_{2}^{2} \pi^{2} y^{3}}{5775}  
-\frac{ 512 n_{2}^{4} \pi^{4} y^{5}}{155925}  
-\frac{ 16384 n_{2}^{6} \pi^{6} y^{7}}{155925}. 
\end{align*} 
Its asymptotic behavior for $y \to 0$ is 
\begin{equation}\label{eq:8.3.4fasym}
-64 \sigma_{-4}(|n_2|)^2 \left( \frac{1}{110 \pi ^2 n_2^2 y^5}-\frac{1}{72 y^3}+O\left(\frac{1}{y^2}\right) \right).
\end{equation}
From this we obtain that 
\begin{equation}\label{eq:8.11convergence}
    \alpha_{-n_2,n_2} = O(|n_2|^{-6}), \quad |n_2| \to \infty.
\end{equation}

\section{$\alpha =3/2$ and $\beta=7/2$}\label{sec:3/2,7/2}

In this section we solve 
\[
(\Delta - \lambda) f(z) = - 30 \zeta(3) \zeta(7)E_{3/2}(z)E_{7/2}(z), \quad z=x+iy \in \mathfrak{H}
\]
 for 
 \[f(z) = \sum_{n \in \Z} \sum_{n_1+n_2=n} \hat{f}_{n_1,n_2}(y) e^{2 \pi i n x}\]
 in terms of $\hat{f}_{n_1,n_2}(y) = \hat{f}^P_{n_1,n_2}(y)+\hat{f}^H_{n_1,n_2}(y).$
  
  For $n_1n_1=0$ but not both zero, 
\begin{equation}\label{eq:9.1afdef}
\hat{f}^P_{0,n}(y) =-16\pi^3 \sigma_{-6}(|n|) |n|^3\sum_{i=0,1 } \nu_{i}(n, y) K_i(2 \pi |n| y) ,
\end{equation} and
\begin{equation}\label{eq:9.1bfdef}
 \hat{f}^P_{n,0}(y) = -8\pi \sigma_{-2}(|n|) |n|\sum_{i=0,1 } \nu_{i}(n, y) K_i(2 \pi |n| y) ,
\end{equation}
for $n_1 n_2\neq 0$ and $n_1+n_2\neq 0$,
\begin{align}\label{eq:9.3fdef}
\hat{f}^P_{n_1,n_2}(y) = - 128 &\pi^4  |n_1||n_2|^3\sigma_{-2}(|n_1|) \sigma_{-6}(|n_2|)\nonumber\\ & \times \sum_{(i,j) \in \{(0,0),(0,1),(1,1)\}}  \eta_{i,j}(n_1, n_2, y) K_i(2 \pi |n_1| y) K_j(2 \pi |n_2| y),
\end{align}
and 
for $n_1 =-n_2$,
\begin{align}\label{eq:9.4fdef}
\hat{f}^P_{-n_2,n_2}(y) = - 128 &\pi^4  |n_2|^4\sigma_{-2}(|n_2|) \sigma_{-6}(|n_2|)\nonumber\\ & \times \sum_{(i,j) \in \{(0,0),(0,1),(1,1)\}}  \mu_{i,j}(n_2, y) K_i(2 \pi |n_2| y) K_j(2 \pi |n_2| y)
\end{align}
where $\nu_{i}, \eta_{i,j}$ and $\mu_{i,j}$ defined below depending on each value of $\lambda$.
 
\subsection{$\lambda =30$} This case corresponds to \cite[C.3.3]{CGPWW2021} with $r=5$.

\subsubsection{$n_1=0$ and $n_2=0$.}

It is not complicated to show that 
\[
\hat{f}_{0,0}(y) =\frac{2143260 y^8 \zeta (3) \zeta (7)+297675 \pi ^2 y^6 \zeta (7)+864 \pi ^6 y^2 \zeta (3)+448
   \pi ^8}{714420 y^3}+\frac{c_1}{y^5}+c_2 y^6
\]
for some $c_1, c_2 \in \C$. Its asymptotic behavior for $y \to 0$ can be described by 
$
\hat{f}_{0,0}(y) = \frac{16 \pi ^8}{25515 y^3}
$
and the leading term of the asymptotic behavior doesn't depend on $c_1$ and $c_2$.

\subsubsection{$n_1=0$ and $n_2 \neq 0$}  For $n=n_2\neq0 $, 
\[
\hat{f}^P_{0,n}(y) =-16\pi^3 \sigma_{-6}(|n|) |n|^3\sum_{i=0,1 } \nu_{i}(n, y) K_i(2 \pi |n| y) ,
\]
where 
\begin{align*}
 \nu_0(n, y)&=\sgn(n)\left(4\zeta(2) \left(\frac{1}{\left(\pi ^3 n^3\right) y^3}+\frac{1}{(6 \pi  n) y}\right)-2\zeta(3)\left(\frac{108}{\left(\pi ^5 n^5\right) y^3}+\frac{30}{\left(\pi ^3 n^3\right) y}+\frac{y}{2 \pi  n}\right)\right)
\\ \nu_1(n, y)&=4\zeta(2)  \left(\frac{1}{\left(\pi ^4 n^4\right) y^4}+\frac{2}{\left(3 \pi ^2 n^2\right) y^2}\right)-2\zeta(3) \left(\frac{84}{\left(\pi ^4 n^4\right) y^2}+\frac{108}{\left(\pi ^6 n^6\right) y^4}+\frac{13}{2 \pi ^2 n^2}\right).
\end{align*}
Its asymptotic behavior is 
\begin{align*}
    -16 \pi^3 \sigma_{-6}(|n|) |n|^3 \left( \frac{2 \left(\pi ^2 \zeta (2) n^2-54 \zeta (3)\right)}{\pi ^7 n^7 y^5}-\frac{2 \left(\pi ^2 \zeta (2) n^2-36 \zeta (3)\right)}{3 \left(\pi ^5 n^5\right) y^3}+O\left(\frac{1}{y^2}\right) \right).
\end{align*}
\subsubsection{$n_1\neq 0$ and $n_2 = 0$}
For $n=n_1\neq 0$,
\[
\hat{f}^P_{n,0}(y)=-8\pi \sigma_{-2}(|n|) |n|\sum_{i=0,1 } \nu_{i}(n, y) K_i(2 \pi |n| y) 
\]
with 
\begin{align*}
\nu_0&=\sgn(n)\left(\frac{8\pi\zeta(6)y^{-3}}{7 n }  -15\zeta(7) \left(\frac{6 y}{\pi ^3 n^3}+\frac{1512}{\left(\pi ^7 n^7\right) y^3}+\frac{420}{\left(\pi ^5 n^5\right) y}+\frac{y^3}{10 \pi  n}\right)\right)
\\ \nu_1&=
\frac{8\zeta(6)y^{-4}}{7 n^{2} }  -15\zeta(7) \left(-\frac{2 y^2}{5 \pi ^2 n^2}+\frac{1176}{\left(\pi ^6 n^6\right) y^2}+\frac{1512}{\left(\pi ^8 n^8\right) y^4}+\frac{90}{\pi ^4 n^4}\right).
\end{align*}
Its asymptotic behavior is 
\[
-8  \sigma_{-2} (|n|) |n|\left( \frac{4 \left(\pi ^6 \zeta (6) n^6-19845 \zeta (7)\right)}{7 \pi ^8 n^9 y^5}-\frac{4 \left(\pi ^6 \zeta (6) n^6-4410 \zeta (7)\right)}{7 \left(\pi ^6 n^7\right) y^3}+O\left(\frac{1}{y^2}\right)  \right).
\]

\subsubsection{$n_1 n_2 \neq 0$ and $n_1+n_2 \neq 0$} 
For $\hat{f}^P_{n_1,n_2}(y)$ as in \eqref{eq:9.3fdef}, we have 
\begin{align*}
\eta_{0,0}=
& \sgn(n_1n_2) \left( \frac{y^{-3}}{\pi^{4} \left(n_{1} + n_{2}\right)^{10}}(36 n_{1} n_{2}^{3}  \left(7 n_{1}^{2} - 10 n_{1} n_{2} + 3 n_{2}^{2}\right))\right. \\
&+\frac{n_{1}\, y^{-1}}{35 n_{2} \pi^{2} \left(n_{1} + n_{2}\right)^{8}} \left(5 n_{1}^{6} + 48 n_{1}^{5} n_{2} + 219 n_{1}^{4} n_{2}^{2} + 664 n_{1}^{3} n_{2}^{3} + 4163 n_{1}^{2} n_{2}^{4}\right.  \\ & \ \ \ \ \ \ \ \ \ \ \ \ \ \ \ \ \ \ \ \ \ \ \ \ \ \ \ \ \ \ \ \ \ \ \ \ \ \ \ \ \ \ \ \ \ \  \ \ \ \ \ \ \left.- 6440 n_{1} n_{2}^{5} + 1085 n_{2}^{6}\right)\\
& \left. +\frac{2 n_{1} n_{2} \,y}{105 \left(n_{1} + n_{2}\right)^{6}} \left(3 n_{1}^{4} + 28 n_{1}^{3} n_{2} + 210 n_{1}^{2} n_{2}^{2} - 420 n_{1} n_{2}^{3} + 35 n_{2}^{4}\right) \right)\end{align*}\newpage
	 
	 \begin{align*}
\eta_{0,1}=
& \sgn(n_1) \left( \frac{y^{-4}}{\pi^{5} \left(n_{1} + n_{2}\right)^{11}}(36 n_{1} n_{2}^{4} \left(- 7 n_{1} + 3 n_{2}\right)) \right. \\
&+\frac{n_{1} \, y^{-2}}{35 n_{2}^{2} \pi^{3} \left(n_{1} + n_{2}\right)^{9}}\left(5 n_{1}^{7} + 53 n_{1}^{6} n_{2} + 267 n_{1}^{5} n_{2}^{2} + 883 n_{1}^{4} n_{2}^{3} + 2377 n_{1}^{3} n_{2}^{4} + 7593 n_{1}^{2} n_{2}^{5}\right.\\ & \ \ \ \ \ \ \ \ \ \ \ \ \ \ \ \ \ \ \ \ \ \ \ \ \ \ \ \ \ \ \ \ \ \ \ \ \ \ \ \ \ \ \ \ \ \  \ \ \ \ \ \ \left. - 13545 n_{1} n_{2}^{6} + 2975 n_{2}^{7}\right)\\
& \left. +\frac{n_{1} }{210 \pi \left(n_{1} + n_{2}\right)^{7}}  \left(27 n_{1}^{5} + 253 n_{1}^{4} n_{2} + 1270 n_{1}^{3} n_{2}^{2} + 6594 n_{1}^{2} n_{2}^{3} - 12145 n_{1} n_{2}^{4} + 1505 n_{2}^{5}\right) \right)\\
\eta_{1,0}= 
& \sgn(n_2) \left(  \frac{y^{-4}}{\pi^{5} \left(n_{1} + n_{2}\right)^{11}}(n_{1}^{2} n_{2}^{3} \cdot \left(252 n_{1} - 108 n_{2}\right)) \right. \\
&+\frac{y^{-2}}{35 n_{2} \pi^{3} \left(n_{1} + n_{2}\right)^{9}}(5 n_{1}^{7} + 53 n_{1}^{6} n_{2} + 267 n_{1}^{5} n_{2}^{2} + 883 n_{1}^{4} n_{2}^{3} + 9237 n_{1}^{3} n_{2}^{4} - 12987 n_{1}^{2} n_{2}^{5} \\ & \ \ \ \ \ \ \ \ \ \ \ \ \ \ \ \ \ \ \ \ \ \ \ \ \ \ \ \  \ \ \ \ \ \ + 875 n_{1} n_{2}^{6} + 35 n_{2}^{7})\\
& \left. +\frac{1}{210 \pi \left(n_{1} + n_{2}\right)^{7}}(24 n_{1}^{6} + 225 n_{1}^{5} n_{2} + 1063 n_{1}^{4} n_{2}^{2} + 7042 n_{1}^{3} n_{2}^{3} - 11970 n_{1}^{2} n_{2}^{4} \right.\\ & \ \ \ \ \ \ \ \ \ \ \ \ \ \ \ \ \ \ \ \ \ \ \ \ \ \ \ \  \ \ \ \ \ \ \left.+ 1085 n_{1} n_{2}^{5} + 35 n_{2}^{6}) \right)\\
\eta_{1,1}=
&\frac{y^{-3}}{35 n_{2}^{2} \pi^{4} \left(n_{1} + n_{2}\right)^{10}}(5 n_{1}^{8} + 58 n_{1}^{7} n_{2} + 320 n_{1}^{6} n_{2}^{2} + 1150 n_{1}^{5} n_{2}^{3} + 3260 n_{1}^{4} n_{2}^{4} + 9970 n_{1}^{3} n_{2}^{5}\\ 
& \ \ \ \ \ \ \ \ \ \ \ \ \ \ \ \ \ \ \ \ \ \ \ \ \ \ \ \  \ \ \ \ \ \ 
- 16732 n_{1}^{2} n_{2}^{6} + 910 n_{1} n_{2}^{7} + 35 n_{2}^{8})\\
&+\frac{y^{-1}}{105 n_{2} \pi^{2} \left(n_{1} + n_{2}\right)^{8}}(12 n_{1}^{7} + 132 n_{1}^{6} n_{2} + 716 n_{1}^{5} n_{2}^{2} + 2806 n_{1}^{4} n_{2}^{3} + 12332 n_{1}^{3} n_{2}^{4}
\\ & \ \ \ \ \ \ \ \ \ \ \ \ \ \ \ \ \ \ \ \ \ \ \ \ \ \ \ \  \ \ \ \ \ \ - 21728 n_{1}^{2} n_{2}^{5} + 1820 n_{1} n_{2}^{6} + 70 n_{2}^{7})\\
&+\frac{y}{105 \left(n_{1} + n_{2}\right)^{6}}(2 n_{1} n_{2}  \left(3 n_{1}^{4} + 28 n_{1}^{3} n_{2} + 210 n_{1}^{2} n_{2}^{2} - 420 n_{1} n_{2}^{3} + 35 n_{2}^{4}\right))
\end{align*}
The leading term of the asymptotic expansion is 
\begin{align}\label{eq:9.2.4fasym}
\hat{f}^P_{n_1,n_2}(y) = & -\sgn(n_1n_2) y^{-5} \frac{32 \pi^4  \sigma_{-2}(|n_1|) \sigma_{-6}(|n_2|)}{35  \pi^{6} \left(n_{1} + n_{2}\right)^{11}} (5 n_{1}^{9} + 63 n_{1}^{8} n_{2} \nonumber\\ &+ 378 n_{1}^{7} n_{2}^{2} + 1470 n_{1}^{6} n_{2}^{3} + 4410 n_{1}^{5} n_{2}^{4} + 13230 n_{1}^{4} n_{2}^{5} \\ & + 17640 n_{1}^{3} n_{2}^{6} \log{\left(|n_{1}/n_{2}| \right)} - 6762 n_{1}^{3} n_{2}^{6}\nonumber  - 7560 n_{1}^{2} n_{2}^{7} \log{\left(|n_{1}/ n_{2}|  \right)}\nonumber \\ & - 15822 n_{1}^{2} n_{2}^{7} + 945 n_{1} n_{2}^{8} + 35 n_{2}^{9})  + o(y^{-5}) .\nonumber\end{align}

Comparing the formula above with the leading terms in the asymptotic expansion of $\hat{f}^H_{n_1,n_2}(y)$ as $y\to 0$, we get that there  exists a unique choice of $\alpha_{n_1,n_2}$ such that $\hat{f}_{n_1,n_2}(y)=  \hat{f}^H_{n_1,n_2}(y)+\hat{f}^P_{n_1,n_2}(y)=o(y^{-5})$. Moreover, it follows that for fixed $n \neq 0$,
\begin{equation}\label{alpha_something_asymptot6.3.sdfsdfsdq3}
    \alpha_{n-n_1,n_1} = O(|n_1|^{-6}), \quad  |n_1|\to \infty.
\end{equation}

\subsubsection{$n_1 = -n_2$}
For $\hat{f}^P_{-n_2,n_2}(y) $ as in \eqref{eq:9.4fdef}, we have 
\begin{align*}
\mu_{0,0}&=
\frac{3 y^{-1}}{55 n_{2}^{2} \pi^{2}}  
+\frac{103 y}{6930}  
-\frac{ 8 n_{2}^{2} \pi^{2} y^{3}}{4851}  
-\frac{ 512 n_{2}^{4} \pi^{4} y^{5}}{72765}  
+\frac{16384 n_{2}^{6} \pi^{6} y^{7}}{218295}  ,
\\
\mu_{0,1}&= \sgn(n_2) \left(  
\frac{6 y^{-2}}{55 n_{2}^{3} \pi^{3}}  
+\frac{89 }{6930 n_{2} \pi}  
+\frac{4 n_{2} \pi y^{2}}{1617}  
+\frac{256 n_{2}^{3} \pi^{3} y^{4}}{72765}  
+\frac{8192 n_{2}^{5} \pi^{5} y^{6}}{218295}  \right),
\\
\mu_{1,1}&=
\frac{3 y^{-3}}{55 n_{2}^{4} \pi^{4}}  
-\frac{ y^{-1}}{495 n_{2}^{2} \pi^{2}}  
-\frac{871 y}{48510}  
-\frac{ 8 n_{2}^{2} \pi^{2} y^{3}}{8085}  
-\frac{ 512 n_{2}^{4} \pi^{4} y^{5}}{218295}  
-\frac{16384 n_{2}^{6} \pi^{6} y^{7}}{218295}.  
\end{align*} 

Its asymptotic behavior is 
\begin{equation}\label{eq:9.2.5fsaym}
-128  \sigma_{-2}(|n_2|) \sigma_{-6}(|n_2|) \left( \frac{3}{220 \pi ^2 n_2^2 y^5}-\frac{1}{36  y^3}+O\left(\frac{1}{y^2}\right) \right),
\end{equation}
and
\begin{equation}\label{eq:9.5convergence}
    \alpha_{n_2,-n_2} = O(|n_2|^{-6}), \quad |n_2| \to \infty.
\end{equation}

\section{Appendix A}\label{sec:Appendix A}

In this section, we will provide the explicit solutions which do not appear in \cite{CGPWW2021}.

\subsection{$\alpha = \beta =3/2$ and $\lambda = 2$}\label{sec:3/2,2}

\subsubsection{$n_1=0$ and $n_2=0$.}\label{section6askjdhgkajhsdk}

Any solution of \eqref{eq:fn1n2} for $n_1=n_2=0$ is equal to 
\[
\hat{f}_{0,0}(y) =\frac{16 \zeta (2)^2-9 \zeta (3)^2 y^4+72 \zeta (2) \zeta (3) y^2+48 \zeta (2)^2 \log (y)}{9 y}+c_2 y^2+\frac{c_1}{y}
\]
for some $c_1, c_2 \in \C$. 

Its asymptotic behavior for $y \to 0$ can be described by 
\[
\hat{f}_{0,0}(y) =\frac{16 \zeta (2)^2 \log (y)}{3 y}+y^{-1}\left(\frac{16 \zeta (2)^2}{9 }+c_1\right)+O(1),
\]
and the leading term of the asymptotic behavior does not depend on $c_1$ and $c_2$. Taking $c_2=0$, the $O(y^2)$-term in the asymptotic expansion of $\hat{f}_{0,0}(y)$ vanishes. However, we refrain from choosing $c_1$ until Section \ref{Section6.askojdoakjsdo}.

\subsubsection{$n_1n_2=0$ but not both zero}
For $\hat f^P_{n,0}$ as in \eqref{eq:6.1fdef}, we have 
\begin{align*}
 \nu_0(n,y)&=\sgn(n) \left( \frac{ \zeta (2)}{3 \pi  n y}-\frac{\zeta (3) y}{2\pi  n} \right),\\ 
 \nu_1(n,y)&=-\frac{\zeta (3)}{2\pi ^2 n^2}.
\end{align*}
We note that the asymptotic behavior of $\hat{f}^P_{0,n}(y)$ for $y\to 0$ is as follows:
\[
\frac{16 \zeta (2) \sigma _2(|n|) \log (y)}{3 n^2 y}+y^{-1}\frac{4 \sigma _2(|n|) \left(3 \zeta (3)+4 \gamma  \pi ^2 \zeta (2) n^2+4 \pi ^2 \zeta (2) n^2 \log (\pi  |n|)\right)}{3 \pi ^2 n^4}+O\left(y\right),
\]
where $\gamma$ is the Euler-Mascheroni constant.
We recall from \eqref{asymptotics_of_homogeneous_solutionnasokdijfok} that 
\begin{equation}\label{fHn1n2smallydfger4eter}
\hat{f}^H_{0,n}(y) = \frac{\alpha_{0, n}}{4 \pi  |n|^{3/2} y}-\frac{\pi  \sqrt{|n|}\  \alpha_{0, n}}{2} y+O\left(y^2\right).
\end{equation}
There is a natural choice of $\alpha_{0,n}$ that will get rid of the $O(y^{-1})$-term in $\hat{f}^H_{0,n}(y)+\hat{f}^P_{0,n}(y)$ (however, $O(y^{-1})$-term is not the leading term):
\[
\alpha_{0,n} =-\frac{16  \sigma_2(|n|) \left(3 \zeta (3)+4 \gamma  \pi ^2 \zeta (2) n^2+4 \pi ^2 \zeta (2) n^2 \log (\pi  |n|)\right)}{3 n^{5/2} \pi}.
\]
\subsubsection{$n_1n_2\neq 0$ and $n_1+n_2\neq 0$.}

For
$\hat{f}^P_{n_1,n_2}(y) $ as in \eqref{eq:6.3fdef}, we have 
\begin{align*}
	\eta_{0,0}= \sgn(n_1 n_2) \frac{2 n_{1} n_{2} y}{3 \left(n_{1} + n_{2}\right)^{2}}, & \quad  
	\eta_{0,1}= \sgn(n_1) \frac{n_{1} \left(n_{1} + 3 n_{2}\right)}{6 \pi \left(n_{1} + n_{2}\right)^{3}}, \\
	\eta_{1,0}= \sgn(n_2)  \frac{n_{2}  \left(3 n_{1} + n_{2}\right)}{6 \pi \left(n_{1} + n_{2}\right)^{3}}, & \quad 
	\eta_{1,1}=\frac{2 n_{1} n_{2} y}{3 \left(n_{1} + n_{2}\right)^{2}}.
\end{align*}
The asymptotic expansion of $\hat{f}^P_{n_1,n_2}(y)$ as $y\to 0$ is 
\begin{align}\label{eq:6.1.3fasym}
    \frac{16 \sigma _2\left(|n_1|\right) \sigma _2\left(|n_2|\right) }{3 n_1^2 n_2^2 } \Big(  \frac{\log(y)}{y}&+ \frac{\gamma +\log(\pi)}{ y}\nonumber\\
    &+ \frac{n_1^3 \log \left(|n_1|\right)+3 n_2 n_1^2 \log \left(|n_1|\right)+3 n_2^2 n_1 \log \left(|n_2|\right)+n_2^3 \log \left(|n_2|\right)}{y \left(n_1+n_2\right){}^3}\Big).
\end{align}
We further note that 
\begin{equation}\label{fHn1n2smallydfger4etsdf3}
\hat{f}^H_{n_1,n_2}(y) = \frac{\alpha_{n_1, n_2}}{4 \pi  |n_1+n_2|^{3/2} y}-\frac{\pi  \sqrt{|n_1+n_2|}\  \alpha_{n_1, n_2}}{2} y+O\left(y^2\right).
\end{equation}
We note that we cannot get rid of the $O(y^{-1} \log(y))$-term in $\hat{f}^H_{n_1,n_2}(y) +\hat{f}^P_{n_1,n_2}(y) $ by choosing an appropriate $\alpha_{n_1,n_2}$, but we can get rid of $O(y^{-1})$ by setting
\begin{align}
\alpha_{n_1,n_2} = - & 4 \pi |n_1+n_2|^{3/2} \cdot  \frac{16 \sigma _2\left(|n_1|\right) \sigma _2\left(|n_2|\right) }{3 n_1^2 n_2^2 } \times \nonumber \\ & \times \left(  \gamma +\log(\pi)+ \frac{n_1^3 \log \left(|n_1|\right)+3 n_2 n_1^2 \log \left(|n_1|\right)+3 n_2^2 n_1 \log \left(|n_2|\right)+n_2^3 \log \left(|n_2|\right)}{\left(n_1+n_2\right){}^3}\right).
\end{align}

Choosing $n=1$ and investigating the asymptotic behavior of $\alpha_{n_1, 1-n_1}$ as $n_1 \to \infty$, we note that the sum $\sum_{n_1+n_2 = 1} \alpha_{n_1,n_2}$ does not converge. As we show in the next section, it might still be reasonable to make such a choice of $\alpha_{n_1,n_2}$.

\subsubsection{ $n_1= -n_2$.}\label{Section6.askojdoakjsdo}
For 
$\hat{f}^P_{-n_2,n_2}(y)$  defined in \eqref{eq:6.4fdef}
we have
\begin{align*}
\mu_{0,0}=\frac{ y}{6}  
+\frac{8 n_{2}^{2} \pi^{2} y^{3}}{9},  \quad 
\mu_{0,1}=  \sgn(n_2) \left( 
\frac{1 }{6 n_{2} \pi}  
+\frac{4 n_{2} \pi y^{2}}{9}\right),  \quad 
\mu_{1,1}=
-\frac{5 y}{18}  
-\frac{ 8 n_{2}^{2} \pi^{2} y^{3}}{9}. \end{align*}
Its asymptotic behavior can be described as 
\begin{equation}\label{eq:6.1.4fsaym}
 8  \frac{\sigma_2(|n_2|)^2}{|n_2|^4} \left( \frac{6 \log (\pi  |n_2|)+6 \log (y)+6 \gamma +5}{9  y}\right) +O\left(1\right).
\end{equation}
Once again, we cannot get rid of the $O(y^{-1} \log(y))$-term by choosing $\alpha_{-n_2,n_2}$ appropriately, but we can get rid of the $O(y^{-1})$-term by setting
\[
\alpha_{-n_2,n_2} := -  \frac{8\sigma_2(|n_2|)\sigma_2(|n_2|)}{9|n_2|^4} \left( 6 \log (\pi  |n_2|)+6 \gamma -5\right).
\]
We note that for such choice of $\alpha_{-n_2,n_2}$, the sum $\sum_{n_1+n_2=0, n_2 \neq 0} \alpha_{n_1,n_2}$ diverges, because the sum 
\[
\sum_{p \text{ prime}} \alpha_{-p,p} = - \frac{8}{9}\sum_{p \text{ prime}} \frac{(p^2+1)^2}{p^4} (6 \log(\pi p) + 6 \gamma - 5)
\]
diverges. However, it is still possible to formally calculate the sum of $\alpha_{-n_2,n_2}$ using the Ramanujan summation \eqref{Ramanujan2} and its derivatives  \eqref{Ramanujan3}.

We note that we can choose $c_1$ from Section \ref{section6askjdhgkajhsdk} in such a way that at least formally 
\begin{equation}\label{eq:6.1.4alphasum}
\sum_{n_2=-\infty}^\infty \alpha_{-n_2,n_2}=0.
\end{equation}

\subsection{$\alpha =3/2, \beta=5/2$ and $\lambda =6$}\label{sec:3/2,5/2,6}

\subsubsection{$n_1=0$ and $n_2=0$.}
We can find a particular solution to be 
\begin{align*}
\hat{f}_{0,0}(y) = &\frac{-6750 y^6 \zeta (3) \zeta (5)+3375 \pi ^2 y^4 \zeta (5)+100 \pi ^4 y^2 \zeta (3)+40
   \pi ^6 \log (y)+8 \pi ^6}{6750 y^2}\\ & + \frac{c_1}{y^2}+c_2y^3
\end{align*}
for some $c_1, c_2 \in \C$. Its asymptotic behavior for $y \to 0$ can be described by 
\begin{equation*}
\hat{f}_{0,0}(y) = \frac{4 \pi ^6 \log (y)}{675 y^2}
\end{equation*}
and the leading term of the asymptotic behavior doesn't depend on $c_1$ and $c_2$.

\subsubsection{$n_1=0$ and $n_2 \neq 0$} 

For $\hat{f}^P_{0,n}(y)$ as in \eqref{eq:7.1afdef}, we have 
\begin{align*}
\nu_0(n, y)&=
\frac{\zeta(2)y^{-2}}{5 n^{2} \pi^{2}}  
-\frac{ 2\zeta(3)}{n^{2} \pi^{2}},  
\ \ \ \  \nu_1(n, y)=
\sgn(n)\left(\frac{n^{2} \pi^{2} \zeta(2) - 6\zeta(3)}{3 n^{3} \pi^{3} y}  
-\frac{ \zeta(3)y}{ n \pi} \right).
 \end{align*}
The asymptotic expansion of $\hat{f}^P_{0,n}(y)$ as $y\to 0$ is 
\begin{align*}
- 8 \pi^2 \sigma_{-4}(|n|)|n|^2  \left(   y^{-2} \left( -\frac{\gamma  \zeta (2)}{5 \pi ^2 n^2}+\frac{\zeta (2)}{6 \pi ^2 n^2}-\frac{\zeta (3)}{\pi ^4 n^4}-\frac{\zeta (2) \log (\pi  |n|)}{5 \pi ^2 n^2}  \right)\right. &-\left.\frac{\zeta (2) \log (y)}{5 \pi ^2 n^2 y^2} \right)\\ &  + O(1).
\end{align*}
There is a unique choice of $\alpha_{0,n}$ that gets rid of the $y^{-2}$-term in the asymptotic expansion of $\hat{f}^P_{0,n}(y)+\hat{f}^H_{0,n}(y)$.

\subsubsection{$n_1\neq 0$ and $n_2 = 0$} 
 For $\hat{f}^P_{n,0}(y)$ as in \eqref{eq:7.1bfdef}, we have
\begin{align*}
\nu_0(n, y)&=
\sgn(n) \left(\frac{2\zeta(4)y^{-2}}{5 n \pi}  
+\frac{3\zeta(5)}{n^{3} \pi^{3}}  
-\frac{ \zeta(5)y^{2}}{2 n \pi}  \right),
 \ \ \ \  \nu_1(n, y)=
\frac{3\zeta(5)y^{-1}}{n^{4} \pi^{4}}  
+\frac{\zeta(5)y}{ n^{2} \pi^{2}}  .
\end{align*}
The asymptotic expansion of $\hat{f}^P_{n,0}(y)$ as $y\to 0$ is 
\begin{align*}
- 8 \pi \sigma_{-2}(|n|)|n| \left( 
y^{-2} \left( \frac{3 \zeta (5)}{2 \pi ^5 |n|^5}-\frac{2 \gamma  \zeta (4)}{5 \pi  |n|}-\frac{2 \zeta (4) \log (\pi  |n|)}{5 \pi  |n|}\right) -\frac{2 \zeta (4) \log (y)}{5 \pi  |n| y^2} \right) + O(1).
\end{align*}
There is a unique choice of $\alpha_{n,0}$ that gets rid of the $y^{-2}$-term in the asymptotic expansion of $\hat{f}^P_{n,0}(y)+\hat{f}^H_{n,0}(y)$.

\subsubsection{$n_1n_2\neq 0$ and $n_1+n_2\neq 0$.}

For
$\hat{f}^P_{n_1,n_2}(y)$ as in \eqref{eq:7.3fdef}, we have
\begin{align*}
	\eta_{0,0}=
	&\sgn(n_1)\left[\frac{1}{10 \pi \left(n_{1} + n_{2}\right)^{4}}(n_{1} \left(n_{1}^{2} + 4 n_{1} n_{2} + 11 n_{2}^{2}\right))\right],\\
	\eta_{0,1}=
	&\sgn(n_1n_2)\Big[\frac{y^{-1}}{10 n_{2} \pi^{2} \left(n_{1} + n_{2}\right)^{5}}(n_{1} \left(n_{1}^{3} + 5 n_{1}^{2} n_{2} + 10 n_{1} n_{2}^{2} + 10 n_{2}^{3}\right))\\
	&\ \ \ \ \ \ \ \  \ \ \ \ \ +\frac{y}{15 \left(n_{1} + n_{2}\right)^{3}}(2 n_{1} n_{2} \left(n_{1} + 5 n_{2}\right))\Big]\\
	\eta_{1,0}=
	&\frac{y^{-1}}{10 \pi^{2} \left(n_{1} + n_{2}\right)^{5}}(n_{2}^{2} \cdot \left(5 n_{1} + n_{2}\right))+\frac{y}{15 \left(n_{1} + n_{2}\right)^{3}}(2 n_{1} n_{2} \left(n_{1} + 5 n_{2}\right)),\\
	\eta_{1,1}=
	&\frac{\sgn(n_2)}{30 \pi \left(n_{1} + n_{2}\right)^{4}}(4 n_{1}^{3} + 19 n_{1}^{2} n_{2} + 44 n_{1} n_{2}^{2} + 5 n_{2}^{3}).
\end{align*}

We note that 
\begin{align}\label{eq:7.1.3fasym}
	\hat{f}^P_{n_1,n_2}(y) =& \frac{\sigma _{-4}\left(|n_2|\right) \sigma _{-2}\left(|n_1|\right)}{15 n_1 n_2^2 \left(n_1+n_2\right){}^5 y^2} \times\nonumber \\ 
	& \times 8 \left(6 n_1^5 \log (y)+30 n_2 n_1^4 \log (y)+60 n_2^2 n_1^3 \log (y)+60 n_2^3 n_1^2 \log (y) \right. 
	\\  & +30 n_2^4 n_1 \log (y)+6 n_2^5 \log (y)+6 \gamma  n_1^5+30 \gamma  n_2 n_1^4-4 n_2 n_1^4+60 \gamma  n_2^2 n_1^3
	 \nonumber 
	\\  &
	-23 n_2^2 n_1^3+60 \gamma  n_2^3 n_1^2-63 n_2^3 n_1^2+30 \gamma  n_2^4 n_1-49 n_2^4 n_1+6 \gamma  n_2^5-5 n_2^5
	 \nonumber 
	\\ &
	+6 n_1^5 \log \left(\pi | n_1|\right)
	+30 n_2 n_1^4 \log \left(\pi | n_1|\right)+60 n_2^2 n_1^3 \log \left(\pi  |n_1| \right)  \nonumber 
	\\ &  \left. +60 n_2^3 n_1^2 \log \left(\pi | n_1|\right)+30 n_2^4 n_1 \log \left(\pi | n_2|\right)+6 n_2^5 \log \left(\pi  |n_2|\right) \right) ) +O(1). \nonumber
\end{align}
We are not able to eliminate the highest term, $O(y^{-2} \log(y))$, by choosing appropriate $\alpha_{n_1, n_2}$. However, we are able to eliminate the $O(y^{-2})$-term.

\subsubsection{ $n_1= -n_2$.}

For $
\hat{f}^P_{-n_2,n_2}(y)$ as in \eqref{eq:7.4fdef}, we have
\begin{align*}
	\mu_{0,0}&= \sgn(n_2) \left( 
	\frac{1 }{10 n_{2} \pi}  
	+\frac{2 n_{2} \pi y^{2}}{75}  
	-\frac{ 64 n_{2}^{3} \pi^{3} y^{4}}{225}  \right),
	\\
	\mu_{0,1}&=
	\frac{y^{-1}}{10 n_{2}^{2} \pi^{2}}  
	-\frac{ y}{75}  
	-\frac{ 32 n_{2}^{2} \pi^{2} y^{3}}{225}  ,
	\\
	\mu_{1,1}&= \sgn(n_2) \left( 
	-\frac{9 }{100 n_{2} \pi}  
	+\frac{2 n_{2} \pi y^{2}}{225}  
	+\frac{64 n_{2}^{3} \pi^{3} y^{4}}{225}  	\right).
	\end{align*}
	The asymptotic expansion of $\hat{f}^P_{-n_2,n_2}(y) $ as $y\to 0$	is 
	\begin{equation}\label{eq:7.1.5fasym}
\hat{f}^P_{-n_2,n_2}(y) = 4  \sigma_{-2}(|n_2|) \sigma_{-4}(|n_2|) \left( 	\frac{20 \log (\pi  |n_2|)+20 \log (y)+20 \gamma +9}{25  y}+O\left(1\right) \right).
	\end{equation}
	There is a choice of $\alpha_{-n_2,n_2}$ that gets rid of the $y^{-1}$-term in the asymptotic expansion of $\hat{f}^P_{-n_2,n_2}(y) +\hat{f}^H_{-n_2,n_2}(y) $. However, manipulating $\alpha_{-n_2,n_2}$ cannot not help us get rid of the leading term.
	Thus, as in Section \ref{Section6.askojdoakjsdo}, we can choose $c_1=\alpha_{0,0}$ from Section~\ref{section6askjdhgkajhsdk} so that the contribution from the homogeneous solutions vanishes, that is, at least formally
\begin{equation}\label{eq:7.1.4alphasum}
\sum_{n_2=-\infty}^\infty \alpha_{-n_2,n_2}=0.
\end{equation}

\subsection{$\alpha = \beta =5/2$ and $\lambda =2$}\label{sec:5/2,2}

\subsubsection{$n_1=0$ and $n_2=0$.}
We note that for $n_1=n_2=0$ 
\[
\hat{f}_{0,0}(y) = \frac{-10125 y^8 \zeta (5)^2+2700 \pi ^4 y^4 \zeta (5)-4 \pi ^8}{20250 y^3}+\frac{c_1}{y}+c_2 y^2
\]
for some $c_1, c_2 \in \C$. Its asymptotic behavior for $y \to 0$ can be described by 
$
\hat{f}_{0,0}(y) = -\frac{2 \pi ^8}{10125 y^3}
$
and the leading term of the asymptotic behavior doesn't depend on $c_1$ and $c_2$.

\subsubsection{$n_1 n_2=0$ but not both zero} For $\hat{f}^P_{0,n}(y) $ as in \eqref{eq:8.1fdef}, we have
\begin{align*}
\nu_0(n, y)&=
-\frac{\zeta(5)y}{ n^{2} \pi^{2}}  ,
\\ \nu_1(n, y)&=\sgn(n)\left(
\frac{2\zeta(4)y^{-2}}{5 n \pi}  
-\frac{ \zeta(5)}{ n^{3} \pi^{3}}  -\frac{\zeta(5)y^{2}}{2 n \pi} \right).
\end{align*}
Its asymptotic behavior can be described by $y \to 0$ as 
\begin{align}
&-8\pi^2 \sigma_{-4}(|n|) |n|^2
\left( \frac{\zeta (4)}{5 \pi ^2 n^2 y^3} \right. \\  & \left. +\frac{-5 \zeta (5)+4 \gamma  \pi ^4 \zeta (4) n^4-2 \pi ^4 \zeta (4) n^4+4 \pi ^4 \zeta (4) n^4 \log (\pi  |n|)+4 \pi ^4 \zeta (4) n^4 \log (y)}{10 \pi ^4 n^4 y}+O\left(1\right)\right) .\nonumber
\end{align}
The leading term of the asymptotic expansion is $O(y^{-3})$, the second term is $O(y^{-1} \log(y))$. Manipulating homogeneous solution, we can get rid of the third term in the asymptotic expansion; that is, to get rid of $O(y^{-1})$.

\subsubsection{$n_1n_2\neq 0$ and $n_1+n_2\neq 0$.}

For $\hat{f}^P_{n_1,n_2}(y) $ as in \eqref{eq:8.3fdef}
we have  
\begin{align*}
\eta_{0,0}=
&\frac{y}{15 \left(n_{1} + n_{2}\right)^{2}}(2 n_{1} n_{2})\\
\eta_{0,1}=
&\sgn(n_2)\frac{1}{30 \pi \left(n_{1} + n_{2}\right)^{3}}(4 n_{1}^{2} + 9 n_{1} n_{2} + 3 n_{2}^{2})\\
\eta_{1,0}=
&\sgn(n_1) \frac{1}{30 \pi \left(n_{1} + n_{2}\right)^{3}}(3 n_{1}^{2} + 9 n_{1} n_{2} + 4 n_{2}^{2})\\
\eta_{1,1}=
&\sgn(n_1 n_2)\left(\frac{y^{-1}}{10 n_{1} n_{2} \pi^{2}}+\frac{y}{15 \left(n_{1} + n_{2}\right)^{2}}(2 n_{1} n_{2})\right).
\end{align*}
We note that 
\begin{align}\label{eq:8.1.3fasym}
    \hat{f}^P_{n_1,n_2}(y) = -   \frac{8 y^{-3} \sigma_{-4}(|n_1|) \sigma_{-4}(|n_2|)}{5}+o(y^{-3}).
\end{align}
We are not capable to make the $O(y^{-3})$-term vanish by manipulating $\alpha_{n_1,n_2}$.

\subsubsection{ $n_1= -n_2$.}

For $\hat{f}^P_{-n_2,n_2}(y)$ as in \eqref{eq:8.4fdef}, we have 
\begin{align*}
\eta_{0,0}&=
-\frac{7 y}{30}  
-\frac{ 8 n_{2}^{2} \pi^{2} y^{3}}{45} ,
\\
\eta_{0,1}&= \sgn(n_2) \left( 
-\frac{2 }{15 n_{2} \pi}  
-\frac{4 n_{2} \pi y^{2}}{45}  \right),
\\
\eta_{1,1}&=
\frac{y^{-1}}{10 n_{2}^{2} \pi^{2}}  
+\frac{23 y}{90}  
+\frac{8 n_{2}^{2} \pi^{2} y^{3}}{45}.
\end{align*}
The asymptotic expansion of $\hat{f}^P_{-n_2,n_2}(y) $ as $y\to 0$ is
\begin{align}\label{eq:8.1.4fasym}
- 8 \pi^2  |n_2|^2&(\sigma_{-4}(|n_2|))^2 \nonumber\\
& \times \left( \frac{1}{5 \pi ^2 |n_2|^2 y^3}+\frac{12 \log (\pi  |n_2|)+12 \log (y)+12 \gamma +1}{9   y}+O\left(1\right) \right).
\end{align}

As in Section \ref{Section6.askojdoakjsdo}, we can choose $c_1=\alpha_{0,0}$ from Section \ref{section6askjdhgkajhsdk} so that the contribution from the homogeneous solutions vanishes, that is, at least formally
\begin{equation}\label{eq:8.1.4alphasum}
\sum_{n_2=-\infty}^\infty \alpha_{-n_2,n_2}=0.
\end{equation}

\subsection{$\alpha = \beta =5/2$ and $\lambda = 12$}\label{sec:5/2,12}

\subsubsection{$n_1=0$ and $n_2=0$.}  

We note that 
\[
\hat{f}_{0,0}(y) =\frac{128 \zeta (4)^2-441 \zeta (5)^2 y^8+784 \zeta (4) \zeta (5) y^4+896 \zeta (4)^2 \log (y)}{392 y^3}+c_2 y^4+\frac{c_1}{y^3}
\]
for some $c_1, c_2 \in \C$. Its asymptotic behavior for $y \to 0$ can be described by 
$
\hat{f}_{0,0}(y) = 
$
and the leading term of the asymptotic behavior doesn't depend on $c_1$ and $c_2$.

\subsubsection{$n_1 n_2=0$ but not both zero}  
For $\hat{f}^P_{0,n}(y) $ as in \eqref{eq:8.1fdef}, we have
\begin{align*}
\nu_0(n_0, y)&=
\frac{2\zeta(4)y^{-3}}{7 n^{2} \pi^{2}}  
+\frac{15\zeta(5)y^{-1}}{n^{4} \pi^{4}}  
+\frac{3\zeta(5)y}{2 n^{2} \pi^{2}} , 
\\ \nu_1(n_0, y)&=
\sgn(n)\left(\frac{( 4\zeta(4) \pi^{4} n^{4} +150\zeta(5))y^{-2}}{10 n^{5} \pi^{5}}  
+\frac{9\zeta(5)}{n^{3} \pi^{3}}  
-\frac{ \zeta(5)y^{2}}{2 n \pi} \right).
\end{align*}

Its asymptotic expansion as $y\to 0$ is 
\begin{align*}
&-8\pi^2 \sigma_{-4}(|n|) |n|^2  \\ 
& \times \left(\frac{525 \zeta (5)-20 \gamma  \pi ^4 \zeta (4) n^4+14 \pi ^4 \zeta (4) n^4-20 \pi ^4 \zeta (4) n^4 \log (\pi  |n|)-20 \pi ^4 \zeta (4) n^4 \log (y)}{70 \pi ^6 n^6 y^3}\right.\\ & \ \ \ \ \ \ \ \ \ \ \ \ \ \ \ \ \ \ \ \ \ \ \left.+O\left(\frac{1}{y}\right) \right).
\end{align*}
The leading asymptotic expansion as $y\to 0$ is $O(y^{-3} \log(y))$. The second leading asymptotic expansion is $O(y^{-3})$ -- that one can be eliminated by manipulating the homogeneous solution.
\subsubsection{$n_1 n_2 \neq 0$, $n_1+n_2\neq 0$}
For $\hat{f}^P_{n_1,n_2}(y) $ as in \eqref{eq:8.3fdef}
we have 
\begin{align*}
	\eta_{0,0}=
	&\frac{y^{-1}}{14 \pi^{2} \left(n_{1} + n_{2}\right)^{6}}(n_{1}^{4} + 6 n_{1}^{3} n_{2} + 50 n_{1}^{2} n_{2}^{2} + 6 n_{1} n_{2}^{3} + n_{2}^{4})\\
	&+\frac{y}{105 \left(n_{1} + n_{2}\right)^{4}}(2 n_{1} n_{2} \cdot \left(7 n_{1}^{2} + 54 n_{1} n_{2} + 7 n_{2}^{2}\right))\\
	\eta_{0,1}= \sgn(n_2) 
	& \left[\frac{y^{-2}}{14 n_{2} \pi^{3} \left(n_{1} + n_{2}\right)^{7}}(n_{1}^{2} \left(n_{1}^{3} + 7 n_{1}^{2} n_{2} + 21 n_{1} n_{2}^{2} + 35 n_{2}^{3}\right)) \right. \\
	&\left. +\frac{1}{210 \pi \left(n_{1} + n_{2}\right)^{5}}(28 n_{1}^{4} + 199 n_{1}^{3} n_{2} + 775 n_{1}^{2} n_{2}^{2} + 145 n_{1} n_{2}^{3} + 21 n_{2}^{4})\right] \\
	\eta_{1,0}= \sgn(n_1) 
	& \left[\frac{y^{-2}}{14 n_{1} \pi^{3} \left(n_{1} + n_{2}\right)^{7}}(n_{2}^{2} \cdot \left(35 n_{1}^{3} + 21 n_{1}^{2} n_{2} + 7 n_{1} n_{2}^{2} + n_{2}^{3}\right)) \right. \\
	&\left. +\frac{1}{210 \pi \left(n_{1} + n_{2}\right)^{5}}(21 n_{1}^{4} + 145 n_{1}^{3} n_{2} + 775 n_{1}^{2} n_{2}^{2} + 199 n_{1} n_{2}^{3} + 28 n_{2}^{4}) \right] \\
	\eta_{1,1}= \sgn(n_1n_2) 
	& \left[ \frac{y^{-1}}{210 n_{1} n_{2} \pi^{2} \left(n_{1} + n_{2}\right)^{6}}(21 n_{1}^{6} + 166 n_{1}^{5} n_{2} + 605 n_{1}^{4} n_{2}^{2} + 1520 n_{1}^{3} n_{2}^{3} \right. \\
	& \ \ \ \ \ \ \ \ \ \ \ \ \ \ \ \ \ \ \ \ \ \ \ \ \ \ \ \ \ \ \ \ \ \ \ \ + 605 n_{1}^{2} n_{2}^{4} + 166 n_{1} n_{2}^{5} + 21 n_{2}^{6})	\\
	& +\left.\frac{2 n_{1} n_{2}\, y}{105 \left(n_{1} + n_{2}\right)^{4}} \left(7 n_{1}^{2} + 54 n_{1} n_{2} + 7 n_{2}^{2}\right)\right].
\end{align*}
We note that the asymptotic expansion is given by 
\begin{align}\label{eq:8.2.3fasym}
\hat{f}^P_{n_1,n_2}(y) =   \frac{16   \sigma_{-4}(|n_1|) \sigma_{-4}(|n_2|)}{7  
}y^{-3}\log(y)+o(y^{-3} \log(y)).
\end{align}

\subsubsection{ $n_1= -n_2$.}

For 
$\hat{f}^P_{-n_2,n_2}(y)$ as in \eqref{eq:8.4fdef}, we have 
\begin{align*}
\mu_{0,0}&=
\frac{y^{-1}}{14 n_{2}^{2} \pi^{2}}  
-\frac{13 y}{294}  
-\frac{8 n_{2}^{2} \pi^{2} y^{3}}{735}  
+\frac{256 n_{2}^{4} \pi^{4} y^{5}}{2205}  ,
\\
\mu_{0,1}&= \sgn(n_2) \left[ 
\frac{y^{-2}}{14 n_{2}^{3} \pi^{3}}  
-\frac{5}{294 n_{2} \pi}  
+\frac{4 n_{2} \pi y^{2}}{735}  
+\frac{128 n_{2}^{3} \pi^{3} y^{4}}{2205}  \right],
\\
\mu_{1,1}&=
-\frac{5 y^{-1}}{588 n_{2}^{2} \pi^{2}}  
+\frac{59 y}{1470}  
-\frac{ 8 n_{2}^{2} \pi^{2} y^{3}}{2205}  
-\frac{ 256 n_{2}^{4} \pi^{4} y^{5}}{2205}.
\end{align*}
The asymptotic expansion is 
\begin{equation}\label{eq:8.2.4fasym}
64 \pi^4  |n_2|^4 (\sigma_{-4}(|n_2|) )^2 \left( \frac{84 \log (\pi  |n_2|)+84 \log (y)+84 \gamma +5}{2352 \pi ^4 |n_2|^4 y^3}+O\left(\frac{1}{y}\right) \right).
\end{equation}

As in Section \ref{Section6.askojdoakjsdo}, we can choose $c_1=\alpha_{0,0}$ from Section \ref{section6askjdhgkajhsdk} so that the contribution from the homogeneous solutions vanishes, that is, at least formally
\begin{equation}\label{eq:8.2.4alphasum}
\sum_{n_2=-\infty}^\infty \alpha_{-n_2,n_2}=0.
\end{equation}

\subsection{$\alpha = 3/2, \beta =7/2$ and  $\lambda =12$}\label{sec:3/2,7/2,12}

\subsubsection{$n_1=0$ and $n_2=0$.}
It is not complicated to show that 
\begin{align*}
\hat{f}_{0,0}(y) =&\frac{-10418625 y^8 \zeta (3) \zeta (7)+4630500 \pi ^2 y^6 \zeta (7)+9408 \pi ^6 y^2 \zeta (3)+4480 \pi ^8 \log (y)+640 \pi ^8}{2778300 y^3}\\ & +c_2 y^4+\frac{c_1}{y^3}
\end{align*}
for some $c_1, c_2 \in \C$. Its asymptotic behavior for $y \to 0$ can be described by 
\begin{equation*}
\hat{f}_{0,0}(y) = \frac{32 \pi ^8 \log (y)}{19845 y^3}
\end{equation*}
and the leading term of the asymptotic behavior doesn't depend on $c_1$ and $c_2$.

\subsubsection{$n_1=0$ and $n_2 \neq 0$}  For $
\hat{f}^P_{0,n}(y)$ as in \eqref{eq:9.1afdef}, we have
\begin{align*}
\nu_0&=
\sgn(n)\left(-2\zeta(3) \left( \frac{ 3 y^{-1}}{n^{3} \pi^{3}}  
+\frac{ y}{2 n \pi}  \right) +4\zeta(2) \left( \frac{y^{-3}}{7 n^{3} \pi^{3}}  
+\frac{y^{-1}}{6 n \pi} \right)\right)
\\ \nu_1&=
-2\zeta(3) \left( \frac{ 3 y^{-2}}{n^{4} \pi^{4}}  
+\frac{ 2 }{n^{2} \pi^{2}} \right)+\frac{22\zeta(2) y^{-2}}{15 n^{2} \pi^{2}} .  
\end{align*}

Its asymptotic behavior is 
\begin{align*}
   & -16  \sigma_{-6}(|n|)  \times \\ 
 &\left( \frac{-315 \zeta (3)-60 \gamma  \pi ^2 \zeta (2) n^2+77 \pi ^2 \zeta (2) n^2-60 \pi ^2 \zeta (2) n^2 \log (\pi  n)-60 \pi ^2 \zeta (2) n^2 \log (y)}{105 \pi ^2 n^2 y^3}\right.\\ &  \ \ \ \ \ \ \ \left.+\, O\left(\frac{1}{y}\right) \right).
\end{align*}

\subsubsection{$n_1\neq 0$ and $n_2 = 0$}
 For $
\hat{f}^P_{n,0}(y)$ as in \eqref{eq:9.1bfdef}, we have
\begin{align*}
\nu_0&=\sgn(n)\left(
\frac{16\zeta(6)y^{-3}}{14 n_{2} \pi}  
- 15\zeta(7) \left( \frac{ 6 y^{-1}}{n_{2}^{5} \pi^{5}}  
+\frac{3 y}{5 n_{2}^{3} \pi^{3}}  
+\frac{y^{3}}{10 n_{2} \pi}  \right)\right)
\\ \nu_1&= 15\zeta(7)  \left( 
-\frac{ 6 y^{-2}}{n_{2}^{6} \pi^{6}}  
-\frac{ 18 }{5 n_{2}^{4} \pi^{4}}  
+\frac{y^{2}}{10 n_{2}^{2} \pi^{2}} \right).
\end{align*}
Its asymptotic behavior is 
\begin{align}
   & -8\sigma_{-2}(|n|)  \times \\ 
 &\left( \frac{-315 \zeta (7)-8 \gamma  \pi ^6 \zeta (6) n^6-8 \pi ^6 \zeta (6) n^6 \log (\pi  n)-8 \pi ^6 \zeta (6) n^6 \log (y)}{7 \pi ^6 n^6 y^3}+O\left(\frac{1}{y}\right) \right).\nonumber
\end{align}

\subsubsection{$n_1n_2\neq 0$ and $n_1+n_2\neq 0$.}

For $\hat{f}^P_{n_1,n_2}(y)$ as in \eqref{eq:9.3fdef}, we have 
    \begin{align*}
    \eta_{0,0}=
    & \sgn(n_1) \sgn(n_2)\Big(\frac{y^{-1}}{7 n_{2} \pi^{2} \left(n_{1} + n_{2}\right)^{6}}(n_{1} \left(n_{1}^{4} + 6 n_{1}^{3} n_{2} + 15 n_{1}^{2} n_{2}^{2} + 20 n_{1} n_{2}^{3} + 22 n_{2}^{4}\right) )\\
   &+\frac{y}{105 \left(n_{1} + n_{2}\right)^{4}}(2 n_{1} n_{2} \cdot \left(3 n_{1}^{2} + 14 n_{1} n_{2} + 35 n_{2}^{2}\right))\Big)\\
    \eta_{0,1} =
    &\sgn(n_1)\Big(\frac{y^{-2}}{7 n_{2}^{2} \pi^{3} \left(n_{1} + n_{2}\right)^{7}}(n_{1} \left(n_{1}^{5} + 7 n_{1}^{4} n_{2} + 21 n_{1}^{3} n_{2}^{2} + 35 n_{1}^{2} n_{2}^{3} + 35 n_{1} n_{2}^{4} + 21 n_{2}^{5}\right))\\
    &+\frac{1}{210 \pi \left(n_{1} + n_{2}\right)^{5}}(n_{1} \cdot \left(27 n_{1}^{3} + 143 n_{1}^{2} n_{2} + 325 n_{1} n_{2}^{2} + 497 n_{2}^{3}\right))\Big)\\
    \eta_{1,0} =
    &\sgn(n_1)\Big(\frac{y^{-2}}{7 \pi^{3} \left(n_{1} + n_{2}\right)^{7}}(n_{2}^{3} \cdot \left(7 n_{1} + n_{2}\right))\\
    &+\frac{1}{210 \pi \left(n_{1} + n_{2}\right)^{5}}(24 n_{1}^{4} + 129 n_{1}^{3} n_{2} + 293 n_{1}^{2} n_{2}^{2} + 511 n_{1} n_{2}^{3} + 35 n_{2}^{4})\Big)\\
    \eta_{1,1} =
    &\frac{y^{-1}}{210 n_{2} \pi^{2} \left(n_{1} + n_{2}\right)^{6}}(24 n_{1}^{5} + 153 n_{1}^{4} n_{2} + 422 n_{1}^{3} n_{2}^{2} + 678 n_{1}^{2} n_{2}^{3} + 822 n_{1} n_{2}^{4} + 77 n_{2}^{5})\\
    &+\frac{y}{105 \left(n_{1} + n_{2}\right)^{4}}(2 n_{1} n_{2} \cdot \left(3 n_{1}^{2} + 14 n_{1} n_{2} + 35 n_{2}^{2}\right))
    \end{align*}

We note that 
\begin{align}\label{eq:9.1.4fasym}
	\hat{f}^P_{n_1,n_2}(y)=  \frac{64   \sigma_{-2}(|n_1|) \sigma_{-6}(|n_2|) \log(y) y^{-3}}{7  }+o(y^{-3} \log(y)).
\end{align}
We cannot eliminate the leading term in the asymptotic expansion of $\hat{f}_{n_1,n_2}(y)$ by manipulating $\alpha_{n_1,n_2}$.

\subsubsection{ $n_1= -n_2$.}

For $\hat{f}^P_{-n_2,n_2}(y) $ as in \eqref{eq:9.4fdef}, we have 
\begin{align*}
	\mu_{0,0}&=
	\frac{y^{-1}}{7 n_{2}^{2} \pi^{2}}  
	+\frac{59 y}{1470}  
	-\frac{8 n_{2}^{2} \pi^{2} y^{3}}{1225}  
	+\frac{256 n_{2}^{4} \pi^{4} y^{5}}{3675}  ,
	\\
	\mu_{0,1}&=
	\sgn(n_2) \left(  \frac{y^{-2}}{7 n_{2}^{3} \pi^{3}}  
	+\frac{17}{735 n_{2} \pi}  
	+\frac{4 n_{2} \pi y^{2}}{1225}  
	+\frac{128 n_{2}^{3} \pi^{3} y^{4}}{3675}  \right),
	\\
	\mu_{1,1}&=
	-\frac{13 y^{-1}}{147 n_{2}^{2} \pi^{2}}  
	-\frac{313 y}{7350}  
	-\frac{ 8 n_{2}^{2} \pi^{2} y^{3}}{3675}  
	-\frac{ 256 n_{2}^{4} \pi^{4} y^{5}}{3675}. 
\end{align*}
The asymptotic behavior is 
\begin{align}\label{eq:9.1.5fsaym}
     - 32 &  \sigma_{-2}(|n_2|) \sigma_{-6}(|n_2|) \left( \frac{-42 \log (\pi  |n_2|)-42 \log (y)-42 \gamma -13}{147   y^3}+O\left(\frac{1}{y}\right)\right).
\end{align}

As in Section \ref{Section6.askojdoakjsdo}, we can choose $c_1=\alpha_{0,0}$ from Section \ref{section6askjdhgkajhsdk} so that the contribution from the homogeneous solutions vanishes, that is, at least formally
\begin{equation}\label{eq:9.1.4alphasum}
\sum_{n_2=-\infty}^\infty \alpha_{-n_2,n_2}=0.
\end{equation}

\section{Appendix B}\label{sec:Appendix B}

\subsection{Convolution formulas for the divisor functions}

We recall the two famous identities on divisor functions: by \cite[Theorem 291]{Hardy}, for $s>1$ and $s-a>1$,
\[
\sum_{n=-\infty, n \neq 0}^\infty \frac{\sigma_a(|n|)}{|n|^s} = 2 \zeta(s)\zeta(s-2),
\]
and 
by \cite[Theorem 305]{Hardy}, for $s>1$, $s-a>1$, $s-b>1$ and $s-a-b>1$, 
\begin{equation}\label{Ramanujan2}
\sum_{n=-\infty, n \neq 0}^\infty \frac{\sigma_a(|n|) \sigma_b(|n|)}{|n|^s} = 2 \frac{\zeta(s) \zeta(s-a) \zeta(s-b) \zeta(s-a-b)}{\zeta(2s-a-b)}.
\end{equation}
We note that the latter identity is sometimes referred to as a Ramanujan identity.
Differentiating \eqref{Ramanujan2} with respect to $s$, we obtain 
\begin{equation}\label{Ramanujan3}
\sum_{n=-\infty, n \neq 0}^\infty \frac{\sigma_a(|n|) \sigma_b(|n|) \log(|n|)}{|n|^s} = - 2 \frac{d}{ds} \left(  \frac{\zeta(s) \zeta(s-a) \zeta(s-b) \zeta(s-a-b)}{\zeta(2s-a-b)} \right).
\end{equation}

\subsection{Bessel functions and relations between them}\label{sec:bessels_and_relations_between_them}

By \cite[8.486(17)]{Gradshteyn},

\begin{equation}\label{eq:K2viaK0K1}
K_2(z) = K_0(z) + 2 \frac{K_1(z)}{z},
\end{equation}
thus 
\[
y K_1( 2 |n_1| \pi y) K_2( 2 |n_2| \pi y) = y  K_1( 2 \pi |n_1| y) \left(  K_0( 2 \pi |n_2| y) + \frac{1}{ \pi |n_2 | y} K_1( 2\pi |n_2|y)\right).
\]
We note that by \cite[8.486(12) and 8.486(13)]{Gradshteyn} and \eqref{eq:K2viaK0K1},
\begin{align*}
K_3(y) = \frac{y K_1(y)+4 K_2(y)}{y} = \frac{\left(y^2+8\right) K_1(y)+4 y K_0(y)}{y^2} = 8 \frac{K_1(y)}{y^2} + K_1(y) + 4  \frac{K_0(y)}{y},
\end{align*}
thus 
\begin{align*}
y K_1(2\pi |n_1| y) K_3( 2 \pi |n_2| y) &= K_1( 2 \pi |n_1| y) \\ & \times \left(  2 \frac{K_1(2 |n_2| y)}{\pi^2 |n_2|^2 y}  + 2  \frac{K_0(2 \pi |n_2| y)}{\pi |n_2|}   + y K_1(2 \pi |n_2| y) \right).
\end{align*}
By  \cite[8.486(17)]{Gradshteyn},
\begin{align*}
    y K_2( 2|n_1| y) K_2( 2|n_2| y) =&  y  \left(  K_0( 2\pi |n_2| y) + \frac{1}{ \pi |n_2 | y} K_1( 2\pi |n_2|y)\right) \\ & \times  \left(  2 K_0( \pi |n_2| y) + \frac{1}{\pi |n_2 | y} K_1( 2\pi|n_2|y)\right).
\end{align*}

\bibliography{fkl_bib}
\bibliographystyle{amsplain}
\end{document}